\documentclass[english,12pt]{amsart}
\usepackage{amssymb}
\usepackage{amsfonts}
\usepackage[foot]{amsaddr}
\usepackage{amscd}
\usepackage[all]{xy}
\usepackage{color}
\usepackage{enumerate}
\usepackage[utf8]{inputenc}

\usepackage[T1]{fontenc}
\usepackage{xr-hyper}
\usepackage[colorlinks=true, breaklinks=true, urlcolor= black, allcolors= black,bookmarksopen=true,linktocpage=true,plainpages=false,pdfpagelabels]{hyperref}
\usepackage{cleveref}

\CompileMatrices

\swapnumbers

\def\deg{{\rm deg}}
\def\rad{\operatorname{rad}}
\def\radop{\rad_{\mathrm{op}}}
\def\radcl{\rad_{\mathrm{cl}}}

\def\ac{{\overline{\rm ac}}}

\def\val{{\mathrm{val}}}
\def\Lval{\cL_\val}
\newcommand{\Lqe}[1][\lambda]{{\cL_{\rm {qe},#1}}}

\def\cLe{\cL'} 
\def\cLas{\cL^{\rm as}}

\def\cLarb{\tilde{\cL}} 

\def\Lbas{\cL_{\QQ,I}}  
\def\Lebas{\cLe_{\QQ,I}} 
\def\LbasLA{\Lbas'}  

\def\eqc{{\mathrm{ecc}}}  
\def\fine{{\mathrm{fine}}}  

\def\omin{{\mathrm{omin}}}  
\def\Tomin{\cT_\omin}
\def\Lomin{\cL_\omin}

\def\heqc{h$^\eqc$}

\def\GL{{\rm GL}}

\def\11{{\mathbf 1}}

\def\NN{{\mathbb N}}

\def\QQ{{\mathbb Q}}
\def\RR{{\mathbb R}}

\def\ZZ{{\mathbb Z}}

\def\cA{{\mathcal A}}
\def\cB{{\mathcal B}}

\def\cF{{\mathcal F}}

\def\cL{{\mathcal L}}
\def\cM{{\mathcal M}}

\def\cO{{\mathcal O}}

\def\cT{{\mathcal T}}
\def\cU{{\mathcal U}}

\def\llp{\mathopen{(\!(}}

\def\rrp{\mathopen{)\!)}}

\mathchardef\alphag="7C0B \mathchardef\betag="7C0C
\mathchardef\gammag="7C0D \mathchardef\deltag="7C0E
\mathchardef\varepsilong="7C22 \mathchardef\varphig="7C27
\mathchardef\psig="7C20 \mathchardef\zetag="7C10
\mathchardef\epsilong="7C0F \mathchardef\rhog="7C1A
\mathchardef\taug="7C1C \mathchardef\upsilong="7C1D
\mathchardef\iotag="7C13 \mathchardef\thetag="7C12
\mathchardef\pig="7C19 \mathchardef\sigmag="7C1B
\mathchardef\etag="7C11 \mathchardef\omegag="7C21
\mathchardef\kappag="7C14 \mathchardef\lambdag="7C15
\mathchardef\mug="7C16 \mathchardef\xig="7C18
\mathchardef\chig="7C1F \mathchardef\nug="7C17
\mathchardef\varthetag="7C23 \mathchardef\varpig="7C24
\mathchardef\varrhog="7C25 \mathchardef\varsigmag="7C26
\mathchardef\Omegag="7C0A \mathchardef\Thetag="7C02
\mathchardef\Sigmag="7C06 \mathchardef\Deltag="7C01
\mathchardef\Phig="7C08 \mathchardef\Gammag="7C00
\mathchardef\Psig="7C09 \mathchardef\Lambdag="7C03
\mathchardef\Xig="7C04 \mathchardef\Pig="7C05
\mathchardef\Upsilong="7C07

\newtheorem{thm}[subsubsection]{Theorem}
\newtheorem{lem}[subsubsection]{Lemma}
\newtheorem{cor}[subsubsection]{Corollary}
\newtheorem{prop}[subsubsection]{Proposition}

\newtheorem{claim}[subsubsection]{Claim}
\newtheorem{addendum}
{Addendum}
\Crefname{prop}{Proposition}{Propositions}
\Crefname{thm}{Theorem}{Theorems}
\Crefname{lem}{Lemma}{Lemmas}
\Crefname{cor}{Corollary}{Corollaries}
\Crefname{conj}{Conjecture}{Conjectures}

\theoremstyle{definition}
\newtheorem{defn}[subsubsection]{Definition}
\newtheorem{conv}[subsubsection]{Convention}
\newtheorem{notn}[subsubsection]{Notation}
\newtheorem{example}[subsubsection]{Example}

\newtheorem{def-prop}[subsubsection]{Proposition-Definition}
\newtheorem{def-theorem}[subsubsection]{Theorem-Definition}
\newtheorem{def-lem}[subsubsection]{Lemma-Definition}
\newtheorem{assumption}[subsubsection]{Assumption}
\Crefname{defn}{Definition}{Definitions}

\theoremstyle{remark}
\newtheorem{remark}[subsubsection]{Remark}

\newtheorem{question}[subsubsection]{Question}

\theoremstyle{plain}

\newenvironment{condition}{\begin{itemize}\refstepcounter{equation}\item[(\theequation)]}{\end{itemize}}

\newcounter{dummy}

\numberwithin{equation}{subsection}

\def\boxit#1#2{\setbox1=\hbox{\kern#1{#2}\kern#1}%
\dimen1=\ht1 \advance\dimen1 by #1 \dimen2=\dp1 \advance\dimen2 by
#1
\setbox1=\hbox{\vrule height\dimen1 depth\dimen2\box1\vrule}%
\setbox1=\vbox{\hrule\box1\hrule}%
\advance\dimen1 by .4pt \ht1=\dimen1 \advance\dimen2 by .4pt
\dp1=\dimen2 \box1\relax}

\renewcommand{\theequation}{\thesubsection.\arabic{equation}}

\mathchardef\alphag="7C0B \mathchardef\betag="7C0C
\mathchardef\gammag="7C0D \mathchardef\deltag="7C0E
\mathchardef\varepsilong="7C22 \mathchardef\varphig="7C27
\mathchardef\psig="7C20 \mathchardef\zetag="7C10
\mathchardef\epsilong="7C0F \mathchardef\rhog="7C1A
\mathchardef\taug="7C1C \mathchardef\upsilong="7C1D
\mathchardef\iotag="7C13 \mathchardef\thetag="7C12
\mathchardef\pig="7C19 \mathchardef\sigmag="7C1B
\mathchardef\etag="7C11 \mathchardef\omegag="7C21
\mathchardef\kappag="7C14 \mathchardef\lambdag="7C15
\mathchardef\mug="7C16 \mathchardef\xig="7C18
\mathchardef\chig="7C1F \mathchardef\nug="7C17
\mathchardef\varthetag="7C23 \mathchardef\varpig="7C24
\mathchardef\varrhog="7C25 \mathchardef\varsigmag="7C26
\mathchardef\Omegag="7C0A \mathchardef\Thetag="7C02
\mathchardef\Sigmag="7C06 \mathchardef\Deltag="7C01
\mathchardef\Phig="7C08 \mathchardef\Gammag="7C00
\mathchardef\Psig="7C09 \mathchardef\Lambdag="7C03
\mathchardef\Xig="7C04 \mathchardef\Pig="7C05
\mathchardef\Upsilong="7C07

\newcommand{\cha}{\operatorname{char}}
\newcommand{\RV}{\mathrm{RV}}
\newcommand{\VF}{\mathrm{VF}}
\newcommand{\VG}{\mathrm{VG}}
\newcommand{\RF}{\mathrm{RF}}
\newcommand{\rv}{\operatorname{rv}}
\newcommand{\res}{\operatorname{res}}
\newcommand{\Th}{\operatorname{Th}}
\newcommand{\valring}{\mathcal{O}_K}

\newcommand{\eq}{^{\mathrm{eq}}}

\newcommand{\ltz}{\mathrel{<\joinrel\llap{\raisebox{-1ex}{$\scriptstyle{0}\mkern8mu$}}}}

\DeclareMathOperator*{\sgn}{{sgn}}

\newcommand{\grad}{\operatorname{grad}}

\def\Jac{\operatorname{Jac}}

\definecolor{immi}{rgb}{0,.6,.1}

\newbox\removebox
\newcommand\remove[1]{%
\setbox\removebox=\ifmmode\hbox{$#1$}\else\hbox{#1}\fi%
\leavevmode
\rlap{\textcolor{blue}{\vrule height0.8ex depth-0.6ex width\wd\removebox}}%
\box\removebox
}
\long\def\bigremove#1{%
\par\setbox\removebox=\vbox{#1}%
\vbox{%
\vbox to0pt{\hbox{\tikz\draw[color=blue,thick] (0,0) -- (\wd\removebox,-\ht\removebox)  (\wd\removebox,0) -- (0,-\ht\removebox);}}
\box\removebox
}
}

\newcommand\tp{\mathrm{tp}}

\newcommand\acl{\mathrm{acl}}
\newcommand\dcl{\mathrm{dcl}}
\newcommand{\restr}[2]{{\left.#1\right|_{#2}}}

\definecolor{orange}{rgb}{1,0.5,0}

\newcommand{\private}[1]{\leavevmode{\scriptsize\color{blue}\marginpar{{\scriptsize Private comment}}#1\par}}
\renewcommand{\private}[1]{}

\thanks{The authors would like to thank Y.~Halevi, J.~Koenigsmann, F.~Loeser, A.~Macintyre, W.~Singhof and F.~Vermeulen for interesting discussions on the topics of the paper, and L.~van den Dries for suggesting the name of Hensel minimality.
The author R.\,C. was partially supported by the European Research Council under the European Community's Seventh Framework Programme (FP7/2007-2013) with ERC Grant Agreement nr.\ 615722 MOTMELSUM, by KU Leuven IF C14/17/083, and thanks the Labex CEMPI  (ANR-11-LABX-0007-01). The author I.\,H. was partially supported by the \emph{SFB~878: Groups, Geometry and Actions}, by the research training group
\emph{GRK 2240: Algebro-Geometric Methods in Algebra, Arithmetic and Topology}, and by the individual research grant No.~426488848, all three funded by the Deutsche Forschungsgemeinschaft (DFG, German Research Foundation). Part of the work has been done while I.\,H. was affiliated to the University of Leeds
The author S.\,R. was partially supported by GeoMod AAPG2019 (ANR-DFG), Geometric and Combinatorial Configurations in Model Theory.}

\title{Hensel minimality I}

\author[Cluckers]{Raf Cluckers}
\address{Raf Cluckers, Univ.~Lille, CNRS, UMR 8524 - Laboratoire Paul Painlev\'e, F-59000 Lille, France, and,
KU Leuven, Department of Mathematics, B-3001 Leu\-ven, Bel\-gium}
\email{Raf.Cluckers@univ-lille.fr}
\urladdr{http://rcluckers.perso.math.cnrs.fr/}

\author[Halupczok]{Immanuel Halupczok}
\address{Immanuel Halupczok, Lehrstuhl f\"ur Algebra und Zahlentheorie, Mathematisches Institut, Universit\"atsstr. 1, 40225 D\"usseldorf, Germany}
\email{math@karimmi.de}
\urladdr{http://www.immi.karimmi.de/en/}

\author[Rideau]{Silvain Rideau-Kikuchi}
\address{Silvain Rideau-Kikuchi, CNRS, Université Paris Diderot, Sorbonne Université, Institut de Mathématiques de Jussieu - Paris Rive Gauche (UMR 7586), 8 place Aurélie Nemours, Boite Courrier 7012, 75205 PARIS Cedex 13, France}
\email{silvain.rideau@imj-prg.fr}
\urladdr{https://webusers.imj-prg.fr/~silvain.rideau/en/}



\begin{document}

\begin{abstract}
We present a framework for tame geometry on Henselian valued fields which we call Hensel minimality.
In the spirit of o-minimality, which is key to real geometry and several diophantine applications,
we develop geometric results and applications for Hensel minimal structures that were previously known only under stronger, less axiomatic assumptions.
We show existence of t-stratifications in Hensel minimal structures and
Taylor approximation results which are key to non-archimedean versions of Pila--Wilkie point counting, Yomdin's parameterization results and to motivic integration.
In this first paper we work in equi-characteristic zero; in the sequel paper, we develop the mixed characteristic case and a diophantine application.
\end{abstract}

\maketitle

\newpage

\setcounter{tocdepth}{1}
\tableofcontents

\section{Introduction}
\label{sec:intro}

\subsection{Background}
Our goal is to understand what really lies behind the tameness of definable sets in valued fields and to describe this axiomatically,
thereby providing a basement for research in several directions of non-archimedean geometry. In the present paper, we introduce the framework called `Hensel minimality' which is very natural and readily verifiable from Ax--Kochen/Ershov type results.
Moreover, it has strong consequences that were previously not expected to follow from such basic axioms.

\medskip

Let us begin by describing some history of definable sets in valued fields and their uses for $p$-adic and motivic integration, geometry, and for diophantine applications.
From the 60s on, various results were obtained providing rather precise control over semi-algebraic definable sets in
the fields $\QQ_p$ of $p$-adic numbers and in fields $k\llp t\rrp$ of formal Laurent series with coefficients in a field $k$ of characteristic $0$ by Ax, Kochen, Ershov, Cohen, Macintyre and others \cite{AK1,AK2,AK3,Ersov,Cohen,Ersov,Mac}.
Based on those, Denef \cite{Denef} solved a rationality question of Serre \cite{Serre}
which established a strong connection between number theory and the geometry of $p$-adic definable sets. This in turn became a
motivation for the further development of the model theory of valued fields \cite{Pas,Basarab}, of motivic integration \cite{DLinvent,DL} (an integration theory on arc spaces of varieties),
and
for the study of singularities and of Berkovich spaces by Hrushovski and Loeser among others \cite{Baty-Moreau,DLBarc,HruL}. The model theoretic approach also brought Fubini and Fourier to motivic integrals \cite{CLoes, CLexp, HK} which was key to applications to the Langlands program \cite{CHL,Casselman-Cely-Hales}. Recently, semi-algebraic and subanalytic $p$-adic geometry has led to point counting results analogous to Pila--Wilkie point counting and Yomdin's parameterization results \cite{CCL-PW,CFL, PW, YY}.  Also applications to Lipschitz stratifications were obtained for semi-algebraic and subanalytic sets from non-archimedean geometry \cite{Halup}. See \cite{Bel} for a more extensive panorama of results. These historical results are all based on concrete settings like semi-algebraic and subanalytic geometry; we open them up
by providing a natural axiomatic and more general framework.

On the field of real numbers, a successful framework exists for tame geometry since work by van den Dries, Pillay, and other  \cite{DriesTarski, vdD, KnightPillSt, PillaySteinhornI}.
Indeed, tameness and geometry of real definable sets is very well captured by the notion of o-minimality, an axiomatic condition about definable sets.
It became a central tool in real algebraic geometry and its generalizations, on the one hand because of its beautiful naturality and on the other hand because of its vast consequences on the geometry of definable sets and its strong diophantine applications, e.g. to the André-Oort Conjecture \cite{Pila}.

\medskip

Soon after the first successes of o-minimality, the open quest for analogous frameworks in valued fields emerged as a central question in model theory. Several notions have been presented so far, each of which with some shortcomings or lack of generality. In analogy to o-minimality there are P-minimality \cite{Haskell}, C-minimality \cite{Macpherson,HM}; more motivated by applications there are V-minimality \cite{HK} and  b-minimality \cite{CLb}; motivated by classification theory there is dp-minimality \cite{Onshuus,Jahnke-Simon}.
In the present paper, our new notion of Hensel minimality overcomes these shortcomings and provides a strong and natural framework for non-archimedean tame geometry.
Hensel minimality is easier to axiomatize and more natural than b-minimality and V-minimality, it applies more broadly than P-minimality, C-minimality, V-minimality, or even dp-minimality, and it has stronger consequences than b-minimality, P-minimality and C-minimality.

In fact, we introduce two key notions of Hensel minimality, namely $1$-h-minimality and $\omega$-h-minimality. Most of our results are developed under the weaker assumption of $1$-h-minimality; only some of the resplendency results from Section~\ref{sec:respl} need $\omega$-h-minimality.
As auxiliary notions we define some intermediary variants of Hensel minimality, namely $\ell$-h-minimality for integers $\ell>1$.
We use the term ``Hensel minimality'' to talk about any of those variants, in an implicit way.
Most classically studied examples turn out to be
$\omega$-h-minimal, but for some interesting (somewhat less classical) examples, we were only able to prove $1$-h-minimality (see Section~\ref{sec:examples}).

In this paper, we work in equi-characteristic zero. The mixed characteristic case is derived from this in the sequel \cite{CHRV}.
Indeed, the mixed characteristic version of Hensel minimality comes almost for free (including many of its consequences) from the equi-characteristic zero situation, by using coarsened valuations. In Section~\ref{sec:examples} we briefly define the mixed characteristic variant of Hensel minimality, for the sake of our examples.

\medskip

We now sketch ideas, results and applications.
In a way, o-minimality assumes that every unary definable set is controlled by a finite set, namely its set of boundary points; in this sense, Hensel minimality is very similar to o-minimality. To make this more precise,
the definition of an o-minimal field $R$ can be phrased as the following condition on definable subsets $X$ of $R$: There exists a finite subset $C$ of $R$ such that, for any $x\in R$, whether $x$ lies in $X$ depends only on the tuple $(\sgn(x-c))_{c\in C}$, where $\sgn$ stands for the sign, which can be $-$, $0$, or $+$. The definition of Hensel minimality is similar,
where the sign function is replaced by a suitable function adapted to the valuation, namely the leading term map $\rv$. However,
whereas in the o-minimal world, $C$ is automatically definable over the same parameters as $X$, in the valued field setting, we need to impose precise conditions on the parameters over which $C$ can be defined.
This is also where the differences between $1$-h-minimality and $\omega$-h-minimality arise.

A large part of this paper consist in proving our main geometric results in Hensel minimal structures, which are similar to those in o-minimal structures, in particular
cell decomposition, dimension theory,
the ``Jacobian Property'' (which plays a key role in constructing motivic integration, and which can be considered as an analogue of the Monotonicity Theorem from the o-minimal context where $\sgn$ is replaced by $\rv$),
as well as higher order and higher dimension versions of the Jacobian Property, which state that definable functions have good approximations by their Taylor polynomials. Based on those properties, various recent results in the model theory of Henselian valued fields readily generalize to arbitrary Hensel minimal valued fields, like those on Lipschitz continuity \cite{CCLLip} and t-stratifications (which were introduced in \cite{Halup} and studied further in \cite{Gar.cones,Gar.powbd,HalYin}). As an extra upshot, Hensel minimality intrinsically has resplendency properties in the spirit of resplendent quantifier elimination, i.e., it is preserved by different kinds of expansions of the structure. In particular, it should be considered as a notion of tameness ``relative to the leading term structure $\RV$''.

As first applications, we show existence of t-stratifications in arbitrary Hensel minimal structures, and we
use our results in valued fields to deduce
a uniform Taylor approximation result in power-bounded o-minimal real closed fields, which
strengthens results from \cite{Halup,HalYin};
here, the connection between valued fields and real closed fields is made using the notion of $T$-convexity from \cite{DL.Tcon1,Dri.Tcon2}.
Furthermore, our geometric results lay the ground for further generalizations of motivic integration from \cite{CLoes, CLexp,HK} and its use
in e.g.~\cite{Bilu-Th,ChaL,CCL,CGH5,CLexp,ForeyDens,HKPoiss,HMartin,Kien:rational,Nicaise-semi,Yin-int,Yimu-t-conv}. Concretely, we show how Hensel minimality relates to the axiomatic frameworks under which motivic integration is developed in \cite{CLbounded} and \cite{HK}.
Our results also lay the ground for $C^r$ parameterizations and for bounds on the number of rational points as in \cite{CCL-PW,CFL}, analogous to results by Yomdin \cite{YY2} and by Pila--Wilkie \cite{PW}.
Motivic integration under Hensel minimality is under development in \cite{CLMellin}, and a first diophantine application to counting of rational points of bounded height
is given in \cite{CHRV}.

\subsection{The notion of Hensel minimality}
\label{sec:main:notion}

We start by giving the
definition of Hensel minimality;
see Section \ref{sec:hmin:def} for more details.

Let $K$ be a non-trivially valued field of equi-characteristic $0$,
considered as an $\cL$-structure for some language $\cL$ containing the language $\Lval = \{+,\cdot,\cO_K\}$ of valued fields (where $\cO_K$ is a predicate for the valuation ring).
We use multiplicative notation for valuations. We denote the value group by $\Gamma^\times_K$ and the valuation map by
$$
|\cdot|\colon K\to \Gamma_K :=  \Gamma_K^\times \cup\{0\};
$$
see Section \ref{sec:not:val} for more detailed definitions.

The analogue to the sign map from the o-minimal context will be the ``leading term map'' $\rv\colon K \to \RV$, defined as follows:

\begin{defn}[Leading term structure $\RV_\lambda$]\label{defn:RVI}
Let $\lambda \le 1$ be an element of $\Gamma_K^\times$ and set $I := \{x \in K \mid |x| < \lambda\}$. We define $\RV_\lambda^\times$ to be the quotient of multiplicative groups $K^\times/(1 + I)$, and we let
\[
\rv_\lambda\colon K \to \RV_\lambda := \RV_\lambda^\times \cup \{0\}
\]
be the map extending the projection map $K^\times\to \RV_\lambda^\times$ by sending $0$ to $0$. We abbreviate $\RV_1$ and $\rv_1$ by $\RV$ and $\rv$, respectively. If several valued fields are around, we may also write $\RV_K$ and $\RV_{K,\lambda}$.
\end{defn}

We can now make precise in which sense a set $X\subset K$ can be controlled by a finite set $C \subset K$.

\begin{defn}[Prepared sets; see Definition \ref{defn:lambda-prepared}]
\label{defn:lambda-prepared:intro}
Let $\lambda \le 1$ be in $\Gamma_K^\times$, let $C$ be a finite non-empty subset of $K$ and let $X \subset K$ be an arbitrary subset.
We say that $C$ \emph{$\lambda$-prepares} $X$ if whether some $x\in K$ lies in $X$ depends only on the tuple $(\rv_\lambda(x-c))_{c\in C}$. In other words, if $x, x' \in K$ satisfy
\begin{equation}\label{eq:xyIball}
\rv_\lambda(x-c) = \rv_\lambda(x'-c)  \mbox{ for each } c\in C,
\end{equation}
then one either has $x\in X$ and $x'\in X$, or, one has $x\not\in X$ and $x'\not\in X$.
\end{defn}

A key ingredient of Hensel minimality is that any $A$-definable set $X \subset K$ (for $A \subset K$) can be $1$-prepared by a finite $A$-definable set $C$. This however is not yet strong enough a condition: We need some precise control of parameters from $\RV$ and $\RV_\lambda$. This is where we obtain different variants of Hensel minimality, the difference consisting only in some number $\ell$ of allowed parameters:

\begin{defn}[$\ell$-h-minimality; see Definition \ref{defn:hmin}]\label{defn:hmin:intro}
Let $\ell\geq 0$ be either an integer or $\omega$, and let $\cT$ be a theory
of valued fields of equi-characteristic $0$, in a language $\cL$ containing the language $\Lval$ of valued fields. We say that $\cT$ is \emph{$\ell$-h-minimal} if every model $K \models \cT$ has the following property:
\begin{condition}\label{cond:l-h-min-intro}
For every $\lambda \le 1$ in $\Gamma_K^\times$, for every set $A\subset K$
and for every set $A' \subset \RV_\lambda$ of cardinality $\#A' \le \ell$, every $(A \cup \RV \cup A')$-definable set $X \subset K$ can be $\lambda$-prepared by a finite $A$-definable set $C \subset K$.
\end{condition}
\end{defn}

In this definition, the parameters from $\RV$ and $\RV_\lambda$ are imaginary elements (see Section~\ref{sec:not:mod}).

The notion of $\ell$-h-minimality seems most natural in the case $\ell = \omega$, and indeed, this case is especially interesting in view of the resplendency results of Section~\ref{sec:respl}. However,
the case $\ell = 1$ plays an important role, too: most of the results in this paper follow from $1$-h-minimality, and we obtain a rather different equivalent characterization of $1$-h-minimality in terms of definable functions (see Theorem~\ref{thm:tame2vf}). This shows that
$1$-h-minimality is of deeply geometric nature,
even though its definition might appear somewhat
contrived.
Similar characterizations for integers $\ell \ge 2$ are developed in \cite{Verm:h-min}, involving properties of definable functions on $\ell$-dimensional sets.

\medskip

Many structures on valued fields whose model theory is known to behave well are
Hensel minimal. In Section~\ref{sec:examples}, we provide examples of Hensel minimal structures on valued fields, of mixed and of equi-characteristic $0$. Some are new and less expected, others incorporate more classical structures.
We show that the following classical structures are Hensel minimal:

\begin{enumerate}
 \item\label{exam:semi-alg}
    The theory of Henselian valued fields of equi-characteristic $0$ in the pure valued field language is $\omega$-h-minimal.
 \item\label{exam:analytic}
   The theory stays $\omega$-h-minimal if we expand the valued field by analytic functions (forming an analytic structure as in \cite{CLip} or \cite{CLips}): For example, we can consider the theory of a field $K = k((t))$ with $\cha k = 0$, in a language containing all analytic functions $\cO_K^n \to K$ for all $n$, where analytic here means that the function is given by evaluating a $t$-adically converging power series in $x\in \cO_K^n$ (see Section 3.4 of \cite{CLips}).
\item\label{it:tcon}
   If $\Tomin$ is the theory of a power-bounded o-minimal expansion of a real closed field, then the theory of $\Tomin$-convex valued fields in the sense of \cite{DL.Tcon1} is $1$-h-minimal.
   Recall that a $\Tomin$-convex valued fields is obtained by taking a (sufficiently big) model $K$ of $\Tomin$ and turning it into a valued field using the convex closure of an elementary substructure $K_0 \prec K$ as valuation ring; see Section~\ref{sec:Tcon} for details.
\item\label{it:resp} Any expansion of the structure by predicates on Cartesian powers of $\RV$ preserves $0$-, $1$- and $\omega$-h-minimality (Theorem~\ref{thm:resp:h}).
 \item\label{it:coars-rob} The notion of $\omega$-h-minimality is preserved under coarsening of the valuation (see Corollary~\ref{cor:coarse}). In \cite[Theorem~2.2.7]{CHRV}, it is proved that $1$-h-min is also preserved  under coarsening of the valuation.
 \end{enumerate}

By combining (\ref{it:coars-rob}) and (\ref{exam:analytic}), one recovers the examples of coarsened analytic structures from \cite{Rid}, and resplendently so by (\ref{it:resp}).

On the other hand, there are also some well-behaved structures on valued fields which are certainly not Hensel minimal: Since the very definition of Hensel minimality implies that every definable subset of $K$ is either finite or contains a ball, one cannot add a section of the residue field or of the value group to the language, and neither can one add an automorphism of $K$.

As a converse to (\ref{exam:semi-alg}), we show that a valued field which is Hensel minimal is automatically Henselian, see Theorem \ref{thm:hens}. In Proposition \ref{prop:vmin} we compare Hensel minimality to V-minimality from \cite{HK}.

\subsection{Description of the main results}
\label{sec:main:res}

We now list some of the consequences of $1$-h-minimality (often in simpler forms than the versions in the main part of the paper).
The first result is the ``Jacobian Property'' which plays a crucial role e.g.\ in motivic integration, both in the Cluckers--Loeser version \cite{CLoes,CLexp} and in the Hrushovski--Kazhdan version \cite{HK}.

\begin{thm}[Jacobian Property; see Corollary~\ref{cor:JP}]\label{thm:JP:intro}
Let $f\colon K \to K$ be definable (with parameters). Then there exists a finite set $C \subset K$ such that the following holds: for every fiber $B$ of the map sending $x \in K$ to the tuple $(\rv(x - c))_{c \in C}$, there exists a $\xi_B \in \RV$ such that
\begin{equation}\label{eq:JP:intro}
\rv\left(\frac{f(x_1) - f(x_2)}{x_1 - x_2}\right) = \xi_B
\end{equation}
for all $x_1, x_2 \in B$, $x_1 \ne x_2$.
\end{thm}

\begin{remark}
This formulation of the Jacobian Property becomes exactly the o-minimal Monotonicity Theorem if one replaces all occurrences of $\rv$ by $\sgn$ and if one adds the condition that $f$ should be continuous on the fibers $B$ (which, in the valued field setting, follows automatically).
\end{remark}

Corollary~\ref{cor:JP} also includes the following strengthenings of Theorem~\ref{thm:JP:intro}:
\begin{itemize}
 \item
The theorem still holds with $\rv$ replaced by $\rv_\lambda$ for any $\lambda \le 1$ in $\Gamma_K^\times$ (still only assuming $1$-h-minimality), and one can moreover choose a single finite set $C$ that works for all such $\lambda$.
 \item If $f$ is definable with parameters from $A \cup \RV$ with $A \subset K$, then $C$ can taken to be $A$-definable. (Corollary~\ref{cor:JP} only speaks about $\emptyset$-definable sets $f$; Remark~\ref{rem:acl} explains how to deduce this more general version.)
\end{itemize}

Another point of view of the Jacobian Property is that on each fiber $B$ (using notation from the Theorem), $f$ has a good approximation by its first order Taylor series.
One of the deepest results of this paper is a similar result for higher order Taylor approximations:

\begin{thm}[Taylor approximations; see Theorem~\ref{thm:high-ord}]\label{thm:high-ord:intro}
Let $f\colon K \to K$ be a definable function and let $r \in \NN$ be given. Then there exists a finite set $C$ such that for every
fiber $B$ of the map sending $x \in K$ to the tuple $(\rv(x - c))_{c \in C}$,
$f$ is $(r+1)$-fold differentiable on $B$ and we have
\begin{equation}\label{eq:t-higher:intro}
 \left|f(x) -  \sum_{i = 0}^r \frac{f^{(i)}(x_0)}{i!}(x - x_0)^i \right|  \le  |f^{(r+1)}(x_0)\cdot (x-x_0)^{r+1}|
\end{equation}
for every $x_0, x \in B$.
\end{thm}
As in Theorem~\ref{thm:JP:intro}, if $f$ is $(A \cup \RV)$-definable for $A \subset K$, then $C$ can be taken $A$-definable (again using Remark~\ref{rem:acl}).

Using that $\Tomin$-convex valued fields are $1$-h-minimal (see (\ref{it:tcon}) above, on p.~\pageref{it:tcon}), this Taylor approximation result implies
a uniform Taylor approximation result in power-bounded o-minimal real closed fields; see Corollary~\ref{cor:arch}.

\medskip

As in the o-minimal context, the preparation of unary sets and the Jacobian Property lend themselves well (by logical compactness) to obtain results about higher dimensional objects (namely, definable subsets of $K^n$ and definable functions on $K^n$). Indeed, this can be pursued, \emph{mutatis mutandis}, in the style of how cell decomposition and dimension theory are built up from b-minimality in \cite{CLb}. In particular, we obtain results about
\begin{itemize}
 \item
almost everywhere differentiability (Subsection~\ref{sec:cont}), \item
cell decomposition in two variants (Subsections~\ref{sec:cd} and ~\ref{sec:cd:classical}),
\item
dimension theory (Subsection~\ref{sec:dim}), and
\item higher dimensional versions
of the Jacobian Property and of Taylor approximations (Subsections~\ref{sec:sjp} and \ref{sec:taylor-box}).
\end{itemize}
For cell decomposition, Subsection~\ref{sec:cd} provides a new approach specific to Hensel minimality: Usually, the notion of cells in valued fields is a lot more technical than in o-minimal structures, partly due to the lack of (certain) Skolem functions. Item~(\ref{it:resp}) above (on p.~\pageref{it:resp}) allows us to add the missing Skolem functions to the language without destroying Hensel minimality, and there are also some tools enabling us to get back to the original language afterwards (see Subsection~\ref{sec:alg:skol}).
This allows us, in Subsection~\ref{sec:cd}, to work with a notion of cell which is much less technical than most previous ones in valued fields.

\medskip

One of our original goals was to deduce the Jacobian Property starting from abstract conditions on unary sets only.
Previous proofs of the Jacobian Property either use piecewise analyticity arguments, even in the semi-algebraic case (like in \cite{CLip}), or they are restricted to rather specific situations (like \cite{HK} for V-minimal valued fields and \cite{Yin.tcon} for power-bounded $T$-convex valued fields). Our new proof is more general and works for arbitrary Henselian valued fields which are $1$-h-minimal, and goes without analyticity arguments.

\medskip

The remainder of this paper is organized as follows:
In Section~\ref{sec:first-properties}, after fixing notation and terminology,  we develop many basic tools that are useful for proofs in Hensel minimal theories,
and we obtain first elementary results like a key ingredient to dimension theory (Lemma \ref{lem:fin-inf}) and
a weak version of the Jacobian Property (Lemma~\ref{lem:gammaLin}). We also show that those two results are essentially equivalent to
$1$-h-minimality (Theorem~\ref{thm:tame2vf}).

The deepest results of this paper are contained in Section~\ref{sec:derJac} about definable functions from $K$ to $K$,
namely almost everywhere differentiability and Taylor approximation.

Section~\ref{sec:respl} is devoted to understanding in which sense Hensel minimality is a notion ``relative to $\RV$'' and to proving that
various modifications of the language preserve Hensel minimality.

The geometric results in $K^n$ (like cell decomposition, dimension theory) are collected in
Section~\ref{sec:compactn}. That section also contains our main application, of t-stratifications, our higher dimensional Taylor approximation results, and Cluckers--Loeser style motivic integration.

Finally, in Section~\ref{sec:examples}, we show that many previously studied structures are Hensel minimal, and we give some new examples as well, namely coarsened valued fields as variants of classical analytic structures. Since some of those examples are of mixed characteristic, at the beginning of Section~\ref{sec:examples}, we briefly define the mixed characteristic variant of Hensel minimality which is treated in full in the sequel \cite{CHRV} to this paper.

As an application we show in Section~\ref{sec:Tcon} how Hensel minimality results yield corresponding results in power-bounded real closed fields. At the end, we compare our notion to V-minimality from \cite{HK}.

\section{Hensel minimality}\label{sec:hmin}\label{sec:first-properties}

\subsection{Notation and terminology for valued fields and balls}\label{sec:not:val}

In the entire paper, $K$ will denote a valued field of characteristic zero, with valuation ring $\cO_K \subset K$ and maximal ideal $\cM_K \subset \cO_K$. In this paper, we only consider non-trivially valued fields, i.e., $\cO_K \ne K$. Moreover, apart from Section \ref{sec:examples}, $K$ will always be of equi-characteristic $0$, meaning that both $K$ and the residue field $\cO_K/\cM_K$ have characteristic $0$.

Note that we allow Krull-valuations (and thus valuations of arbitrary rank), that is,  we allow $\cO_K$ to be an arbitrary (proper, by the non-triviality) subring of $K$ such that for every element $x\in K^\times$, at least one of $x$ or $x^{-1}$ belongs to $\cO_K$.
The value group is then defined to be the quotient $\Gamma_K^\times := K^\times / \cO_K^\times$ of multiplicative groups.

We denote the valuation map by $|\cdot|\colon K^\times \to \Gamma_K^\times$ and use multiplicative notation for the value group. We write $\Gamma_K$ for the disjoint union $\Gamma_K^\times \cup \{0\}$, we extend the valuation map to $|\cdot|\colon K \to \Gamma_K$ by setting $|0| := 0$,
and we define the order on $\Gamma_K$ in such a way that $\cO_K = \{x \in K \mid |x| \le 1\}$ and $|x|<|y|$ whenever  $x/y\in \cM_K$  for  $x$ and nonzero $y$ in $K$.

For $x = (x_1, \dots, x_n) \in K^n$, we set $|x| := \max_i |x_i|$.

We use the (generalized) leading term structures $\RV_\lambda$ (for $\lambda \le 1$ in $\Gamma_K^\times$) that have already been introduced in Definition~\ref{defn:RVI},
and we also denote the natural map $\RV_\lambda \to \Gamma_K$ by $|\cdot|$.
Note that $\RV_\lambda$ is a semi-group for multiplication.

\begin{remark} Write $\RV^\times$ for $\RV\setminus \{0\}$.
Recall that one has a natural short exact sequence of multiplicative groups $(\cO_K/\cM_K)^\times \to \RV^\times \to \Gamma^\times_K$. (So
$\RV$ combines information from the \underline{r}esidue field and \underline{v}alue group.)
\end{remark}

\begin{example}
In the case $K = k((t))$, the above short exact sequence naturally splits, giving an isomorphism $\RV^\times \to  (\cO_K/\cM_K)^\times \times \Gamma^\times_K$,
which, for $a = \sum_{i=N}^\infty a_i t^i \in K^\times$ with $a_N \ne 0$, sends $\rv(a)$ to $(a_N, N)$.
\end{example}

\begin{remark}
For $a, a' \in K$, one has $\rv_\lambda(a) = \rv_\lambda(a')$ if and only if either $a = a' = 0$, or $|a - a'| < \lambda\cdot |a|$.
\end{remark}

We consider several kinds of balls:

\begin{defn}[Balls]\label{defn:balls}
\begin{enumerate}
 \item
We call a subset $B \subset K^n$ a \emph{ball} if $B$ is infinite, $B \ne K^n$, and for all $x,x'\in B$ and all $y\in K^n$ with $|x-y|\leq |x-x'|$ one has $y\in B$.
\item
By an \emph{open ball}, we mean a set of the form
$$
B = B_{<\gamma}(a) := \{x\in K^n\mid | x - a | < \gamma \}
$$
for some $a\in K^n$ and some $\gamma \in \Gamma_K^\times$.
\item
By a \emph{closed ball}, we mean a set of the form
$$
B = B_{\le\gamma}(a) := \{x\in K^n\mid | x - a | \leq  \gamma \}
$$
for some $a\in K^n$ and some $\gamma \in \Gamma_K^\times$.
\item
For $B$ as in (2) or (3), we call $\gamma$ the \emph{radius} of the open (resp.~closed) ball and denote it by $\radop(B)$ (resp.~$\radcl(B)$).
\end{enumerate}
\end{defn}

We call the valuation on $K$ discrete if there is a uniformizing element $\varpi$ in $\cO_K$, namely satisfying $\varpi\cO_K=\cM_K$.

\begin{remark}\label{rem:rad-op}
When $\Gamma_K^\times$ is discrete, a ball $B$ can be open and closed at the same time, but with $\radcl(B) = |\varpi|\cdot\radop(B) < \radop(B)$.
\end{remark}

Note also that for any $\lambda \le 1$ in $\Gamma_K^\times$ and for any $\xi \in \RV_\lambda \setminus \{0\}$, the preimage $\rv_{\lambda}^{-1}(\xi)$ is an open ball satisfying
$\radop(\rv_{\lambda}^{-1}(\xi)) = |\xi| \cdot \lambda$.

\begin{defn}[$\lambda$-next balls]\label{defn:next}
Fix $\lambda \le 1$ in $\Gamma^\times_K$.
\begin{enumerate}
 \item
We say that a ball $B$ is $\lambda$-next to an element $c \in K$
if
$$
B= \{x\in K\mid \rv_\lambda(x-c) = \xi \}
$$
for some (nonzero) element $\xi$ of $\RV_\lambda$.
\item
We say that a ball $B$ is $\lambda$-next to a finite non-empty set $C\subset K$
if $B$ equals $\bigcap_{c\in C} B_c$
with $B_c$ a ball $\lambda$-next to $c$ for each $c\in C$.
\end{enumerate}
\end{defn}

\begin{remark}
Using that the intersection of the finitely many balls $B_c$ is either empty or equal to one of the $B_c$, one deduces that
every ball $B$ which is $\lambda$-next to $C$ is in particular $\lambda$-next to one element $c \in C$.
\end{remark}

\begin{remark}
Given a finite non-empty set $C \subset K$,
the fibers of the map $x \mapsto (\rv_\lambda(x -  c))_{c \in C}$ are exactly the singletons consisting of one element of $C$ and the balls $\lambda$-next to $C$. In particular, the balls $\lambda$-next to $C$ form a partition of $K \setminus C$, and a subset $X$ of $K$
is $\lambda$-prepared by $C$ (as in Definition~\ref{defn:lambda-prepared}) if and only if every ball $B$ which is $\lambda$-next to $C$ is either contained in $X$ or disjoint from $X$.
\end{remark}

\begin{example}\label{ex:1prep}
The balls $1$-next to an element $c \in K$ are exactly the maximal balls in $K$ not containing $c$. From this, one deduces that
a ball $1$-next to a finite set $C$ is exactly a maximal ball disjoint from $C$. This means that a set $X \subset K$ is $1$-prepared by $C$ if and only if every ball disjoint from $C$ is either contained in $X$ or disjoint from $X$. Note again how closely ``every definable set can be $1$-prepared'' resembles o-minimality (where ``balls disjoint from $C$'' becomes ``intervals disjoint from $C$'').
\end{example}

Given a subset $A$ of a Cartesian product $B\times C$ and an element $b\in B$, we write $A_b$ for the fiber $\{c\in C\mid (b,c)\in A\}$. Also, for a function $g$ on $A$, we write $g(b,\cdot)$ for the function on $A_b$ sending $c$ to $g(b,c)$.

\subsection{Model theoretic notations and conventions}
\label{sec:not:mod}

As already stated, in the entire paper, $K$ is a non-trivially valued field of characteristic zero, and outside of Section \ref{sec:examples}, $K$ is moreover of equi-characteristic $0$.
In the entire paper, we fix a language $\cL$ containing the language $\Lval$ of valued fields, and we consider the valued field $K$ as an $\cL$-structure. More precisely,
as ``language of valued fields'', it suffices for us to take $\Lval = \{+,\cdot, \cO_K\}$, where $\valring$ is a predicate for the valuation ring; in any case, we only care about which sets are definable (and not how they are definable). For that reason, we will often specify languages only up to interdefinability.

If not specified otherwise, ``definable'' always refer to the fixed language $\cL$. As usual, ``$\cL$-definable'' means definable (in $\cL$) without additional parameters, ``$A$-definable'' means $\cL(A)$-definable, and
``definable'' means $\cL(K)$-definable (i.e., with arbitrary parameters).

Sometimes, we will also consider $K$ as a structure in other languages (e.g.\ $\cL'$); in that case, we may specify the language as an index, writing e.g.\ $\Th_{\cL'}(K)$ for the theory of $K$ considered as an $\cL'$-structure.

Almost in the entire paper (more precisely, everywhere except in parts of Section~\ref{sec:fixI}), $\cL$ will be a one-sorted language. Nevertheless, we often
work with imaginary sorts of $K$, i.e., quotients $K^n/\mathord{\sim}$ where $\sim$ is a $\emptyset$-definable equivalence relation. In particular, we consider imaginary definable sets and imaginary elements.
As usual, this can be made formal either by working in the $\cL\eq$-structure $K\eq$ (see e.g.\ \cite[Proposition~8.4.5]{TentZiegler}), or, equivalently, by ``simulating'' imaginary objects in $\cL$, namely as follows:
\begin{itemize}
 \item By a ``definable subset $X$ of $K^n/\mathord{\sim}$'', we really mean its preimage in $K^n$, i.e., a definable set $\tilde X \subset K^n$ which is a union of $\sim$-equivalence classes.
 \item If, in a formula $\phi(x, \dots)$, the variable $x$ runs over an imaginary sort $K^n/\mathord{\sim}$, this means that we really have an $n$-tuple $\tilde{x}$ of variables
 (running over $K^n$) and that the truth value of $\phi(\tilde{x}, \dots)$
 only depends on the equivalence class of $\tilde x$ modulo $\sim$.
 \item If $A$ is a set of potentially imaginary elements, then by ``$A$-definable'', we mean definable in the expansion of $K$ by predicates for the equivalence classes (in $K^n$ for some $n$) corresponding to the imaginary elements $a \in A$; the corresponding extension of the language $\cL$ by predicate symbols is denoted by $\cL(A)$.
\item A (potentially imaginary) element $b$ is said to be in the definable closure of a set $A$ (of potentially imaginary elements) if the equivalence class in $K^n$ corresponding to $b$
  is $\cL(A)$-definable. If $X$ is any imaginary sort (or even more generally an arbitrary set of imaginary elements), we write $\dcl_X(A)$ for the set of elements from $X$ that are
  in the definable closure of $A$. Being in the algebraic closure, with notation $\acl_X(A)$, is defined accordingly.
\end{itemize}

The value group $\Gamma_K$ is of course an imaginary sort. In general, $\RV_\lambda$ (for $\lambda \le 1$ in $\Gamma^\times_K$) is, by itself, not an imaginary sort, since the equivalence relation used to define it may not be $\emptyset$-definable. However, the disjoint union
\[
\RV_\bullet := \bigcup_{\lambda \le 1} \RV_\lambda
\]
is an imaginary sort and $\RV_\lambda$ is a definable subset of
$\RV_\bullet$; in particular, it makes sense to use elements from $\RV_\lambda$ as parameters.

\subsection{Hensel minimality}\label{sec:hmin:def}

We first restate the definitions of prepared sets and of Hensel minimality from the introduction. We use $\RV_\lambda$ and $\rv_\lambda$ from Definition \ref{defn:RVI}.

\begin{defn}[Prepared sets]\label{defn:lambda-prepared}
Let $K$ be a non-trivially valued field of equi-characteristic $0$.
Let $\lambda \le 1$ be in $\Gamma_K^\times$, let $C$ be a finite non-empty subset of $K$ and let $X \subset K$ be an arbitrary subset.
We say that $C$ \emph{$\lambda$-prepares} $X$ if
there exists a set $\Xi \subset\RV_\lambda^{\#C}$
such that $X$ is the preimage of $\Xi$ under the map $K\to \RV_\lambda^{\#C}, x \mapsto (\rv_\lambda(x-c))_{c\in C}$.
\end{defn}

The condition that $C$ is non-empty is essentially irrelevant, since one could always add $0$ to $C$, but it will sometimes avoid pathologies.

\begin{example}
A subset of $K$ is $\lambda$-prepared by the set $C = \{0\}$ if and only if it is the preimage under $\rv_\lambda$ of a subset of $\RV_\lambda$.
\end{example}

\begin{defn}[$\ell$-h-minimality]\label{defn:hmin}
Let $\ell\geq 0$ be either an integer or $\omega$, and let $\cT$ be a (possibly non-complete) theory
of valued fields of equi-characteristic $0$, in a language $\cL$ containing the language $\Lval$ of valued fields. We say that $\cT$ is \emph{$\ell$-h-minimal} if every model $K \models \cT$ has the following property:
\begin{condition}\label{cond:l-h-min}
For every $\lambda \le 1$ in $\Gamma_K^\times$, for every set $A\subset K$
and for every set $A' \subset \RV_\lambda$ of cardinality $\#A' \le \ell$, every $(A \cup \RV \cup A')$-definable set $X \subset K$ can be $\lambda$-prepared by a finite $A$-definable set $C \subset K$.
\end{condition}
\end{defn}

Note that a theory $\cT$ is $\ell$-h-mininimal if and only if each of its completions is. Conversely, compactness arguments show that most of our results hold uniformly in a theory $\cT$ if and only they hold in every completion (see \Cref{prop:uniform} for example). For this reason we will often work with complete theories in this paper.

\begin{remark}
In the case $\ell = 0$, $\lambda$ plays no role and the definition simplifies to: Any $(A \cup \RV)$-definable set (for $A \subset K$) can be $1$-prepared by a finite $A$-definable set.
\end{remark}

\begin{remark}
In the case $\ell = \omega$, we can more generally allow $X$ to use parameters $\xi_i \in \RV_{\lambda_i}$ for different $\lambda_i$ (using that we have surjections $\RV_\lambda \to \RV_{\lambda'}$ for $\lambda' > \lambda$);
in that case, $C$ is required to $\lambda$-prepare $X$ for $\lambda := \min_i \lambda_i$.
\end{remark}

\begin{remark}
The assumption that $C$ is definable using only the parameters from $K$ will enable us
to simultaneously prepare families of sets parametrized by $\RV$ (and $\RV_\lambda$).
This plays a central role to make Hensel minimality independent of the structure on $\RV$.
\end{remark}

\subsection{Basic model theoretic properties of Hensel minimality}

Let us now verify that the notion of Hensel minimality (in all its variants) has some basic
properties one would expect from any good model theoretic notion.

First of all, note that it is preserved under expansions of the structure by constants; more precisely:

\begin{lem}[Adding constants]\label{lem:addconst}
Let $\ell\geq 0$ be an integer or $\omega$,
suppose that $\Th(K)$ is $\ell$-h-minimal, and let $A$ be any subset of $K \cup \RV\eq$. Then $\Th_{\cL(A)}(K)$ is also $\ell$-h-minimal.
\end{lem}

Here, by $\RV\eq$, we mean imaginary sorts of the form $\RV^n/\mathord{\sim}$, for some $n$ and some $\emptyset$-definable equivalence relations $\sim$. In particular, the lemma allows us to add constants from $\Gamma_K$.
Note however that adding parameters from other sorts than $K$ and $\RV\eq$ may destroy Hensel minimality.

The lemma should be clear in the case $A \subset K \cup \RV$; we mainly give the following proof to show that parameters from $\RV\eq$ are not a problem either.

\begin{proof}[Proof of Lemma~\ref{lem:addconst}]
We verify Definition~\ref{defn:hmin}:
Let $K'\models \Th_{\cL(A)}(K)$ and let $X\subset K'$ be $\cL(A\cup A'\cup\RV_{K'}\cup A'')$-definable, for some
$A' \subset K'$ and some $A'' \subset \RV_{K',\lambda}$ satisfying $\#A'' \le \ell$, with $\lambda \le 1$ in $\Gamma_{K'}$.
Choose $\tilde A \subset \RV_{K'}$ such that every element of $A \cap \RV_{K'}\eq$ is $\cL(\tilde A)$-definable.
Then $X$ is $\cL((A \cap K) \cup \tilde{A}\cup A'\cup\RV_{K'}\cup A'')$-definable,
so by $\ell$-h-minimality of $\Th_{\cL}(K')$, $X$ can be $\lambda$-prepared by a finite $\cL((A \cap K) \cup A')$-definable set $C$.
In particular, \(C\) is $\cL(A\cup A')$-definable, as desired.
\end{proof}

With this lemma in mind, many results in this paper are formulated for $\emptyset$-definable sets; those results then automatically also hold for $A$-definable sets, when $A \subset K \cup \RV\eq$,
and using a compactness argument (given in Remark~\ref{rem:compact}), one then often obtains family versions of the results.

As so often in model theory, it is sufficient to verify Hensel minimality in sufficiently saturated models. To see this, we first prove that preparation is a first order property in the following sense:

\begin{lem}[Preparation is first order]\label{lem:prep:def}
Let $X_q$ and $C_q$ be $\emptyset$-definable families of subsets of $K$, where $q$ runs over a $\emptyset$-definable subset $Q$ of an arbitrary possibly imaginary sort. Suppose that $C_q$ is finite for every $q \in Q$. Then the
set of pairs $(q,\lambda) \in Q \times \Gamma^\times_K$ with
$\lambda \le 1$ such that $C_q$ $\lambda$-prepares $X_q$ is $\emptyset$-definable.
\end{lem}

\begin{proof}
If $\phi(x, q)$ defines $X_q$ and $\psi(z, q)$ defines $C_q$, then
the above set of pairs $(q,\lambda)$ is defined by the following $\cL$-formula:
\[
\forall x,x'\colon \big(\underbrace{\forall z\colon (\psi(z, q) \to \rv_\lambda(x - z) = \rv_\lambda(x' - z))}_{\text{i.e., }(\rv_\lambda(x-c))_{c\in C_q} = (\rv_\lambda(x'-c))_{c\in C_q}} \to (\phi(x, q) \leftrightarrow \phi(x', q))\big)
\]
\end{proof}

\begin{lem}[Saturated models suffice]\label{lem:emin-sat} Let $\ell\geq 0$ be either an integer or $\omega$ and suppose that $K$ is $\aleph_0$-saturated. Then the theory $\Th(K)$ of $K$ is $\ell$-h-minimal if and only if $K$ satisfies Condition~(\ref{cond:l-h-min}) from Definition~\ref{defn:hmin}.
\end{lem}

\begin{proof}
We need to show that if $K$ satisfies (\ref{cond:l-h-min}), then so does any
other model $K'$ of $\Th(K)$.

Suppose for contradiction that $K'$ is a model not satisfying (\ref{cond:l-h-min}), i.e., there exist a $\lambda \le 1$ in $\Gamma^\times_{K'}$, tuples $a \in (K')^n$, $\zeta \in \RV_{K'}^{n'}$, $\xi \in \RV_{K',\lambda}^{n''}$ with $n, n'$ arbitrary and $n'' \le \ell$, and an $(a, \zeta, \xi)$-definable set
$X = \phi(K', a, \zeta, \xi) \subset K'$ such that no finite non-empty $a$-definable set $C \subset K$ $\lambda$-prepares $X$.

For fixed $\phi$, the non-existence of $C$ can be expressed by an infinite conjunction of $\cL$-formulas
in $(\lambda, a, \zeta, \xi)$. Indeed,
for every formula $\psi(z, a)$ that could potentially define $C$ and for every integer $k \ge 1$, there is (by Lemma~\ref{lem:prep:def}) an $\cL$-formula $\chi_\psi(\lambda, a, \zeta, \xi)$ expressing
``$\psi(K', a)$ has cardinality $k$ and $\psi(K', a)$ does not $\lambda$-prepare $\phi(K', a, \zeta, \xi)$.''

The fact that this partial type $\{\chi_\psi \mid \psi\text{ as above}\}$ is realized in $K'$ implies that it is also realized in $K$, so that Condition~(\ref{cond:l-h-min}) fails in $K$.
\end{proof}

Whether a structure is o-minimal can also be characterized via its $1$-types. For $0$-h-minimality, we have a similar characterization:

\begin{lem}[$0$-h-minimality in terms of types]\label{lem:type-0-h-min}
Suppose that $K$ is $\aleph_0$-saturated.
The theory $\Th(K)$ is $0$-h-minimal
if and only if, for every parameter set $A \subset K$ and every ball $B \subset K \setminus \acl_K(A)$, any two elements of $B$ have the same type over $A \cup \RV$.
\end{lem}

\begin{remark}
One could also formulate similar conditions for $\ell$-h-minimality, but it would be more technical.
\end{remark}

\begin{proof}[Proof of Lemma~\ref{lem:type-0-h-min}]
``$\Rightarrow$'': Suppose for contradiction that $B$ contains two elements $x, x'$ with $\tp(x/A \cup \RV) \ne \tp(x'/A \cup \RV)$. This means that there exists an $(A \cup \RV)$-definable set $X$ containing $x$ but not $x'$. By $0$-h-minimality, there exists a finite $A$-definable set $C$ $1$-preparing $X$. In particular, $C \subset \acl_K(A)$ and hence $C \cap B = \emptyset$.
However, by Example~\ref{ex:1prep}, this implies that $B$ is either contained in $X$ or disjoint from $X$, contradicting the properties of $x$ and $x'$.

``$\Leftarrow$'': Let $X$ be $(A \cup \RV)$-definable, for some parameter set $A \subset K$ which we may assume to be finite, and suppose that no finite $A$-definable $C \subset K$ $1$-prepares $X$.
This means (using Example~\ref{ex:1prep} again) that for every finite $A$-definable $C \subset K$, there exists an open ball $B \subset K$ that is disjoint from $C$ and such that neither $B \subset X$ nor $B \subset K \setminus X$.
Taking all those conditions on $B$ together (for all finite $A$-definable $C$), we obtain a (finitely satisfiable) type in an imaginary variable running over the open balls in $K$. A realization of this type is a ball $B$ that is disjoint from $\acl_K(A)$ on the one hand, but on the other hand contains elements $x$ and $x'$ that satisfy $x \in X$ and $x' \notin X$ and hence have different types over $A \cup \RV$.
\end{proof}

Yet another way to characterize o-minimality is: Every unary definable set is already quantifier free definable in the language $\{<\}$.
Using Lemma~\ref{lem:prep eq qf}, one obtains a similar kind of characterization of $0$-h-minimality and of $\omega$-h-minimality.

\subsection{Basic properties under weaker assumptions}
\label{sec:weak-h}

Recall (from the beginning of Section \ref{sec:not:val})  that $K$ is a non-trivially valued field of equi-characteristic zero, considered as a structure in a language $\cL \supset \Lval = \{+, \cdot,\cO_K\}$.
In this Section \ref{sec:weak-h}
we assume $\Th(K)$ to be ``Hensel minimal without control of parameters'', namely:
\begin{assumption}\label{ass:no-ctrl}
For every $K' \equiv K$, every definable (with any parameters) subset of $K'$ can be $1$-prepared by a finite set $C$. (We do not impose definability conditions on $C$.)
\end{assumption}

Note that this assumption is preserved under adding constants to $\cL$ (even from arbitrary imaginary sorts), so below, every occurrence of ``$\emptyset$-definable'' can also be replaced by ``$A$-definable''.

\begin{lem}[$\exists^\infty$-elimination]\label{lem:finite}
Under Assumption~\ref{ass:no-ctrl},
every infinite definable set $X \subset K$ contains an (open) ball. In particular, if $\{X_q\mid q\in Q\}$ is a $\emptyset$-definable family of subsets $X_q$ of $K$, for some $\emptyset$-definable set $Q$ in an arbitrary possibly imaginary sort, then
the set $Q' := \{q\in Q \mid X_q\text{ is finite}\}$ is a $\emptyset$-definable set, and there exists a uniform bound $N \in \NN$ on the cardinality of $X_q$ for all $q \in Q'$.
\end{lem}

\begin{proof}
If $X$ is infinite, then it is not contained in the finite set $C$ preparing it, which implies that it contains a ball.
The definability of $Q'$ then follows since the condition that a set contains an open ball can be expressed by a formula. The existence of the bound $N$ then follows by compactness: If no bound would exist, then in a sufficiently saturated model, we would find a $q \in Q'$ with $\#X_q > N$ for every $N$.
\end{proof}

\begin{lem}[Finite sets are $\RV$-parametrized]\label{lem:average}
Under Assumption~\ref{ass:no-ctrl}, let $C_q\subset K$ be a $\emptyset$-definable family of finite sets, where $q$ runs over some $\emptyset$-definable set $Q$ in an arbitrary possibly imaginary sort.
Then there exists a $\emptyset$-definable family of injective maps $f_q\colon C_q \to \RV^k$ (for some $k$).
\end{lem}

\begin{proof}[Proof of Lemma~\ref{lem:average}]
Using Lemma~\ref{lem:finite}, we can bound the cardinality $\#C_q$ and then assume that it is constant.
We do an induction over $\#C_q$.

If $C_q$ is always a singleton or empty, we can define $f_q$ to always be constant. Otherwise,
the lemma is obtained by repeatedly taking averages of the elements of $C_q$ and subtracting. More precisely, setting $a_q := \frac1{\#C_q}\sum_{x \in C_q} x$, we get that the map $\hat f_q\colon C_q\to \RV, x \mapsto \rv(x - a_q)$ is not constant on $C_q$.
Therefore, each fiber $\hat f^{-1}_q(\xi)$ of $\hat f$ (for $\xi \in \hat f_q(C_q)$) has cardinality less than $C_q$, so
by induction, we obtain a definable family of
injective maps $g_{q,\xi}\colon \hat f^{-1}_q(\xi) \to \RV^k$. Now set $f_q(x) := (\hat f_q(x), g_{q,\hat f_q(x)}(x))$.
\end{proof}

The family of balls $1$-next to some finite set $C \subset K$ can be parameterized by $\RV$-variables; more precisely (and more generally), we have the following:

\begin{lem}[$\lambda$-next balls as fibers]\label{lem:InextFam}
Under Assumption~\ref{ass:no-ctrl}, let $\lambda \le 1$ be an element of $\Gamma^\times_K$,
let $A$ be any set of possibly imaginary parameters containing $\lambda$, and
let $C\subset K$ be a finite non-empty $A$-definable set.
Then there exists an $A$-definable map $f\colon K \to \RV^k \times \RV_\lambda$ (for some $k$) such that each nonempty fiber of $f$ is either a singleton contained in $C$ or contained
in a single ball $\lambda$-next to $C$. In the case $\lambda = 1$,
we may even obtain that each fiber of $f$ which is not a singleton is equal to a ball $1$-next to $C$.
\end{lem}

\begin{proof}
Given $x \in K$, let $\mu(x) := \min\{|x-c| \mid c \in C\}$ be the minimal distance to elements of $C$, let $C(x) := \{c \in C \mid |x-c| = \mu(x)\}$ be the set $c \in C$ realizing that distance,
and let $a(x) := \frac1{\#C(x)}\sum_{c \in C(x)} c$ be the average of those elements.
Note that the map $a\colon K \to K$ has finite image.
Using Lemma~\ref{lem:average}, we find an injective map $\alpha$ from the image of $a$ to $\RV^k$.
If $\lambda < 1$, we define
\[
f(x) := (\alpha(a(x)), \rv_\lambda(x - a(x))).
\]
In the case $\lambda = 1$,
to obtain the more precise statement, we define
\[
f(x) := \begin{cases}
(\alpha(a(x)), \rv(x - a(x))) & \text{if } |x - a(x)| \ge \mu(x)\\
(\alpha(a(x)), \rv(0))        & \text{if } |x - a(x)| < \mu(x).
        \end{cases}
\]
It is now just a computation to check that the lemma holds. For the first part, suppose that $f(x_1) = f(x_2) = (\zeta, \xi) \in \RV^k \times \RV_\lambda$ for some $x_1, x_2 \in K$; our aim is to show that either $x_1 = x_2 \in C$, or that they both lie in the same ball $\lambda$-next to $C$.
Set $a_0 := a(x_1) = a(x_2)$.

If $\xi = 0$, then either $x_i = a_0$ or $|x_i - a_0| < \mu(x_i)$. The latter is only possible if $a_0 \notin C$ and it implies that $x_i$ and $a_0$ lie in the same ball $1$-next to $C$. This shows: Either $x_1 = x_2 \in C$, or $\lambda = 1$ and both lie in the same ball $1$-next to $C$.

If $\xi \ne 0$, then the fact that $|x_i - a(x)| \le \mu(x_i)$ shows that $\rv_\lambda(x_i - a(x))$ determines $\rv_\lambda(x_i - c)$ for each $c \in C$, so we are done with the first part of the lemma.

For the second part, pick $x_1, x_2$ in the same ball $\lambda$-next to $C$. Then one obtains $\mu(x_1) = \mu(x_2) =: \mu_0$, $C(x_1) = C(x_2) =: C_0$ and $a(x_1) = a(x_2) =: a_0$. If $|x_i - a_0| = \mu$, then $\rv(x_i - a_0)$ is determined by $\rv(x_i - c)$ for any $c \in C_0$, so $f(x_1) = f(x_2)$. Otherwise, i.e., if $|x_i - a_0| < \mu$,
then $f(x_i) = (\alpha(a_0), 0)$ for $i = 1,2$ by definition.
\end{proof}

\begin{remark}\label{rem:InextFam}
In Lemma~\ref{lem:InextFam}, we can also find a map $f$ with codomain $\RV_\lambda^{k+1}$ instead of $\RV^k \times \RV_\lambda$; indeed, in Lemma~\ref{lem:average} and its proof, $\RV$ can be replaced by $\RV_\lambda$ everywhere.
\end{remark}

\subsection{Preparing families}
\label{sec:Prepfam}

In this subsection, we show that $\ell$-h-minimality does not only imply that we can
prepare definable subsets of $K$ (by finite sets $C$), but also
various other kinds of definable objects that ``live in $K \times \RV^k \times \RV_\lambda^\ell$''.

\begin{defn}[Preparing families]\label{defn:uniform}
Let $C \subset K$ be a non-empty set and let $\lambda \le 1$ be an element of $\Gamma_K^\times$.
\begin{enumerate}
 \item\label{it:uniform}
Suppose that $W$ is, up to permutation of coordinates, a subset of $K \times Q$, where $Q$ is a Cartesian product of some (possibly imaginary) sorts.
We say that $C$ \emph{$\lambda$-prepares} $W$
if for every ball $B \subset K$ $\lambda$-next to $C$, the fiber
$W_x \subset Q$ does not depend on $x$, when $x$ runs over $B$.
(The terminology will only be applied when exactly one of the coordinates of $W$ runs over $K$; the fibers $W_x$ always are over that coordinate.)

We synonymously also say that $C$ \emph{$\lambda$-prepares} the family $(W_\xi)_{\xi \in Q}$ of subsets of $K$. Indeed, the above condition is equivalent to $C$ $\lambda$-preparing $W_\xi \subset K$ for each $\xi \in Q$.
\item If $f$ is a definable function whose graph $W$
lives in a Cartesian product as in (\ref{it:uniform}), we say that $C$ \emph{$\lambda$-prepares} $f$ if it $\lambda$-prepares $W$.
\item
We say that $C$ \emph{prepares} a set
\[W \subset K \times Q \times \bigcup_{\lambda \in \Gamma_K^\times, \lambda \le 1}\RV_\lambda^{k},\]
\emph{uniformly in $\lambda$}, where $Q$ is a product of sorts and $k \ge 0$ is an integer,
if it $\lambda$-prepares
$W \cap (K \times Q \times \RV_\lambda^{k})$ for each $\lambda \le 1$ in $\Gamma^\times_K$. As in (\ref{it:uniform}), we allow the coordinates of $W$ to be in a different order, but making sure that this is unambiguious.
\end{enumerate}
\end{defn}

\begin{prop}[Preparing families]\label{prop:uniform}
Assume that $\Th(K)$ is $\ell$-h-minimal and that $A$ is a subset of $K$.
For any integer $k\ge 0$
and any $(A \cup \RV)$-definable set
$$
W\subset K\times \RV^k\times \bigcup_{\lambda \le 1}\RV_\lambda^{\ell},
$$
there exists a finite non-empty $A$-definable set $C$
which prepares $W$ uniformly in $\lambda$.
\end{prop}

\begin{proof}
For each $\lambda$ and each $\xi\in \RV^k\times \RV_{\lambda}^\ell$, let $C_\xi$ be a finite $A$-definable set $\lambda$-preparing $W_\xi$. By a usual compactness argument (see Remark~\ref{rem:compact} below), we may suppose that $C := \bigcup_\xi C_\xi$ is finite and  $A$-definable. It prepares each $W_\xi$ for each $\xi$ and hence also $W$.
\end{proof}

\begin{remark}\label{rem:compact}
In the above proof, we used a compactness argument which we will be using (in variants) many times in this paper. We give some details once: First recall that by Lemma~\ref{lem:prep:def}, ``preparing'' is a definable condition. In particular, the set
\[
\Xi_{\xi} :=\{\xi' \in  \RV^k\times  \bigcup_{\lambda \le 1}\RV_{\lambda}^\ell \mid
C_\xi\text{ $\lambda$-prepares $W_{\xi'}$}
\}
\]
is $A$-definable.
Since $\xi \in \Xi_{\xi}$, the union of all $\Xi_{\xi}$ is equal to $\RV^k\times \bigcup_{\lambda \le 1}\RV_\lambda^{\ell}$, and then compactness implies that finitely many sets $\Xi_{\xi_i}$ suffice to cover everything. Now \(C := \bigcup_i C_{\xi_i}\) is a finite \(A\)-definable set which prepares every \(W_\xi\).
\end{remark}

\begin{remark}\label{rem:non:compl}
A variant of the compactness argument shows that Propsosition~\ref{prop:uniform} holds even more uniformly, namely: Given an $\ell$-h-minimal theory $\cT$ (possibly non-complete)
and a formula $\phi$ defining a set
\[
W_{K,a} := \phi(K,a) \subset K\times \RV_K^{k}\times  \bigcup_{\lambda \le 1}\RV_\lambda^{\ell}
\]
for $K \models \cT$ and $a \in K^\nu \times \RV_K^{\nu'}$,
there exists a formula $\psi$
defining a set $C_{K,a} := \psi(K,a) \subset K$ such
that for each model $K \models \cT$ and each $a \in K^\nu \times \RV_K^{\nu'}$, this set $C_{K,a}$ uniformly prepares $W_{K,a}$.
\end{remark}

\begin{remark}\label{rem:finiterange:set}
In Proposition~\ref{prop:uniform}, we can also replace $\RV^k$ by
any Cartesian product $Q$ of sorts from $\RV\eq$.
Indeed, in that case, just apply the original version of the proposition to the preimage of $W$ under some quotient map
$\RV^k \to Q$.
In particular, $W$ can additionally use (arbitrarily many) $\Gamma_K$-coordinates, since $\Gamma_K$ is a quotient of $\RV$. (One can of course also allow $\RF$-coordinates, given that $\RF$ can be considered as a $\emptyset$-definable subset of $\RV$.)
\end{remark}

In almost all applications of Proposition~\ref{prop:uniform}, we will only need the following special case:

\begin{cor}[Preparing families]\label{cor:prep}
Assume that $\Th(K)$ is $0$-h-minimal and that $A$ is a subset of $K$. For any $k>0$ and any $(A \cup \RV)$-definable set
$$
W\subset K\times \RV^k,
$$
there exists a finite non-empty $A$-definable set $C \subset K$
$1$-preparing $W$ (in the sense of Definition~\ref{defn:uniform}).
\end{cor}

\begin{proof}[Proof of Corollary~\ref{cor:prep}]
Clear.
\end{proof}

Recall that we set $\RV_\bullet := \bigcup_{\lambda \le 1} \RV_\lambda$.

\begin{cor}[$\RV$-unions stay finite]\label{cor:finiterange:set}
\begin{enumerate}
 \item
Assume that $\Th(K)$ is $1$-h-minimal.
For any $k \ge 0$ and any definable (with parameters) set $W\subset \RV_\bullet^k \times K$ such that the fiber $W_\xi\subset K$ is finite for each $\xi\in \RV_\bullet^{k}$, the union $\bigcup_{\xi} W_\xi$ is also finite.
\item
Under the (weaker) assumption that $\Th(K)$ is $0$-h-minimal, the previous statement still holds if we assume $W\subset \RV^k \times K$.
\end{enumerate}
\end{cor}

\begin{proof}
We suppose that $W\subset \RV^k \times \RV_\bullet^{k'} \times K$
and proceed by induction on $k'$. If $k' \le 1$, we obtain both claims of the corollary from Proposition~\ref{prop:uniform} applied with $\ell = k'$, namely: We find a finite set \(C\) such that, for all \(\lambda\) and \(\xi\in \RV^k\times\RV_\lambda^{k'}\), $W_\xi$ is $\lambda$-prepared by \(C\). Since \(W_\xi\) is finite, we have \(W_\xi\subset C\), and hence \(\bigcup_\xi W_\xi \subset C\) is finite.

The case of $k'>1$ now follows by induction on $k'$ and the case $k'=1$.
\end{proof}

\begin{cor}[Finite image in $K$]\label{cor:finiterange}
\begin{enumerate}
 \item
Assume that $\Th(K)$ is $1$-h-minimal. Then the image of any definable (with parameters) function $f\colon  \RV_\bullet^{k} \to K$ for any $k \ge 0$ is finite.
\item
Under the (weaker) assumption that $\Th(K)$ is $0$-h-minimal, the previous statement still holds for functions $f\colon  \RV^{k} \to K$.
\end{enumerate}
\end{cor}

\begin{proof}
Apply Corollary \ref{cor:finiterange:set} to the graph of $f$.
\end{proof}

\begin{remark}\label{rem:prep}
Remark~\ref{rem:finiterange:set} also applies to Corollaries \ref{cor:prep}, \ref{cor:finiterange:set}, and \ref{cor:finiterange}. In particular, we can additionally allow (arbitrarily many) $\Gamma_K$-coordinates in $W$
(in \ref{cor:prep}, \ref{cor:finiterange:set}) and in the domain of $f$ (in \ref{cor:finiterange}).
\end{remark}

\begin{cor}[Removing $\RV$-parameters]\label{cor:acl}
Assume that $\Th(K)$ is $0$-h-minimal. For any $A \subset K$ and any finite $(A \cup \RV\eq)$-definable set $C \subset K$, there exists a finite $A$-definable set $C' \subset K$ containing $C$. In other words, $\acl_K(A \cup \RV\eq) = \acl_K(A)$.
\end{cor}

\begin{proof}
Add constants for $A$ to the language.
We have $C = W_ {\xi_0}$ for some $\emptyset$-definable $W \subset K \times \RV^k$ and some $\xi_0 \in \RV^k$. We may assume  that all fibers $W_{\xi}$ have cardinality at most the cardinality of $C$. Let $C'$ be their union, which is finite by Corollary~\ref{cor:finiterange:set}.
\end{proof}

\begin{remark}\label{rem:acl}
Many results in this paper are stated in the form:
\begin{itemize}
 \item[$(\star)$]
For every $\emptyset$-definable object $X$ of a certain kind, there exists a finite $\emptyset$-definable set $C \subset K$ which ``prepares'' $X$ in some sense (depending on the context). \end{itemize}
By Lemma~\ref{lem:addconst}, we get for free that
$(\star)$ holds more generally, namely if $X$ is $(A \cup \RV)$-definable, for $A \subset K$, we get an $(A \cup \RV)$-definable $C$ ``preparing'' $X$.
By applying Corollary~\ref{cor:acl}, we then may even assume that $C$ is $A$-definable. (It will always be the case that if $C$ prepares $X$, then so does any set containing $C$.) Finally, using that the notions of preparation under consideration will be definable, we can apply compactness to get for free that this works uniformly in all models of a non-complete theory, in the same style as in Remark~\ref{rem:non:compl}.
\end{remark}

We end this subsection by noting that $\RV$ is stably embedded in a strong sense (namely with the $\RV$-parameters being in the definable closure of the original parameters):

\begin{prop}[Stable embeddedness of $\RV$]\label{prop:stab}
Assume that $\Th(K)$ is $0$-h-minimal. Then
$\RV$ is stably embedded in the following strong sense: Given any $A \subset K$, every $A$-definable set $X \subset \RV^n$
is already $\dcl_{\RV}(A)$-definable.
\end{prop}

\begin{proof}
We may assume that $A$ is finite; we do an induction on the cardinality of $A$.

Let $A = \hat A \cup \{a\}$. Then we have an $\hat A$-definable set $Y \subset K \times \RV^n$ such that
$X$ is equal to the fiber $Y_a \subset \RV^n$.
By applying Corollary~\ref{cor:prep} to $Y$, we find a finite $\hat A$-definable set $C \subset K$ such that
either $a \in C$, or for every $a' \in K$ in the same ball $1$-next to $C$ as $a$, we have $Y_{a'} = Y_a$.
Using Lemma~\ref{lem:InextFam}, we find an $\hat A$-definable map $f\colon K \to \RV^k$ (for some $k$)
whose fibers are exactly the elements of $C$ and the balls $1$-next to $C$. In particular, the set $X = Y_a$ is definable
using $\hat A$ and $f(a)$ as parameters. Thus we have $X = Z_{f(a)}$ for some $\hat A$-definable set
$Z \subset \RV^k \times \RV^n$. By induction, $Z$ is $\dcl_{\RV}(\hat A)$-definable, so $X$ is $\dcl_{\RV}(\hat A)\cup \{f(a)\}$-definable and hence $\dcl_{\RV}(A)$-definable.
\end{proof}

\begin{remark}\label{rem:stab}
If $\Th(K)$ is $\omega$-h-minimal, there are also various variants of Proposition~\ref{prop:stab} involving $\RV_\lambda$, with similar proofs. For example, building on Remark~\ref{rem:InextFam} instead of Lemma~\ref{lem:InextFam}, one obtains that
any $A$-definable subset of $\prod_i \RV_{\lambda_i}$ (for $A \subset K$) is
$\dcl_{\RV_{\lambda}}(A)$-definable with $\lambda = \min_i \lambda_i$.
\end{remark}

\subsection{Henselianity of the valued field $K$}\label{sec:hens}

As an analogue of o-minimal fields being real closed, in this subsection, we prove that any equi-characteristic zero valued field which is Hensel minimal (in any language containing $\Lval$) is Henselian. This is one reason why we call our notion ``Hensel minimality''.

A collection of balls is called nested if for any two balls in the collection, one is contained in the other.

\begin{lem}[Definable spherical completeness]\label{lem:intersection}
Assuming $0$-h-minimality, let $\{B_q \mid q\in Q\}$ be a definable family of nested balls in $K$,
for some non-empty definable set $Q$ in an arbitrary, possibly imaginary sort. Then the intersection $\bigcap_{q \in Q} B_q$ is non-empty.
\end{lem}
\begin{proof}
First we suppose that $Q \subset \Gamma^\times_K$ and that each $B_q$ is an open ball of radius $q$.
By Corollary~\ref{cor:prep} (and Remark~\ref{rem:prep}), there exists a finite set $C$ $1$-preparing the family of balls $B_q$.
Now one checks that at least one of the following two situations occurs.
Firstly: For each $q\in Q$, the intersection of $C$ with $B_q$ is non-empty. Secondly: The set $Q$ has a minimum $q_0\in Q$. (Indeed, suppose that the intersection of $C$ with $B_{q_0}$ is empty for some $q_0\in Q$; then $Q$ contains no $q < q_0$, since $B_q$ would not be $1$-prepared by $C$.) In both situations, the lemma follows.

Finally, we reduce the general case to the case that $Q\subset \Gamma^\times_K$ and that each $B_q$ is an open ball of radius $q$. To this end, for $\gamma\in\Gamma^\times_K$, let $B'_{\gamma}$ be the (necessarily unique) open ball of radius $\gamma$ containing some $B_q$ ($q \in Q$) if such a ball exists, and let it be the empty set otherwise. Then it is clear that the non-empty $B'_{\gamma}$ form a nested definable family of open balls. Moreover, the intersection of the non-empty $B'_\gamma$ equals the intersection of the $B_q$ (since each $B_q$ is equal to the intersection of all open balls containing $B_q$).
\end{proof}

\begin{thm}[Hensel minimality implies Henselian]\label{thm:hens}
Suppose that $K$ is a valued field of equi-characteristic $0$ with $0$-h-minimal theory (in a language $\cL \supset \Lval$). Then $K$ is Henselian.
\end{thm}
If $\cL$ is the pure valued field language, Corollary~\ref{cor:ex:equi} implies the converse. Combining, we have, for $K$ of equi-characteristic $0$:
$K$ is Henselian if and only if $\Th_{\Lval}(K)$ is $0$-h-minimal, if and only if
$\Th_{\Lval}(K)$ is $\omega$-h-minimal.

\begin{proof}[Proof of Theorem~\ref{thm:hens}]
Let $P\in\cO_K[X]$ be a polynomial such that $P(0)\in\cM_K$ and $P'(0)\in\cO^{\times}_K$; we need to prove that $P$ has a root in $\cM_K$. (Note that uniqueness of such a root then follows automatically.)
The idea is to use ``Newton approximation'' as in the usual proof of Hensel's lemma for complete discretely valued fields, but where complete and discretely valued is replaced by definably spherically complete.

To make this formal, we suppose that $P$ has no root in $\cM_K$ and we set $B_x := B_{\le |P(x)|}(x)$ for $x \in \cM_K$. Note that the $B_x$ form a definable family of balls, parameterized by $x \in \cM_K$. We will prove that (a) all these balls form a chain under inclusion and (b) that an element in the intersection of all those balls (which is non-empty by Lemma~\ref{lem:intersection}) is, after all, a root of $P$.

(a) Let $x_1, x_2 \in \cM_K$ be given and set $\varepsilon := x_2 - x_1$. To see that the balls $B_{x_1}$ and $B_{x_2}$ are not disjoint, we verify that $|\epsilon| \le \max\{|P(x_1)|, |P(x_2)|\}$. Taylor expanding $P$ around $x_1$ yields
\begin{equation}\label{eq:tay}
|P(x_1+\varepsilon) - P(x_1) - \varepsilon P'(x_1)| \le |\varepsilon^2|,
\end{equation}
which, together with $|P'(x_1)| = 1$, implies $|\epsilon| \le \max\{|P(x_1)|, |P(x_1 + \epsilon)|\}$.

(b) Let $x_1$ be in the intersection $\bigcap_{x \in \cM_K} B_x$ and suppose that $P(x_1) \ne 0$. Then (\ref{eq:tay}) with $\varepsilon := -  \frac{P(x_1)}{P'(x_1)}$ implies $|P(x_1 + \varepsilon)| \le |\varepsilon^2| < |\varepsilon|$ and hence $x_1 \notin B_{x_1 + \varepsilon}$, contradicting our choice of $x_1$. \end{proof}

\begin{remark}
Lemma~\ref{lem:intersection} implies a ``definable Banach Fixed Point Theorem'' (exactly in the form of \cite[Lemma~2.32]{Halup}, and with the same proof).
The above proof of Theorem~\ref{thm:hens} can be considered as applying that Fixed Point Theorem to
the map $\cM_K \to \cM_K, x \mapsto x - P(x)/P'(x)$.
\end{remark}

\private{This subsection directly works in finitely ramified mixed char. With infinite ramification, I'm unsure.}

\subsection{Definable functions}
\label{sec:fctn}

We continue assuming that $K$ is an equi-characteristic $0$ valued field, and we now assume that $\Th(K)$ is $1$-h-minimal (unless specified otherwise). Under those assumptions, we now prove first basic properties of definable functions in one variable; in particular, we already obtain a weak version of the Jacobian Property  (Lemma~\ref{lem:gammaLin}) and simultaneous domain and image preparation (Proposition~\ref{prop:range}).

\private{In case we want to check again how much $1$-h-min we really need: In this section, $1$-h-min is used only as follows: (1) in \ref{lem:fin-inf}; Use \ref{cor:finiterange:set};
(2) in \ref{lem:balltosmall}: use \ref{cor:prep}.}

The first result is key to dimension theory (though in our proofs of dimension theory,
this will only be used indirectly, namely in the proof of Proposition \ref{prop:b-min}).
\begin{lem}[Basic preservation of dimension]\label{lem:fin-inf}
Assume (as convened for the whole Section \ref{sec:fctn}) that  $\Th(K)$ is $1$-h-minimal.
Let $f\colon K \to K$ be a definable function.
Then there are only finitely many function values which are taken infinitely many times.
\end{lem}

\begin{proof}
We may assume $f$ to be $\emptyset$-definable (say, after adding enough parameters from $K$ to the language).
Suppose $f$ takes infinitely many values $y$ infinitely many times.
Then for each such $y$,
$f^{-1}(y)$ contains a ball (by Lemma~\ref{lem:finite}). Thus, letting $X \subset K$ be the set of points where $f$ is locally constant, $f(X)$ is still infinite.

Let $W\subset K\times \Gamma^\times_K$ consist of those $(x,\lambda)$ such that $f$ is constant on $B_{<\lambda}(x)$. This set $W$ is $\emptyset$-definable, so we find a finite set $C$ $1$-preparing $W$ (by Corollary \ref{cor:prep}).
By enlarging $C$, we may moreover assume that $C$ also $1$-prepares $X$.

From the fact that $f(X)$ is infinite, we can deduce that there exists a ball $B_0 \subset X$ $1$-next to $C$ such that $f(B_0)$ is still infinite. Indeed, letting $g\colon K \to \RV^k$ be a map whose fibers are the singletons in $C$ and the balls $1$-next to $C$ (using Lemma~\ref{lem:InextFam}), if $f(g^{-1}(\xi))$ would be finite for every $\xi \in g(X)$, then so would be $f(X)$ (by Corollary~\ref{cor:finiterange:set}).

Choose $x \in B_0$ and $\lambda_0 \in \Gamma^\times_K$ such that $f$ is constant on $B_{<\lambda_0}(x)$. Since $C$ $1$-prepares $W$, $f$ is constant on $B_{<\lambda_0}(x')$ for every $x' \in B_0$.

Set $\lambda_1 := \lambda_0/\radop(B_0)$. Then the family $F_1$ of open balls of radius $\lambda_0$ contained in $B_0$ can be definably parametrized by a subset of $\RV_{\lambda_1}$
(using some parameters). Indeed, if we fix $c \in K$ such that $B_0$ is $1$-next to $c$, then each member of $F_1$ is of the form $c + \rv_{\lambda_1}^{-1}(\xi)$ for some
$\xi \in \RV_{\lambda_1}$. We define $F_2$ to be the family of $f(B)$, for $B$ in $F_1$. Then each family member of $F_2$ is a singleton, yet their union is infinite, contradicting Corollary \ref{cor:finiterange:set}.
\end{proof}

Using this, we obtain that definable functions are (in a strong sense) locally constant or injective:

\begin{lem}[Piecewise constant or injective]\label{lem:loc-inj}
For every $\emptyset$-definable map $f\colon K \to K$, there exists a finite $\emptyset$-definable set $C$ such that
for every ball $B$ $1$-next to $C$, $f$ is either constant or injective on $B$.
\end{lem}

\begin{proof}
First, consider the set $Y_\infty$ of $y \in K$ such that $f^{-1}(y)$ is infinite. By Lemma~\ref{lem:fin-inf}, this set $Y_\infty$ is finite, so we can find
a finite $\emptyset$-definable set $C$ $1$-preparing $f^{-1}(y)$ for each $y \in Y_\infty$. Indeed, for each $y \in Y_\infty$, one finds a finite $y$-definable
set $C_y$ $1$-preparing $f^{-1}(y)$ (uniformly in $y$), and one lets $C$ be the union of the sets $C_y$.

For each $y \in K \setminus Y_\infty$, the set $f^{-1}(y)$ is finite, so by Lemma~\ref{lem:average},
there exists a $\emptyset$-definable
family of injective functions $g_y\colon f^{-1}(y) \to \RV^k$ for some $k \ge 1$. For convenience, we set $g_y(x) := 0$ if $y \in Y_\infty$ and $x \in f^{-1}(y)$,
so that we can define a function $h\colon K \to \RV^k$ by $h(x) := g_{f(x)}(x)$.
We then enlarge our above set $C$ (using  Corollary~\ref{cor:prep}) so that it also $1$-prepares (the graph of) $h$.

We claim that this set $C$ is as desired, so let $B$ be a ball $1$-next to $C$.
If $f(B) \cap Y_\infty \ne \emptyset$, then $B$ is contained in one of the sets $f^{-1}(y)$, for $y \in Y_\infty$, and hence $f$ is constant on $B$.
Otherwise, we use that $h$ is constant on $B$ to deduce that $f$ is injective on $B$. Indeed, $f(x_1) = f(x_2) = y$ for some $y \in K \setminus Y_\infty$
implies $h(x_1) = g_{y}(x_1) = g_{y}(x_2) = h(x_2)$,
so injectivity of $h$ on $f^{-1}(y)$ implies $x_1 = x_2$.
\end{proof}

The next lemma says that a definable function sends most open balls either to points or to open balls.

\begin{lem}[Images of balls are balls]\label{lem:open-to-open} Let $f\colon K \to K$ be a $\emptyset$-definable function.
There exists
a $\emptyset$-definable finite set $C \subset K$ such that for every open ball $B$ disjoint from $C$,
$f(B)$ is either a point or an open ball.
\end{lem}

\begin{proof}
Define $W \subset K \times \Gamma^\times_K$ to consist of those $(x, \lambda)$ for which the open ball $B_{<\lambda}(x)$ ``has a good image'', i.e., $f(B_{<\lambda}(x))$ is a singleton or an open ball.
By  Corollary \ref{cor:prep}, we find a finite $\emptyset$-definable subset $C_0$ of $K$ such that for any ball $B$ $1$-next to $C_0$, the fiber $W_x \subset \Gamma^\times_K$ does not depend on $x \in B$.

Fix a ball $B_0$ $1$-next to $C_0$. We first prove that any open ball $B$ strictly contained in $B_0$ has a good image.
Suppose otherwise, i.e., that $B_1$ is an open ball strictly contained in $B_0$ with bad image.
Then the fact that $W_x$ does not depend on $x \in B_0$ implies that every translate of $B_1$ contained in $B_0$ also has bad image.
We can find an infinite definable (with parameters) family $F_1$ consisting of such translates of $B_1$ and parameterized by
a subset of $\RV$. Indeed, the sets $x_0 + \rv^{-1}(\xi)$ form such a family for a suitable $x_0 \in K$ and when $\xi$ runs over a suitable subset of $\RV$.

Consider the family $F_2$ of the sets $f(B)$ for $B$ in $F_1$,
use Corollary \ref{cor:prep} to find a finite set $D \subset K$ $1$-preparing the family $F_2$, and let $g\colon K \to \RV^k$ be a function whose fibers are the balls $1$-next to $D$ and the individual points of $D$ (as obtained using Lemma~\ref{lem:InextFam}).
Since none of the balls $B$ in $F_1$ are good, none of the sets $f(B)$ in $F_2$ are exactly a point or an open ball, so each $f(B)$ consists of several fibers of $g$; in other words, $g \circ f\colon K \to \RV^k$ is non-constant on every $B$ in $F_1$.

Now we get a contradiction by applying Corollary \ref{cor:prep} to the graph of $g \circ f$. Indeed, any set $C'$ $1$-preparing
that graph would have to contain at least one point in each ball $B$ from $F_1$ (since $g \circ f$ is not constant on $B$), so $C'$ cannot be finite.
This finishes the proof that balls strictly contained in $B_0$ have good image.

The only problematic open balls left (i.e., which are disjoint from $C_0$ and might have bad image) are the ones $1$-next to $C_0$. To get hold of those, we run a similar argument as above: We let $F_1$ be the family of balls $1$-next to $C_0$; we find a finite set $D$ $1$-preparing $f(B)$ for each
$B$ in $F_1$, and we define $g\colon K \to \RV^k$ as before, so that in particular
if $g \circ f$ is constant on a ball $B$ $1$-next to $C_0$, then that ball $B$ has good image.

Now we find a finite set $C'$ $1$-preparing the graph of $g \circ f$ and we set $C := C_0 \cup C'$. In this way, among the balls $1$-next to $C_0$, all those which have bad image are not disjoint from $C$.
Note that since this time, $F_1$ is $\emptyset$-definable, and hence so are $D$, $g$ and $C'$.
\end{proof}

The following is a key technical lemma, which serves later in the proof of the Jacobian Property. The main point of the statement is that a definable map cannot scale all small balls by one factor and all big balls by a different factor.

\begin{lem}[Preservation of scaling factor]\label{lem:balltosmall} 
Let $B$ be either $\cO_K$ or $\cM_K$ and
let $f\colon B\to K$ be a $\emptyset$-definable function.
Suppose that there are $\alpha$ and $\beta$ in $\Gamma^\times_K$
with
$\alpha<1$
such that for every open ball $B' \subset B$ of radius $\alpha$, $f(B')$ is contained in an open ball of radius $\beta$. Then the following hold.
\begin{enumerate}
 \item
The image $f(B)$ is contained in a finite union of closed balls of radius $\beta/\alpha$.
\item
If we moreover assume that $B = \cM_K$ and that $f(B)$ is an open ball, then $\radop f(B) \le \beta/\alpha$.
\end{enumerate}
\end{lem}

\begin{proof}
(1)
The family of all open balls of radius $\alpha$ contained in $B$ can be (definably) parametrized by a definable set $\Lambda \subset \RV_\alpha$, namely, set $\Lambda := \rv_\alpha(B)$, and for $\xi \in \Lambda$, let $B_\xi$ be the open ball of radius $\alpha$ containing $\rv_\alpha^{-1}(\xi)$ (which exists, since $\radop(\rv_\alpha^{-1}(\xi)) \le \alpha$).
By Proposition~\ref{prop:uniform}, there exists a finite set $C$ such that for each $\xi\in\Lambda$, the set $f(B_\xi)$ is $\alpha$-prepared by $C$. This implies that for
every open ball $B'$ which
is $\alpha$-next to $C$ and every $\xi\in\Lambda$, one has
$$
\text{either } B'\subset f(B_\xi), \quad\text{or } B'\cap  f(B_\xi)=\emptyset.
$$
Hence, if the radius of the open ball $B'$ is larger than $\beta$, then $B'\cap  f(B_\xi)$ is empty for all $\xi\in\Lambda$. (Indeed, by assumption, $f(B_\xi)$ is contained in an open ball of radius $\beta$, so $f(B_\xi)$ cannot contain the larger ball $B'$.)

Thus, we find that $f(B)$ is contained in the union of $C$ with those open balls $\alpha$-next to $C$ that have radius at most $\beta$. This union
equals the union over $c\in C$ of the closed balls $B_{\le \beta/\alpha}(c)$. This proves (1).

(2)
If the value group is dense, then (2) follows from (1) using that the largest open balls contained in a finite
union of closed balls of radius $\beta/\alpha$ have radius $\beta/\alpha$.

If the value group is discrete (see just above Remark \ref{rem:rad-op}), we apply Part (1) to $g(x) := f(\varpi x)\colon \cO_K \to K$, where $\varpi \in K$ is a
uniformizing element. Since $g$ sends open balls of radius $|\varpi|^{-1} \cdot \alpha$ to open balls of radius $\beta$,
we obtain that
$f(\cM_K) = g(\cO_K)$ is contained in a finite union of closed balls of radius $|\varpi| \cdot \beta/\alpha$, which
is the same as a finite union of open balls of radius $\beta/\alpha$. Now the claim follows from the assumption that
$f(\cM_K)$ is itself an open ball.
\end{proof}

By combining various previous lemmas, we already obtain a weak form of the Jacobian Property.

\begin{lem}[Valuative Jacobian Property]\label{lem:gammaLin}
For every $\emptyset$-definable function $f\colon K \to K$, there exists a finite $\emptyset$-definable set $C \subset K$ such that
for every ball $B$ $1$-next to $C$, we have the following:
Either $f$ is constant on $B$, or there exists a $\mu_B \in \Gamma^\times_K$ such that
\begin{enumerate}
 \item
for every open ball $B' \subset B$, $f(B')$ is an open ball of radius $\mu_B \cdot\radop(B')$;
 \item for every $x_1, x_2 \in B$, we have $|f(x_1) - f(x_2)| = \mu_B\cdot |x_1 - x_2|$.
\end{enumerate}
\end{lem}

\begin{proof}
Choose $C$ $\emptyset$-definable such that for every ball $B$ $1$-next to $C$, we have:
\begin{itemize}
 \item $f$ is constant or injective on $B$ (using Lemma~\ref{lem:loc-inj});
 \item $f(B')$ is a point or an open ball for every open ball $B' \subset B$ (using Lemma~\ref{lem:open-to-open}).
\end{itemize}
Moreover, we assume (using Corollary~\ref{cor:prep}) that $C$ $1$-prepares the graph of the function $r\colon K \times \Gamma^\times_K \to \Gamma_K$
defined by
\[r(x,\lambda) = \begin{cases}
                  \radop(f(B_{<\lambda}(x))) & \text{if } f(B_{<\lambda}(x)) \text{ is an open ball}\\
                  0 & \text{otherwise}.
                 \end{cases}
\]
We claim that this $C$ is as desired, so fix a ball $B$ $1$-next to $C$ for the remainder of the proof.
Note that for $x$ running over $B$, the function $r(x,\cdot)$ is independent of $x$,
so from now on we write $r(\lambda)$ instead of $r(x, \lambda)$.

If $f$ is constant on $B$, there is nothing to do, so we may assume that $f$ is injective on $B$.
We claim that the lemma holds with $\mu_B := \radop(f(B))/\radop(B)$. (Note that indeed, $f(B)$ is an open ball.)

(1) Fix $\alpha \in \Gamma_K$ with $0 < \alpha < \radop(B)$.
Any open ball $B' \subset B$ of radius $\alpha$ is sent to an open ball $f(B')$ of radius $r(\alpha)$.
By applying Lemma~\ref{lem:balltosmall} (to a suitably rescaled function), we deduce that the radius of the open ball $f(B)$ is at most
$\radop(B)\cdot r(\alpha)/\alpha$; this implies $r(\alpha)/\alpha \ge \mu_B$.

To get the other inequality, namely $r(\alpha)/\alpha \le \mu_B$, we apply the same argument to the inverse function
$f^{-1}\colon f(B) \to B$. This inverse is well-defined since $f$ is injective on $B$, so it remains
to verify that $f^{-1}$ sends open balls $B'' \subset f(B)$ of radius $r(\alpha)$ to open balls of radius $\alpha$.
Indeed, choose any $x \in f^{-1}(B'')$; then $f(B_{< \alpha}(x)) = B''$
and hence $f^{-1}(B'') = B_{< \alpha}(x)$.

(2) For every $x_1, x_2 \in B$, (1) implies $|f(x_1) - f(x_2)| \le \mu_B\cdot |x_1 - x_2|$. Indeed,
applying (1) to a ball of the form $B_{<\alpha}(x_1)$ with $|x_1 - x_2| < \alpha \le \radop(B)$
yields $|f(x_1) - f(x_2)| < \mu_B\cdot \alpha$ for every $\alpha > |x_1 - x_2|$, and hence
$|f(x_1) - f(x_2)| \le \mu_B\cdot |x_1 - x_2|$.
The same argument with $f$ replaced by $f^{-1}\colon f(B) \to B$ yields the other inequality.
\end{proof}

Using Lemma~\ref{lem:gammaLin}, we deduce that the domain and the image of a definable function can be prepared simultaneously and in a compatible way.

\begin{prop}[Domain and range preparation]\label{prop:range}
Let $f\colon K \to K$ be a $\emptyset$-definable function
and let $C_0 \subset K$ be a finite, $\emptyset$-definable set.
Then there exist finite, $\emptyset$-definable sets $C, D \subset K$ with $C_0 \subset C$ such that $f(C) \subset D$, and for every $\lambda \le 1$ in $\Gamma^\times_K$ and every ball $B \subset K$ that is $\lambda$-next to $C$, the image $f(B)$ is either a singleton in $D$ or a ball $\lambda$-next to $D$.
\end{prop}

The role of the $C_0$ is to make it possible to combine this proposition with other preparation results for functions (notably Theorem~\ref{thm:high-ord}): First apply the other results to get a set $C_0$; then apply Proposition~\ref{prop:range} to enlarge $C_0$ to $C$ and to obtain $D$.

Note however that the proposition cannot be combined very well with itself:
Given $f_1, f_2\colon K \to K$ as in the proposition, there are, in general, no $C, D_1, D_2$ such that $(f_i, C, D_i)$ are as in the proposition for both $i = 1, 2$, as the following example shows.

\begin{example}
Fix $r \in K$ of negative valuation. For $x \in \cO_K$, we set $f_1(x) = f_2(x) = f_1(x + r) = x$ and $f_2(x + r) = x+1$, with the exception that $f_2(0) = r$. Extend $f_1$ and $f_2$ by $0$ outside of $\cO_K \cup (\cO_K + r)$. Then one successively deduces: $0 \in C \Rightarrow 0 \in D_1 \Rightarrow r \in C \Rightarrow 1 \in D_2 \Rightarrow 1 \in C \Rightarrow 1 \in D_1 \Rightarrow \dots$; this shows that $C$ cannot be finite.
\end{example}

In Addendum \ref{add:cd:alg:range} to the cell decomposition Theorem~\ref{thm:cd:alg:skol}, we will state a version of simultaneous preparation of the domain and image which avoids this problem (and thus works for several functions simultaneously), by working piecewise.

\begin{proof}[Proof of Proposition~\ref{prop:range}]
First, enlarge $C_0$ to a set $C_1$ using Lemma~\ref{lem:gammaLin}, so that on each ball $B$ $1$-next to $C_1$, $f$ sends open balls to open balls and the quotient
\begin{equation}\label{eq:quoco}
|f(x_1) - f(x_2)|/|x_1 - x_2|
\end{equation}
is constant. Then let $D$ be a set containing $f(C_1)$ and preparing the family $f(B)$, where $B$ runs over the balls $1$-next to $C_1$. (This family is parametrized by $\RV$-variables, by Lemma~\ref{lem:InextFam}, so Corollary~\ref{cor:prep} applies.) Now
denote by $X \subset K$ the locus where $f$ is not locally constant, and let $C$ be $\big(X\cap f^{-1}(D)\big) \cup C_1$. Then clearly $C$ is finite. We claim that these $C$ and $D$ are as desired.

Note that we have $C \supset C_1$, so the only thing that is not clear from the construction is that balls $\lambda$-next to $C$ are sent to elements of $D$ or to balls
$\lambda$-next to $D$. We first deal with the case $\lambda = 1$.

Let $B$ be a ball $1$-next to $C$ and let $B_1$ be the ball $1$-next to $C_1$ containing $B$. If $f(B_1)$ is a singleton, then it is an element of $D$ and we are done. If $f(B_1)$ is not a singleton, then $B_1\subset X$ by the constancy of (\ref{eq:quoco}) on $B_1$.
If furthermore $f(B_1) \cap D = \emptyset$, then $B_1 \cap C = \emptyset$, so $B = B_1$ and $f(B_1)$ is a ball $1$-next to $D$.
Finally, suppose that $f(B_1) \cap D \ne \emptyset$. Then $B_1 \cap C \ne \emptyset$, so $B$ is $1$-next to some $c \in C \cap B_1$. Then using that (\ref{eq:quoco}) is constant on $B_1$, we obtain that $f(B)$ is a ball $1$-next to $f(c) \in D$.

Now suppose that $\lambda < 1$. Any ball $B$ $\lambda$-next to $C$ is contained in a ball $B'$ $1$-next to $C$. If $f$ is constant on $B'$, we are done; otherwise, $f(B')$ is a ball $1$-next to $D$, and we deduce that $f(B)$ is $\lambda$-next to $D$, using once more that (\ref{eq:quoco}) is constant on $B'$.
\end{proof}

\subsection{An equivalent condition to $1$-h-minimality}\label{sec:ta}

The conclusions of Lemmas~\ref{lem:fin-inf} and \ref{lem:gammaLin} together are actually equivalent to $1$-h-minimality.
More precisely, we have the following equivalence, where Condition~(T1) is slightly weaker than Lemma~\ref{lem:gammaLin}:

\begin{thm}[Criterion for $1$-h-minimality]\label{thm:tame2vf}
Let $\cL$ be a language containing the language $\Lval$ of valued fields, and let $\cT$ be an
$\cL$-theory whose models are valued fields of equi-characteristic $0$.
Then $\cT$ is $1$-h-minimal if and only if, for every model $K$ of $\cT$, for every set $A \subset K \cup \RV$ and for every
$\cL(A)$-definable $f\colon K \to K$, we have the following:
\begin{enumerate}[(T1)]
\item There exists a finite $\cL(A)$-definable $C \subset K$ such that for every ball $B \subset K$ $1$-next to $C$,
there exists a $\mu_B \in \Gamma_K$ such that for every $x_1, x_2 \in B$, we have $|f(x_1) - f(x_2)| = \mu_B \cdot |x_1 - x_2|$.
\item The set $\{d \in K \mid f^{-1}(d)$ is infinite$\}$ is finite.
\end{enumerate}
\end{thm}

Note that in the above, we intentionally allow $C$ to use the parameters from $\RV$ in its definition, in contrast to the condition in the definition of $1$-h-minimality.

\begin{remark}\label{rem::tame2vf}
The conditions given in Theorem \ref{thm:tame2vf} are very closely related to the tameness notion of Definition~2.1.6 of \cite{CCL-PW} and its variant from Section~2.1 of \cite{CFL}, one difference being that we do not assume the existence of an angular component map.
\end{remark}

By Lemmas~\ref{lem:gammaLin} and \ref{lem:fin-inf}, (T1) and (T2) follow from $1$-h-minimality, so in the remainder of Subsection~\ref{sec:ta}, we assume that $\cT$ satisfies (T1) and (T2), our goal being to prove $1$-h-minimality.
We also assume throughout the subsection that $K$ is a model of $\cT$.

\begin{remark}\label{rem:nocontrol}
By applying (T1) to the characteristic function of an $A$-definable set $X \subset K$ (for $A \subset K \cup \RV$), we find a finite $A$-definable set $C \subset K$ $1$-preparing $X$. In particular, Assumption~\ref{ass:no-ctrl} is satisfied, so we may apply Lemmas~\ref{lem:finite} ($\exists^\infty$-elimination), \ref{lem:average} (existence of injective functions from finite sets to $\RV$) and \ref{lem:InextFam} (the balls $1$-next to $C$ are $\RV$-parametrized).
\end{remark}

\begin{lem}\label{lem:tame2finite}
Assume (as convened for the remainder of Subsection~\ref{sec:ta}) that $\cT$ satisfies (T1) and (T2) from Theorem \ref{thm:tame2vf} and that $K$ is a model.
We have the following:
\begin{enumerate}
\item Any definable (with parameters) function $f\colon\RV^k \to K$ has finite image.
\item If $C_\xi \subset K$ is a definable (with parameters) family of finite sets, parametrized by $\xi \in \RV^k$, then
the union $\bigcup_{\xi \in \RV^k} C_\xi$ is still finite.
\end{enumerate}
\end{lem}

\begin{proof}
We first prove both claims for $k = 1$.

(1) The composition $f \circ \rv\colon K \to K$ is locally constant everywhere except possibly at $0$, so the claim follows from
(T2).

(2)
By Lemma~\ref{lem:finite}, the cardinality of $C_\xi$ is bounded independently of $\xi$,
so we may assume that the cardinality is constant, say, equal to $m$.
Let $\sigma_1, \dots, \sigma_m \in \ZZ[x_1, \dots, x_m]$ be the elementary symmetric functions in $m$ variables, which
we consider as functions on the set of $m$-element subsets of $K$.
By (1), for each $i$, the function $f_i(\xi) := \sigma_i(C_\xi)$ has finite image. Since $\sigma_1(C), \dots, \sigma_m(C)$ together determine $C$,
there are only finitely many different sets $C_\xi$, which implies the claim.

Now we deduce (2) for arbitrary $k$ by induction: Given a definable family $C_{\xi, \xi'} \subset K$, for $\xi \in \RV$ and $\xi' \in \RV^k$, we first obtain that $D_{\xi'} := \bigcup_{\xi \in \RV} C_{\xi, \xi'}$ is finite for every $\xi'$ and then that the entire union is finite.

Finally, we obtain (1) for arbitrary $k$ by applying (2) to the family of singletons $C_\xi := \{f(\xi)\}$.
\end{proof}

\private{All the stuff above the following lemma does not use any control of parameters, i.e., in (1), we wouldn't need to require $C$ to be $A$-definable.}

\begin{lem}\label{lem:tame-0h}
Given an $A$-definable function $f\colon K \to K$, with $A \subset K \cup \RV$, we can find $C$ as in Condition (T1) of Theorem~\ref{thm:tame2vf} which is moreover
$(A \cap K)$-definable. In particular, $\Th(K)$ is $0$-h-minimal.
\end{lem}

\begin{proof}
Using (T1), we first find an $A$-definable set $C$. We consider it as a member of an $(A\cap K)$-definable family of sets $C_\xi$, parametrized by $\xi \in \RV^k$.
By Lemma~\ref{lem:tame2finite}, the union $C' := \bigcup_{\xi \in \RV^k} C$ is still finite. Moreover, it is $(A\cap K)$-definable, and since it contains $C$, it satisfies the requirements of (T1).

For the in-particular part, apply this improved (T1) to the characteristic function of an $A$-definable set $X \subset K$ (where $A \subset K \cup \RV$), to find a finite $(A \cap K)$-definable $C$ $1$-preparing $X$.
\end{proof}

Note that we can now use the results from Section~\ref{sec:Prepfam} with $\ell = 0$; we will in particular use Corollary~\ref{cor:prep} several times to prepare $\RV$-parametrized families of subsets of $K$.

The next lemma is a variant of Lemma~\ref{lem:gammaLin} with different assumptions.

\begin{lem}\label{lem:gammaLinTame}
Given an $A$-definable function $f\colon K \to K$, for $A \subset K \cup \RV$, we can find \(C\) which is  \((A\cap K)\)-definable, such that Condition (T1) of Theorem~\ref{thm:tame2vf} holds, and such that moreover, for $B \subset K$ $1$-next to $C$, for $\mu_B$ as in (T1) and for $B' \subset B$ an open ball, the image $f(B')$ is either a singleton (when $\mu_B = 0$) or an open ball of radius $\radop(B') \cdot \mu_B$.
\end{lem}

\begin{proof}
Let $C$ be an $(A\cap K)$-definable set satisfying (T1).
 By Lemma~\ref{lem:InextFam}, the family of balls $B$ $1$-next to $C$ can be parametrized using $\RV$-parameters, so using
Corollary~\ref{cor:prep}, we can prepare the family of $f(B)$ using a finite $(A\cap K)$-definable set $D$.
Similarly, we find a finite $(A\cap K)$-definable set $C_0$ preparing the family $f^{-1}(B_1)$, where $B_1$ runs over the balls $1$-next to $D$.
We claim that the set $C' := C \cup C_0$ does the job.

So let $B' \subset K$ be $1$-next to $C'$, let $B$ be the ball $1$-next to $C$ containing $B'$,
and let $\mu_B \in \Gamma_K$ be as in (T1), i.e., such that we have
\begin{equation}\label{eq:valeq}
|f(x_1) - f(x_2)| = \mu_B \cdot |x_1 - x_2|
\end{equation}
for all $x_1, x_2 \in B$.
We suppose that $\mu_B \ne 0$ (otherwise $f(B')$ is a singleton) and we have to show that for every open ball $B'' \subset B'$,
$f(B'')$ is an open ball of radius $\radop(B'')\cdot \mu_B$.

Let $B''_1$ be the smallest ball containing $f(B'')$.
Using (\ref{eq:valeq}), one obtains that $B''_1$ is an open ball with radius $\radop B''_1 = \radop(B'')\cdot \mu_B$, so it remains to show that $f(B'')$ is equal to the entire ball $B''_1$.

By our definition of $C_0$, $f(B')$ (and hence also $B''_1$) is contained in a ball $B_1$ that is $1$-next to $D$, and by definition of $D$
(and since $B_1$ is certainly not disjoint from $f(B)$), $B_1$ is contained in $f(B)$. In particular, the entire ball $B''_1$ is contained in $f(B)$. However, using (\ref{eq:valeq}) again, we obtain that no element of $B \setminus B''$ can be sent
into $B''_1$, so we deduce $f(B'') = B''_1$, as desired.
\end{proof}

\begin{proof}[Proof of Theorem~\ref{thm:tame2vf}, ``$\Leftarrow$'']
We are given $\lambda \le 1$ in $\Gamma^\times_K$ and an $(A \cup \{\xi\})$-definable set
$X \subset K$, with $A \subset K\cup\RV$ and $\xi \in \RV_\lambda$. We need to find a finite $(A \cap K)$-definable set $C \subset K$ which $\lambda$-prepares $X$. We may assume that $\lambda$ is $A$-definable.

By Lemma~\ref{lem:tame2finite} (2), it is enough to find a $C$ which is $A$-definable. Indeed, then we can get rid of the $\RV$-parameters in the same way as in Lemma~\ref{lem:tame-0h}.

We consider $X$ as a member of an $A$-definable family $X_y$, where $y$ runs over $K$, with
$X = X_y$ for every $y \in Y := \rv_\lambda^{-1}(\xi)$.
(This uses $\lambda \in \dcl_{\Gamma_K}(A)$.)

By Remark~\ref{rem:nocontrol}, we find, for each $y \in K$, a finite $(A \cup \{y\})$-definable set $\hat C_y$
that $1$-prepares $X_y$. Using Lemma \ref{lem:average}, we find an $(A \cup \{y\})$-definable
injective map $h_y\colon \hat C_y \to \RV^k$.
By compactness, we may assume that those $\hat C_y$ and $h_y$ form $A$-definable
families parametrized by $y \in K$. This allows us to write the set $\bigcup_{y \in K} \{y\} \times \hat C_y \subset K^2$
as a disjoint union of graphs of functions $g_\eta\colon Y_\eta \subset K \to K$, which
form an $A$-definable family parameterized by $\eta \in \RV^k$, namely:
$g_\eta(y) := h_y^{-1}(\eta)$, where the domain $Y_\eta$ of $g_\eta$ is the set of those
$y \in K$ for which $\eta$ is in the image of $h_y$. (For some $\eta$, $Y_\eta$ may be empty.)

For each $\eta$, we find a finite $(A\cup \{\eta\})$-definable set $D_\eta \subset K$ $1$-preparing $Y_\eta$ and $Z_\eta := \{y \in Y_\eta \mid g_\eta(y) \in X_y\}$
and also preparing $g_\eta$ in the sense of Lemma~\ref{lem:gammaLinTame} (by applying the lemma to $g_\eta$ extended by $0$ outside of $Y_\eta$).
Using compactness again, we suppose that $D_\eta$ is an $A$-definable family parametrized by $\eta$, so that the union $D := \{0\} \cup \bigcup_\eta D_\eta$ is $A$-definable. Note that by Lemma~\ref{lem:tame2finite} (2), $D$ is finite.

If $D \cap Y \ne \emptyset$, we obtain a finite $A$-definable set
$\bigcup_{y \in D } \hat C_y$
 which $\lambda$-prepares (even $1$-prepares) $X$ and hence we are done. So now assume that $D \cap Y = \emptyset$.
In that case, we claim that the following set $C$ $\lambda$-prepares $X$:
By Lemma~\ref{lem:InextFam}, the balls $1$-next to $D$ form an $A$-definable family $(B_\xi)_\xi$, where $\xi$ runs over $\RV^k$ for some $k \ge 1$.
We let (using Corollary~\ref{cor:prep}) $C \subset K$ be a finite $A$-definable set which $1$-prepares $g_\eta(B_\xi)$ for each $(\eta, \xi)$ satisfying $B_\xi \subset Y_\eta$.

To prove that this $C$ indeed $\lambda$-prepares $X$, we need to verify: For every $x \in X$ and every $x' \in K \setminus X$, there exists a $c \in C$ such that $|x - c|\cdot \lambda \le |x - x'|$.
Let $B_1 \subset K$ be the smallest (closed) ball containing $x$ and $x'$.

Fix a $y_0 \in Y$. Since $\hat C_{y_0}$ $1$-prepares $X$, $\hat C_{y_0} \cap B_1$ is non-empty, so there exists an $\eta$ such that $y_0 \in Y_\eta$ and $g_\eta(y_0) \in B_1$. We fix such an $\eta$ for the remainder of the proof.

Recall that $Y$ is a ball satisfying $Y \cap D = \emptyset$,
and let $B_0$ the ball $1$-next to $D$ containing $Y$.
Our choice of $D$ in particular implies that $g_\eta$ is defined on all of $B_0$ and that for $\mu_{B_0}$ as in Lemma~\ref{lem:gammaLinTame}, $g_\eta(Y)$ and $g_\eta(B_0)$ are open balls satisfying
\begin{equation}\label{eq:th1}
\radop(g_\eta(Y)) = \radop(Y) \cdot \mu_{B_0}
\quad \text{and} \quad
\radop(g_\eta(B_0)) = \radop(B_0) \cdot \mu_{B_0}.
\end{equation}
Moreover, since $D$ also $1$-prepares $\{y \in Y_\eta \mid g_\eta(y) \in X_y\}$, and since $X = X_y$ for all $y \in Y$, we obtain that $g_\eta(Y)$ is either contained in $X$ or disjoint from $X$. By our choice of $\eta$, we have $g_\eta(Y) \cap B_1 \ne \emptyset$. However, $B_1$ is neither contained in $X$ nor disjoint from $X$, so we deduce $B_1 \not\subset g_\eta(Y)$. This implies
\begin{equation}\label{eq:th2}
\radop(g_\eta(Y)) \le \radcl(B_1) = |x - x'|.
\end{equation}
Since $C$ $1$-prepares $g_\eta(B_0)$, there exists a $c \in C$ such that
\begin{equation}\label{eq:th3}
|g_\eta(y_0) - c| \le \radop(g_\eta(B_0)).
\end{equation}
Finally, recall that $Y$ is a fiber of $\rv_\lambda$.
Since $0 \notin B_0$ (which we ensured by putting $0$ into $D$), we deduce $\lambda\cdot \radop(B_0) \le \radop(Y)$.
Putting this together with (\ref{eq:th1}), (\ref{eq:th2}) and (\ref{eq:th3}), we obtain:
\[
|g_\eta(y_0) - c| \cdot \lambda \le \mu_{B_0}\cdot \radop(B_0)\cdot \lambda \le  \mu_{B_0}\cdot \radop(Y) \le  |x - x'|.
\]
Now recall that \(g_\eta(y_0) \in B_1\), so we obtain the desired result:
\[
|x -c| \cdot \lambda \leq \max \{\underbrace{|x-g_\eta(y_0)|}_{\le |x - x'|},|g_\eta(y_0) - c|\}\cdot\lambda \leq |x - x'|.
\qedhere\]
\end{proof}

\subsection{A criterion for preparability}\label{sec:trans}

We already proved that various kinds of objects ``living over $K$'' can be prepared by finite subsets of $K$ in various different senses.
We finish Section~\ref{sec:first-properties} with a criterion simplifying such proofs (Lemma~\ref{lem:fullyT}). Since we want to apply this to quite different
notions of ``being prepared'', we do this at an abstract level, using the following definition.

\begin{defn}[Preparing bad balls]
Let $\cB$ be a set of closed balls in $K$ (the ``bad balls''). We say that a finite set $C \subset K$ \emph{prepares} $\cB$
if for every $B \in \cB$, the intersection $C \cap B$ is non-empty.
\end{defn}

\begin{example}
Given a definable subset $X \subset K$, we can let $\cB$ be the set of all closed balls $B$ that are neither disjoint from $X$ nor contained in $X$.
Then a finite set $C \subset K$ $1$-prepares $X$ if and only if it prepares $\cB$.
Indeed, the implication ``$\Rightarrow$'' is clear. For the other implication, suppose that $B'$ is a ball
$1$-next to $C$ containing both, a point $x \in X$ and a point $x' \in K \setminus X$; then the smallest (closed) ball containing $x$ and $x'$ is disjoint from $C$ but lies in $\cB$.
\end{example}

Note that we can (and will) assume $\cB$ to be an imaginary definable set.
As such, in the above example, $\cB$ is definable over the same parameters as $X$.

By a compactness argument (as in the proof of Lemma~\ref{lem:type-0-h-min}), a $\emptyset$-definable $\cB$ can be prepared by a finite $\emptyset$-definable $C \subset K$ if and only if, in a sufficiently saturated model $K$, no ball $B \in \cB$ is disjoint from $\acl_K(\emptyset)$. The following lemma provides an even weaker condition that needs to be checked if one wants to prove that $\cB$ can be prepared:

\begin{lem}[Criterion for preparability]\label{lem:fullyT}
Suppose that $\Th(K)$ is $0$-h-minimal and that $K$ is $|\cL|^+$-saturated.
Let $\cB$ be a $\emptyset$-definable set of closed balls in $K$. We call an arbitrary ball $B' \subset K$ a ``bad ball'' if it contains a ball $B \in \cB$ as a subset.

Suppose that for every $a \in K$, there is no open ball $B' \subset K \setminus \acl_K(a)$ which is $1$-next to $a$ and bad.
Then there exists a finite $\emptyset$-definable set $C \subset K$ preparing $\cB$.
\end{lem}

\begin{remark}
It might sound more natural to let the balls in $\cB$ be open instead of closed, but that would make the lemma false: We use, in the proof, that any open bad ball already contains a closed bad subball.
\end{remark}

\begin{proof}[Proof of Lemma~\ref{lem:fullyT}]
Suppose that no finite $\emptyset$-definable set $C$ prepares $\cB$.
Then, by saturation, we find a single ball $B \in \cB$ which is disjoint from all finite $\emptyset$-definable $C$; in other words, this ball satisfies $B \cap \acl_K(\emptyset) = \emptyset$. We fix $B$ for the remainder of the proof.

If $B$ is $1$-next to some $a \in \acl_K(\emptyset)$, then we also have $B \cap \acl_K(a) = \emptyset$, and we get a contradiction to the assumption that such a $B$ should not be bad. Therefore, if for $a \in \acl_K(\emptyset)$ we denote by $B_a$ the ball $1$-next to $a$ containing $B$, none of the sets $B_a \setminus B$ is empty. By saturation, also the intersection
$\bigcap_{a \in \acl_K(\emptyset)} B_a \setminus B$ is non-empty.
Let $B'$ be the smallest ball containing $B$ and any chosen element of that intersection. This ball has the following properties: It is closed, it is disjoint from $\acl_K(\emptyset)$, and it is strictly bigger than $B$.
Set $\gamma := \radcl(B)$ and $\mu := \radcl(B')$.

All elements of $B'$ have the same type over $\RV$; indeed, if two elements of $B'$ could be distinguished by a formula with parameters from $\RV$, then any finite set $C$ preparing the $\RV$-parametrized family of subsets of $K$ defined by that formula would have to contain points of $B'$; but then $C$ cannot be $\emptyset$-definable, since $B' \cap \acl_K(\emptyset) = \emptyset$.

We deduce that all the open balls $B_{<\mu}(b) \subset B'$ (for $b \in B'$) are bad. Indeed, for $b \in B$, we have
$B_{<\mu}(b) \supset B_{\le \gamma}(b) \in \cB$, and since any other $b' \in B'$ has the same type as $b$ over $\RV$ (and $\cB$ is $\emptyset$-definable), we also have $B_{\le \gamma}(b') \in \cB$,
witnessing that $B_{<\mu}(b')$ is bad.

Now fix an arbitrary $a \in B'$. By saturation, there exists a $b' \in B'$ such that $B'' := B_{<\mu}(b')$ is disjoint from $\acl_K(a)$. This however contradicts the assumption of the lemma, since $B''$ is bad and $1$-next to $a$.
\end{proof}

\section{Derivation, the Jacobian Property and Taylor approximation}\label{sec:derJac}

We continue assuming that $K$ is a valued field of equi-characteristic $0$,
considered as a structure in a fixed language $\cL \supset \Lval$, and we
assume that $\Th(K)$ is $1$-h-minimal.
The goal of this section is to prove various preparation results for definable functions $f\colon K \to K$,
starting with Theorem~\ref{thm:der} which states that $f$ is almost everywhere differentiable, and
culminating in Theorem~\ref{thm:high-ord}, which states that away from a finite set, $f$
has good approximations by its Taylor polynomials.

\subsection{Derivation, strict and classical}
\label{sec:derivative}

We define the derivative of functions $f\colon K \to K$ as usual:

\begin{defn}[Classical derivative]
We say that the \emph{(classical) derivative} of a function $f\colon K \to K$ exists at $x\in K$ if there exists $a\in K$ such that for each $\varepsilon$ in $\Gamma^\times_K$ there exists $\delta$ in $\Gamma^\times_K$ such that for all $y\in K$ with $|x-y|<\delta$ and $x\ne y$ one has
$$
|\frac{f(x)-f(y)}{x-y} - a | < \varepsilon;
$$
we then write $f'(x)$ for $a$.
\end{defn}

In the context of totally disconnected fields, one sometimes needs to put a stronger condition on the existence of derivatives, namely:

\begin{defn}[Strict derivative]\label{defn:strict-der}
We say that the \emph{strict derivative} of a function $f\colon K \to K$ exists at $x\in K$ if there exists $a\in K$ such that for each
$\varepsilon$ in $\Gamma^\times_K$ there exists $\delta$ in $\Gamma^\times_K$ such that for all $y_i\in K$
with $|x-y_i|<\delta$ for $i=1,2$ and $y_1 \not=y_2$ one has
$$
|\frac{f(y_1)-f(y_2)}{y_1-y_2} - a | < \varepsilon.
$$
\end{defn}

\begin{remark}
As usual,
one easily verifies that the set of $x$ where the classical derivative of a definable function $f$ exists is definable over the same parameters as $f$, and similarly for the strict derivative. Moreover, the derivative of $f$ is definable over the same parameters as $f$.
\end{remark}

\begin{thm}[Existence of derivatives]\label{thm:der}
Let $K$ be a valued field of equi-characteristic $0$ such that $\Th(K)$ is $1$-h-minimal, and let $f\colon K \to K$ be a definable (with parameters) function.
Then the strict derivative of $f$ exists almost everywhere, i.e., the set of $x \in K$ such that the strict derivative of $f$ does not exist at $x$ is finite.
\end{thm}

Most of the remainder of this subsection is devoted to the proof of this theorem, so from now on we fix a definable function $f\colon K \to K$.
By Lemma~\ref{lem:addconst}, we may as well impose $f$ to be $\emptyset$-definable.

We also fix a handy notation for difference quotients: Given unequal $x_1,x_2$ in $K$, we set
\[
q_f(x_1, x_2) := \frac{f(x_1) -f(x_2)}{x_1 - x_2}.
\]
For convenience we put $q_f(x, x):=0$ for any $x\in K$.

We start by proving an auxiliary lemma, which is a weak form of Theorem~\ref{thm:der}, stating that in some sense, one has finitely-valued derivatives.

\begin{lem}\label{lem:f'fin}
For $x \in K$, let $A_x$ be the set of accumulation points of
\begin{equation}\label{eq:def_A_x}
\lim_{\substack{y_1,y_2\to x\\y_1\ne y_2}}q_f(y_1, y_2).
\end{equation}
There is a finite $\emptyset$-definable set $C\subset K$ such that, for every $x\in K\setminus C$, $A_x$ is finite and we moreover have
\begin{equation}\label{eq:limmin}
\lim_{\substack{y_1,y_2\to x\\y_1\ne y_2}}
\min_{a\in A_x} \left|  q_f(y_1, y_2) - a  \right| = 0.
\end{equation}
\end{lem}

\begin{proof}
Clearly, the set $C$ of $x$ where $A_x$ is not as desired is $\emptyset$-definable (using Lemma~\ref{lem:finite} to express finiteness of $A_x$). We need to show that $C$ is finite.

We may assume that $K$ is sufficiently saturated and we suppose that $C$ is infinite. Then $C$ contains a transcendental element $x_0$ (i.e., $x_0 \notin \acl_K(\emptyset)$).
We will prove that for transcendental $x_0 \in K$,
there exists a finite set $A$ satisfying (\ref{eq:limmin}).
One easily checks that this implies that $A_{x_0} \subset A$ and that then, $A_{x_0}$ also satisfies (\ref{eq:limmin}), contradicting $x_0 \in C$. Constructing $A$ is done in several steps.

\medskip

Step 1:
Using (once more) that $K$ is sufficiently saturated, we find an entire open ball $B_1 := B_{<\mu}(x_0)$ that is disjoint from $\acl_K(\emptyset)$.
We fix these $\mu$ and $B_1$ once and for all and more generally define $B_\lambda :=  B_{<\lambda\cdot\mu}(x_0)$, for $\lambda \le 1$.

\medskip

For Steps 2--4, we fix $\lambda \le 1$ in $\Gamma^\times_K$, and
$\zeta$ will always be an element of $\RV$ satisfying $|\zeta| < \lambda \cdot \mu$ (so that $\zeta = \rv(x - x')$ for some $x, x' \in B_\lambda$).

\medskip

Step 2: For $\zeta$ as above, the subset of $\RV_\lambda$ defined by $Q_{x,\zeta} := \{\rv_\lambda(q_f(x, x')) \mid x' \in x+ \rv^{-1}(\zeta)\}$ is independent of $x$ when $x$ runs over $B_\lambda$.

\medskip

Proof: We $\lambda$-prepare the family $(Q_{x,\zeta})_{x,\zeta}$ using Proposition~\ref{prop:uniform} (i.e., we consider $Q_{x,\zeta}$ as a fiber of a definable set $Q \subset  K \times \RV \times \RV_\lambda$ and $\lambda$-prepare $Q$). The finite $\emptyset$-definable
set $C'$ obtained in this way is disjoint from $B_1$ (since $B_1 \cap \acl_K(\emptyset) = \emptyset$), so $B_\lambda$ is contained in a ball $\lambda$-next to $C'$. This implies Step 2.
\qed2

\medskip

Set $Q_\zeta := Q_{x,\zeta}$.

\medskip

Step 3: $Q_\zeta \subset Q_{2\zeta}$ (where $2\zeta$ is $\rv(2)\cdot\zeta$).

\medskip

Proof: Fix any $\xi\in Q_\zeta$. (We need to show that $\xi \in  Q_{2\zeta}$.)
Choose $x_1 \in B_\lambda$ witnessing $\xi \in Q_{x_0,\zeta}$, i.e., such that $\rv(x_0 - x_1) = \zeta$ and
$\rv_\lambda(q_f(x_0, x_1)) = \xi$; similarly choose $x_2 \in B_\lambda$ such that $\rv(x_1 - x_2) = \zeta$ and $\rv_\lambda(q_f(x_1, x_2)) = \xi$. Then the following computation shows that $\rv_\lambda(q_f(x_0, x_2)) = \xi$ (which implies $\xi \in Q_{2\zeta}$):
Set $r_i := x_{i-1} - x_i$ and $s_i := f(x_{i-1}) - f(x_i)$ for $i = 1,2$. We have
\begin{align*}
|q_f(x_0,x_2) - q_f(x_0, x_1)| &=
\left|\frac{s_1 + s_2}{r_1 + r_2} - \frac{s_1}{r_1}\right| =
\left|\left(\frac{s_2}{r_2} - \frac{s_1}{r_1}\right)\cdot\frac{r_2}{r_1+r_2}\right|\\
&= \underbrace{\left|q_f(x_1,x_2) - q_f(x_0,x_1)\right|}_{< \lambda\cdot|q_f(x_0,x_1)|}\cdot\underbrace{\left|\frac{r_2}{r_1+r_2}\right|}_{=1}
\end{align*}
and hence $\rv_\lambda(q_f(x_0,x_2)) = \rv_\lambda(q_f(x_0,x_1))$.
\qed3

\medskip

Applying Step 3 repeatedly shows: $Q_\zeta \subset Q_{2^n\zeta}$ for every integer $n  \ge 0$.

\medskip

Step 4: $Q_\zeta$ is a singleton for every $\zeta$.

\medskip

Proof: Suppose otherwise, i.e., $\xi_1, \xi_2 \in Q_\zeta$ with $\xi_1 \ne \xi_2$. The set
$$
X := \{x \in B_\lambda \mid \rv_\lambda(q_f(x_0, x)) = \xi_1\};
$$
is definable, and hence it can be
$1$-prepared by a finite set $D$. (We do not care about the parameters needed to define $X$ and $D$.)
However,
for each of the (disjoint) balls $B_n := \{x\in B_\lambda\mid \rv (x_0 - x) = 2^n \zeta\}$ (where $n$ runs over the non-negative integers), we have neither $X \cap B_n = \emptyset$ (since $\xi_1 \in Q_{2^n\zeta}$) nor $B_n \subset X$ (since $\xi_2 \in Q_{2^n\zeta}$); hence for $D$ to $1$-prepare $X$, we would need $D \cap B_n \ne \emptyset$ for every $n$, contradicting the finiteness of $D$.
\qed4

\medskip

Let us reformulate what we obtained until now in a slightly different way: Given any $\zeta \in \RV$ and any $\lambda \in \Gamma_K$ satisfying $|\zeta|/\mu < \lambda \le 1$, Step~4 states that the entire set
$\tilde{Q}_{\zeta,\lambda} := \{q_f(x, x') \mid x, x' \in B_\lambda, \rv(x - x') = \zeta\}$
is contained in a single ball $\tilde B_{\zeta,\lambda}$ $\lambda$-next to $0$.

\medskip

Step 5: Using Corollary~\ref{cor:prep}, choose a finite set $A \subset K$ $1$-preparing the family $(\tilde{Q}_{\zeta,\lambda})_{\zeta,\lambda}$, for $\zeta \in \RV$ and \(\lambda \in \Gamma_K^\times\) satisfying $|\zeta|/\mu < \lambda \leq 1$. (Again, we do not care about parameters.)

\medskip

The last step consists in showing that the set $A$ satisfies (\ref{eq:limmin}), as desired. More precisely, we show:

\medskip

Step 6: There exists a constant $\kappa \in \Gamma_K$ such that for every $\lambda < 1$ and every $x, x' \in B_\lambda$ with $x \ne x'$, we have $\min_{a\in A} |  q_f(x,x') - a  | \le \lambda \cdot \kappa$.

\medskip

Proof: The constant is $\kappa := \max \{|a| \mid a \in A\}$. Let $x, x' \in B_\lambda$ be given. For $\zeta := \rv(x - x')$, we have $q_f(x,x') \in \tilde{Q}_{\zeta,\lambda} \subset \tilde B_{\zeta,\lambda}$. Note that $\tilde B_{\zeta,\lambda}$ is an open ball of radius $\lambda\cdot |q_f(x, x')|$.
Since $A$ $1$-prepares $\tilde{Q}_{\zeta,\lambda}$, there exists an $a \in A$ such that $|q_f(x,x') - a| \le \radop(\tilde B_{\zeta,\lambda}) = \lambda\cdot |q_f(x, x')|$. Since $\lambda < 1$, this in particular implies $|q_f(x,x')| = |a|$, so the right hand side is $\lambda \cdot |a| \le \lambda \cdot \kappa$ and we are done.
\qed6

\medskip

This finishes the proof of the lemma.
\end{proof}

Using the lemma, we can now prove that the derivative of $f$ exists almost everywhere.

\begin{proof}[Proof of Theorem~\ref{thm:der}]
Let $A_x$ (for $x \in K$) and $C$ be as in Lemma \ref{lem:f'fin}.
We need to show that, for almost all $x\in K \setminus C$, the set $A_x$ is a singleton.
Suppose otherwise, i.e., that $A_x$ is not a singleton for infinitely many $x$. Then as usual, we can find a ball $B$ such that $A_x$ consists of several elements for every $x \in B$.

We first shrink $B$ in such a way that the cardinality $\#A_x$ is constant for $x \in B$, and then we further shrink it to make
the map $x \mapsto A_x$ ``approximately constant'' on $B$ in the following sense: There is
a $\mu \in \Gamma^\times_K$ such that for every $x, x' \in B$, the relation $a \sim_\mu a' :\iff |a - a'| < \mu$ defines a bijection between $A_x$ and $A_{x'}$.
This shrinking of $B$ is possible as follows: By $1$-preparing the (graph of the) function $B \to \Gamma^\times_K, x \mapsto \min_{a_1, a_2 \in A_{x}, a_1 \ne a_2} |a_1 - a_2|$, we find a subball (which we will now call $B$) on which this minimum is constant equal to some $\mu \in \Gamma^\times_K$. We then choose any $x_0 \in B$, and using the definition (\ref{eq:def_A_x}) of $A_{x_0}$, we replace $B$ by an even smaller ball around $x_0$
such a way that for every $x_1, x_2 \in B$, we have
\[
\min_{a \in A_{x_0}}|\frac{f(x_1) - f(x_2)}{x_1 - x_2} - a| < \mu.
\]
This implies that $\sim_\mu$ defines bijections as desired.

Fix any $a_0 \in \bigcup_{x \in B}A_x$. We
apply Lemma~\ref{lem:gammaLin} (2) (the Valuative Jacobian Property) to the $\cL(a_0)$-definable function $\tilde f(x) := f(x) - a_0x$. This allows
us to further shrink $B$ in such a way that the quotient $|\tilde f(x_1) - \tilde f(x_2)|/|x_1 - x_2|$ is constant for $x_1, x_2 \in B, x_1 \ne x_2$.

This now leads to a contradiction, as follows.
Fix any $x \in B$. For every $a \in A_x$, there exist $x_1, x_2 \in B$ with $|\frac{f(x_1) - f(x_2)}{x_1 - x_2} - a| < \mu$ (by definition of $A_x$). Since
\[
\frac{\tilde f(x_1) - \tilde f(x_2)}{x_1 - x_2} =
\frac{f(x_1) - f(x_2)}{x_1 - x_2} - a_0,
\]
this implies $|\frac{\tilde f(x_1) - \tilde f(x_2)}{x_1 - x_2}| < \mu$ if $a \sim_\mu a_0$ and
$|\frac{\tilde f(x_1) - \tilde f(x_2)}{x_1 - x_2}| \ge \mu$ if $a \not\sim_\mu a_0$.
Since $A_x$ contains elements $a$ of both kinds (with and without $a \sim_\mu a_0$), this contradicts $|\frac{\tilde f(x_1) - \tilde f(x_2)}{x_1 - x_2}|$ being constant.
\end{proof}

Using the existence of derivatives, we can now reformulate the Valuative Jacobian Property (Lemma~\ref{lem:gammaLin}) in a nicer way:

\begin{cor}[Valuative Jacobian Property]\label{cor:gammaLin}
For every $\emptyset$-definable function $f\colon K \to K$, there exists a finite $\emptyset$-definable set $C \subset K$
such that for every ball $B$ $1$-next to $C$, $f'$ exists on $B$, $|f'|$ is constant on $B$, and we have the following:
\begin{enumerate}
 \item For every $x_1, x_2 \in B$, we have $|f(x_1) - f(x_2)| = |f'(x_1)|\cdot |x_1 - x_2|$.
   (In particular, $f$ is constant on $B$ if $f' = 0$ on $B$.)
 \item
   If $f' \ne 0$ on $B$, then for every open ball $B' \subset B$, $f(B')$ is an open ball of radius $|f'(x)| \cdot\radop(B')$
   for any $x \in B$.
\end{enumerate}
\end{cor}

\begin{proof}
The only difference between this and Lemma~\ref{lem:gammaLin} is that the factor called $\mu_B$ in Lemma~\ref{lem:gammaLin} is now claimed
to be equal to $|f'(x)|$ for any $x \in B$. But indeed, by definition of the derivative, for $x \in B$ and $x'$ sufficiently close to $x$, we have
$|f'(x)| = \frac{|f(x') - f(x)|}{|x' - x|} = \mu_B$.
\end{proof}

\subsection{The Jacobian Property and Taylor approximations}
\label{sec:Jac}

We now come to one of the central results of this paper: Every definable function $K \to K$ is, away from a finite set $C$,
well approximated by its Taylor polynomials. Here follows the precise statement and various variants and corollaries. The proof will be given in the next subsection, and higher dimensional variants will be deduced in Subsections~\ref{sec:sjp} and \ref{sec:taylor-box}.

\begin{defn}[Taylor Polynomial]\label{defn:taylor}
Let $f\colon X \subset K \to K$ be a function which is $r$-fold differentiable at $x_0$, for some $r \in \NN$ and some $x_0 \in K$.
Then we write
\begin{equation}
T^{<r+1}_{f,x_0}(x) := T^{\le r}_{f,x_0}(x) := \sum_{i = 0}^r \frac{f^{(i)}(x_0)}{i!}(x - x_0)^i
\end{equation}
for the Taylor polynomial of degree $r$ of $f$ around $x_0$. (Here, $f^{(i)}$ denotes the $i$-th derivative of $f$.) Similarly, when $f:X\subset K^m\to K$ is $r$-fold differentiable at $x_0\in X$, then we write  \begin{equation}\label{eq:taylor-m}
T^{<r+1}_{f,x_0}(x):=T^{\le r}_{f,x_0}(x) := \sum_{i \in\NN^m,\ |i|\le r} \frac{f^{(i)}(x_0)}{i!}(x - x_0)^i
\end{equation}
for the Taylor polynomial of degree $r$ of $f$ around $x_0$ (where (\ref{eq:taylor-m}) uses multi-index notation, with $|i|=\sum_{j=1}^m i_j$).
\end{defn}

\begin{thm}[Taylor approximations]\label{thm:high-ord}
Let $f\colon K \to K$ be a $\emptyset$-definable function and let $r \in \NN$ be given. Then there exists a finite $\emptyset$-definable set $C$ such that for every open ball $B$ $1$-next to $C$, $f$ is $(r+1)$-fold differentiable on $B$, $|f^{(r+1)}|$ is constant on $B$, and we have:
\begin{equation}\label{eq:t-higher}
 |f(x) -  T^{\le r}_{f,x_0}(x) | \leq |f^{(r+1)}(x_0)\cdot (x-x_0)^{r+1}|
\end{equation}
for every $x_0, x \in B$.
\end{thm}

\begin{remark}
As explained in Remark~\ref{rem:acl}, here (in Theorem~\ref{thm:high-ord}) and in the following (Corollaries~\ref{cor:high-ord}, \ref{cor:high-ord-S}, \ref{cor:JP}), we can as well allow $f$ to be $(A \cup \RV)$-definable for $A \subset K$ and obtain an $A$-definable set $C$.
\end{remark}
\begin{remark}\label{rem:high-ord}
In fact, Theorem \ref{thm:high-ord} also comes in a variant with equality instead of weak inequality in (\ref{eq:t-higher}). Indeed, this variant follows from Corollary \ref{cor:high-ord-S} for $r+1$ and the non-archimedean triangular property.
\end{remark}

As in the reals, we can also express the error in terms of the $r$-th derivative. We mention two different
such variants of the theorem. (Both will also play a role in the proof of the theorem.)

\begin{cor}[Taylor approximations]\label{cor:high-ord}
Let $f\colon K \to K$ be a $\emptyset$-definable function and let $r \in \NN$ be given. Then there exists a finite $\emptyset$-definable set $C$ such that for every open ball $B$ $1$-next to $C$, $f$ is $r$-fold differentiable on $B$, $|f^{(r)}|$ is constant on $B$, and we have:
\begin{equation}\label{eq:c-higher}
 |f(x) -  T^{\le r}_{f,x_0}(x) | \leq |f^{(r)}(x_0)\cdot (x-x_0)^{r+1}|/\radop(B)
\end{equation}
for every $x_0, x \in B$.
\end{cor}

\begin{proof}
Apply (a) Theorem~\ref{thm:high-ord} to $f$, (b) Corollary~\ref{cor:gammaLin} to $f^{(r)}$ and (c) Corollary~\ref{cor:prep} to $\rv\circ f^{(r)}$,
and let $C$ be the union of the sets obtained in these ways.
Then we have, for $B$ $1$-next to $C$ and for $x_0 \in B$,
\[
|f^{(r+1)}(x_0)| \cdot \radop(B) \overset{(b)}{=} \radop(f^{(r)}(B)) \overset{(c)}{\le} |f^{(r)}(x_0)|.
\]
Using this, (\ref{eq:t-higher}) implies (\ref{eq:c-higher}).
\end{proof}

\begin{cor}[Taylor approximations]\label{cor:high-ord-S}
Let $f\colon K \to K$ be a $\emptyset$-definable function and let $r \in \NN$ be given. Then there exists a finite $\emptyset$-definable set $C$ such that for every open ball $B$ $1$-next to $C$, $f$ is $r$-fold differentiable on $B$, $|f^{(r)}|$ is constant on $B$, and we either have
\begin{equation}\label{eq:c-higher-S}
 |f(x) -  T^{\le r}_{f,x_0}(x) | < |f^{(r)}(x_0)\cdot (x-x_0)^{r}|
\end{equation}
or $f(x) =  T^{\le r}_{f,x_0}(x)$ (if $|f^{(r)}(x_0)| = 0$)
for every $x_0, x \in B$.
\end{cor}

\begin{proof}
This follows directly from Corollary~\ref{cor:high-ord}, using that $|x - x_0| < \radop(B)$.
\end{proof}

Note that in the case $r = 1$, we in particular get back the $\RV_\lambda$-Jacobian Property, studied in an analytic setting in \cite[Theorem 6.3.7, Remark 6.3.16]{CLip}.
More precisely, the following version of the Jacobian Property is even uniform for all $\lambda \le 1$.

\begin{cor}[Jacobian Property]\label{cor:JP}
Let $f\colon K \to K$ be a $\emptyset$-definable function. Then there exists a finite $\emptyset$-definable set $C$ such that for every $\lambda \le 1$ in $\Gamma^\times_K$, for every ball $B$ $\lambda$-next to $C$ and for every $x_0$ and $x$ in $B$ with $x\not=x_0$, we have:
\begin{enumerate}
 \item The derivative $f'$ exists on $B$ and $\rv_\lambda\circ f'$ is constant on $B$.
 \item $\rv_\lambda(\frac{f(x) - f(x_0)}{x - x_0}) = \rv_\lambda(f')$
 \item For every open ball $B' \subset B$, $f(B')$ is either a point or an open ball.
\end{enumerate}
\end{cor}

\begin{defn}[Jacobian Property]\label {defn:JP}
If Conditions (1)--(3) of Corollary~\ref{cor:JP} hold, we say that $f|_B$ has the \emph{$\RV_\lambda$-Jacobian Property} (or just \emph{Jacobian Property}, in the case $\lambda = 1$).
\end{defn}

\begin{proof}[Proof of Corollary~\ref{cor:JP}]

Choose $C$ using Corollary~\ref{cor:high-ord} (applied to $f$, with $r = 1$); this will ensure (2), as we shall see below.

To obtain that $\rv_\lambda\circ f'$ is constant on balls $\lambda$-next to $C$, enlarge $C$ using Proposition~\ref{prop:uniform}: Consider the graph of $\rv_\lambda \circ f'$ as a subset of $K \times \RV_\lambda \subset K \times \RV_\bullet$ and take the union of all of those as the set $W$ in the proposition.

Finally, to obtain that $f(B')$ is either a point or a ball, enlarge $C$ once more using Lemma~\ref{lem:open-to-open}.

It remains to verify (2), so fix $\lambda \le 1$, fix a ball $B$ $\lambda$-next to $C$, fix $x_0, x \in B$, and let $B'$ denote the ball $1$-next to $C$ containing $B$. Then we have
\begin{equation}\label{eq:to-jac}
|f(x) -  f(x_0) - f'(x_0)(x - x_0) | \leq |f'(x_0)\cdot (x-x_0)^2|/\radop(B').
\end{equation}
Since $\radop(B) = \lambda\cdot \radop(B')$, we have
$|x - x_0| < \lambda\cdot \radop(B')$,
so after dividing both sides of (\ref{eq:to-jac}) by $|x - x_0|$, we obtain
\[
|\frac{f(x) -  f(x_0)}{x - x_0} - f'(x_0)|< |f'(x_0)|\cdot \lambda,
\]
which is equivalent to (2).
\end{proof}

\subsection{Proof of Taylor approximation}

Before we dive into the proof of Theorem~\ref{thm:high-ord}, here is a somewhat technical lemma that will be needed for one special case.

\begin{lem}[Preparing equivalence relations]\label{lem:equivalence}
Suppose that $\Gamma_K$ is discrete. Fix an open ball $B \subset K$ which is
disjoint from $\acl(\emptyset)$. Suppose that $\sim$ is an $\emptyset$-definable binary relation on $K$ with the following properties:
(i) the restriction of $\sim$ to $B \times B$ is an equivalence relation with finitely many equivalence classes; (ii) every ball $B'$ strictly contained in $B$ (i.e., $B' \subsetneq B$) is contained in a single equivalence class.
Then the entire ball $B$ is a single equivalence class.
\end{lem}

\begin{proof}
We write $k := \cO_K/\cM_K$ for the residue field, $\res\colon \cO_K \to k$ for the residue map, and we
choose a linear bijection $f\colon B \to \cO_K$. (Note that since $\Gamma_K$ is discrete, the open ball $B$ is also closed.) Then $\sim$ induces an equivalence relation $\sim_k$ on $k$
(satisfying $x \sim x' \iff \res(f(x)) \sim_k \res(f(x'))$ for $x, x' \in B$).

For $x \in K$, set $W_x := \{\rv(x' - x) \mid x' \in K, x' \sim x\}$.
Since $B$ is disjoint from $\acl_K(\emptyset)$, there exists a single set $W \subset \RV$ such that $W_x = W$ for all $x \in B$
(by Corollary~\ref{cor:prep}). This set $W$ yields a set $T \subset k$ such that $z \sim_k z'$ if and only if $z - z' \in T$ (for $z, z' \in k$).

A straight-forward computation shows: The fact that ``$z - z' \in T$'' defines an equivalence relation on $k$ implies that $T$ is a $\ZZ$-submodule of $k$. The number of equivalence classes is equal to the cardinality of the $\ZZ$-module quotient $k/T$, so to finish the proof of the lemma\footnote{The authors would like to thank W.~Singhof for providing the argument at the end the proof of \Cref{lem:equivalence}.}, suppose that this quotient has finite cardinality $N > 1$.
Then for any $z \in k$, we have $Nz \in T$. However, this fails if we choose any $z' \in k \setminus T$ and set $z := \frac1Nz'$.
\end{proof}

We now come to the proof of Theorem~\ref{thm:high-ord}. It goes by induction on $r$, where Corollaries~\ref{cor:high-ord} and \ref{cor:high-ord-S} are used as intermediate steps. More precisely, Theorem~\ref{thm:high-ord} follows from the following four lemmas:

\begin{lem}\label{lem:h:t1}
Theorem~\ref{thm:high-ord} holds for $r = 0$.
\end{lem}

\begin{lem}\label{lem:h:ts}
For any $r \ge 1$, if Theorem~\ref{thm:high-ord} holds for $r-1$ (for every $1$-h-minimal theory)
and Corollary~\ref{cor:high-ord-S} holds for all $r' < r$,
then Corollary~\ref{cor:high-ord-S} holds for $r$.
\end{lem}

\begin{lem}\label{lem:h:sc}
For any $r \ge 1$, if Corollary~\ref{cor:high-ord-S} holds for $r$ (for every $1$-h-minimal theory), then Corollary~\ref{cor:high-ord} holds for $r$.
\end{lem}

\begin{lem}\label{lem:h:ct}
For any $r \ge 1$, if Corollary~\ref{cor:high-ord} holds for $r$
(for every $1$-h-minimal theory), then Theorem~\ref{thm:high-ord} holds for $r$.
\end{lem}

For all of the proofs, first note that the condition about sufficient differentiability of $f$ on $K\setminus C$ appearing in the theorem and in the corollaries is easily obtained using Theorem~\ref{thm:der}; this will not be further mentioned. Similarly, $|f^{(r)}|$ (or $|f^{(r+1)}|$) can easily be made constant on balls $1$-next to $C$ using Corollary~\ref{cor:prep}.

\begin{proof}[Proof of Lemma~\ref{lem:h:t1}]
In the case $r = 0$, (\ref{eq:t-higher}) becomes $|f(x) - f(x_0)| \le |f'(x_0)\cdot (x - x_0)|$,
so the theorem follows from Corollary~\ref{cor:gammaLin}.
\end{proof}

The proofs of the other three lemmas are long and technical, but the main idea is simple (and the same in all three lemmas): Given a function $f$ which we want to control on a ball $B$, we define
$g(x) := f(x) - ax^r$ in such a way that the $r$-th derivative of $g$ is small on $B$. In this way, applying the inductive assumption to $g$ yields a particularly strong bound, which will be good enough to obtain the desired bound on $f$.

The difficulty with this approach is that we need to control the parameters over which $g$ is definable. A powerful ingredient for this is
Lemma~\ref{lem:fullyT}, which allows us to use points near $B$ as parameters. Nevertheless, various additional technical tricks are needed (different for each of the lemmas) to make the proofs work.

\begin{proof}[Proof of Lemma~\ref{lem:h:ts}]Namely:\\
From $\underbrace{|f -  T^{\le r-1}_{f,x_0} | \leq |f^{(r)}\cdot (x-x_0)^{r}|}_{(\ref{eq:t-higher})\text{ for }r-1}$  to
$\underbrace{|f -  T^{\le r}_{f,x_0} | < |f^{(r)}\cdot (x-x_0)^{r}|}_{(\ref{eq:c-higher-S})}$.
\medskip

We assume that $K$ is sufficiently saturated.
Then it suffices to prove that on every open ball $B \subset K$ which is disjoint from $\acl(\emptyset)$, either (\ref{eq:c-higher-S}) holds, or we have $f(x) = T^{\leq r}_{f,x_0}(x)$ (by the compactness argument
from Lemma~\ref{lem:type-0-h-min}), so fix such a $B$ and suppose that $x_0, x \in B$ violate (\ref{eq:c-higher-S}).

We may assume that $f^{(r)} \ne 0$ on $B$, since
otherwise, Theorem~\ref{thm:high-ord} for $r-1$ implies $f(x) = T^{\leq r-1}_{f,x_0}(x) = T^{\leq r}_{f,x_0}(x)$ for any $x_0, x \in B$ and we are done.

\medskip

Case 1: There exists an open ball $B' := B_{<\delta}(x_0)$ containing $x$ which is strictly smaller than $B$.

\medskip

We fix the above $\delta$ for the remainder of the proof of Case 1.
Also, fix any $a \in B$ and choose any open ball $B'' = B_{< \delta}(x_0') \subset B$ of radius $\delta$ disjoint from $\acl(a)$. (Such a $B''$ exists
by saturation.) Since $x_0$ and $x_0'$ have the same type over $\Gamma_K$ (by Lemma~\ref{lem:type-0-h-min}), there exists
an $x' \in B''$ such that $x'_0, x'$ violate (\ref{eq:c-higher-S}).

We apply Theorem~\ref{thm:high-ord} for $r-1$ to the function
\[g(x) := f(x) - \frac{f^{(r)}(a)}{r!}\cdot x^r.
\]
Since $g$ is $a$-definable and $B''$ is disjoint from $\acl(a)$, we obtain
\begin{equation}\label{eq:gpt}
|g(x') - T^{\le r-1}_{g,x'_0}(x')| \le  | g^{(r)}(x'_0)(x'-x'_0)^r|.
\end{equation}
The definition of $g$ has been chosen such
that $g^{(r)}(x) = f^{(r)}(x) - f^{(r)}(a)$ for all $x \in B$, so using that
$\rv(f^{(r)})$ is constant on $B$ (this uses Corollary~\ref{cor:prep} and that $B \cap \acl_K(\emptyset) = 0$), we deduce
\begin{equation}\label{eq:g<f}
|g^{(r)}(x)| < |f^{(r)}(x)|
\end{equation}
for all $x \in B$.
From this, we obtain the following (explanation of $(\star)$ below):
\begin{equation}\label{eq:get-chS}
\begin{aligned}
 |f(x') -  T^{\le r}_{f,x'_0}(x') | &\overset{(\star)}{=} |g(x') -  T^{\le r}_{g,x'_0}(x') |\\
 &= |g(x') -  T^{\le r-1}_{g,x'_0}(x') + \frac{1}{r!}g^{(r)}(x_0')(x'-x'_0)^{r}|\\
 &\leq \max\{|g(x') -  T^{\le r-1}_{g,x'_0}(x')|,|\frac{1}{r!}g^{(r)}(x_0')(x'-x'_0)^{r}|\}\\
 & \overset{(\ref{eq:gpt})}{\le}  | g^{(r)}(x'_0)(x'-x'_0)^r| \overset{(\ref{eq:g<f})}{<} |f^{(r)}(x'_0)\cdot (x'-x'_0)^{r}|
\end{aligned}
\end{equation}

$(\star)$: $f$ and $g$ differ by a polynomial of degree $r$, so their Taylor approximations differ by the same polynomial.

Thus (\ref{eq:c-higher-S}) holds for $x_0', x'$, which
is a contradiction and hence finishes the proof of Case~1.

\medskip

Case 2: The only open ball $B' \subset B$ containing both $x_0$ and $x$ is already $B$ itself.

\medskip

Note that this can only happen if $\Gamma_K$ is discrete (otherwise the radius of $B'$ can be taken strictly between
$|x_0 - x|$ and $\radop(B)$), and it also means that $|x_0 - x| = \radcl(B) =: \delta$.

By Case~1, we may assume that (\ref{eq:c-higher-S}) holds on every proper subball of $B$. Moreover,
we apply Corollary~\ref{cor:high-ord-S} inductively to derivatives of $f$ to get, for $1 \le i \le r$:
\begin{equation}\label{eq:der}
  |f^{(i)}(x) -  T^{\le r-i}_{f^{(i)},x_0}(x) | < |f^{(r)}|\cdot |(x-x_0)^{r-i}| \le  |f^{(r)}|\cdot \delta^{r-i}
\end{equation}
for every $x_0, x \in B$. (Note that $|f^{(r)}|$ is constant on $B$; here and in the sequel, we just write $|f^{(r)}|$
instead of $|f^{(r)}(x)|$ for $x \in B$.)

\medskip

Step 2.1: Set
\[
d(x_0, x) :=T^{\le r}_{f,x_0}(x) - f(x_0) = \sum_{i = 1}^r \frac{f^{(i)}(x_0)}{i!}(x - x_0)^i
\]
for $x_0, x \in B$. Then, for any $x_1, x_2, x_3 \in B$, we have
\begin{equation}\label{eq:additive}
|d(x_1, x_2) + d(x_2, x_3) - d(x_1, x_3)| < |f^{(r)}| \cdot |x_2-x_3|\cdot \delta^{r-1}.
\end{equation}

\medskip

Proof:
This is a straight-forward computation, consisting in using (\ref{eq:der}) to approximate the derivatives appearing in $d(x_2,x_3)$
by their Taylor series around $x_1$. The details are as follows. In the computation,
``$\approx$'' means that the norm of the difference is less than $ |f^{(r)}| \cdot |x_2-x_3|\cdot \delta^{r-1}$.
\begin{align*}
d(x_2, x_3) &= \sum_{i = 1}^r \frac{f^{(i)}(x_2)}{i!}(x_3 - x_2)^i \approx
\sum_{i = 1}^r\sum_{j = 0}^{r-i} \frac{f^{(i+j)}(x_1)}{i!\cdot j!}(x_2 - x_1)^j(x_3 - x_2)^i
\\
&=\sum_{k = 1}^r \underbrace{\sum_{\substack{i+j=k,\\i>0}}\frac{k!}{i!\cdot j!}(x_2 - x_1)^j(x_3 - x_2)^i}_{= ((x_2 - x_1) + (x_3 - x_2))^k - (x_2 - x_1)^k(x_3 - x_2)^0}\frac{f^{(i+j)}(x_1)}{k!}
\\
&=\sum_{k = 1}^r \big((x_3-x_1)^k - (x_2 - x_1)^k\big)\frac{f^{(k)}(x_1)}{k!} = d(x_1, x_3)  - d(x_1, x_2)
\end{align*}
\qed2.1

\medskip

For $x_1, x_2 \in K$, define
\[
x_1 \sim x_2 :\iff |f(x_2) - f(x_1) - d(x_1, x_2)| < |f^{(r)}| \cdot \delta^r
\]
(with some arbitrary convention if derivatives at $x_1 \notin B$ do not exist).
This relation is definable using only the parameter $|f^{(r)}| \cdot \delta^r \in \Gamma$.
Our aim is to verify the prerequisites of Lemma~\ref{lem:equivalence} (in the language where $|f^{(r)}| \cdot \delta^r$ has been added as a constant), but before that, let us verify that this then finishes the proof:
The lemma then implies that all elements of $B$ are $\sim$-equivalent. In particular,
for our original $x_0, x \in B$ satisfying $|x_0 - x| = \delta$, we obtain
\begin{equation}\label{eq:lem-finishes}
|f(x) - T^r_{f,x_0}(x)| =
|f(x) - f(x_0) - d(x_0, x)| < |f^{(r)}| \cdot \delta^r,
\end{equation}
i.e., (\ref{eq:c-higher-S}) holds, as desired.

\medskip

Step 2.2: The restriction of $\sim$ to $B \times B$ is an equivalence relation.

\medskip

Proof: Fix $x_0 \in B$ and define, for $x \in B$:
\begin{equation}\label{eq:ftilde}
\tilde f(x) := f(x) - d(x_0, x).
\end{equation}
Then, using ``$\approx$'' to mean that the difference has norm less than $|f^{(r)}|\cdot \delta^r$, we have:
\begin{equation}\label{eq:with_x0}
\begin{aligned}
 \tilde f(x_2) - \tilde f(x_1)
&= f(x_2) - f(x_1) - d(x_0, x_2) + d(x_0, x_1)\\
&\overset{(\ref{eq:additive})}{\approx} f(x_2) - f(x_1) - d(x_1, x_2)
\end{aligned}
\end{equation}
So $x_1 \sim x_2$ if and only if $\tilde f(x_2) \approx \tilde f(x_1)$, which is clearly an equivalence relation.
\qed2.2

\medskip

Step 2.3: Each proper subball $B' \subsetneq B$ is contained in a single equivalence class of $\sim$.

\medskip

Proof: This follows from our assumption that $f$ satisfies (\ref{eq:c-higher-S}) on $B'$ (using a similar computation as in (\ref{eq:lem-finishes})).
\qed2.3

\medskip

Step 2.4: $\sim$ has only finitely many equivalence classes on $B$.

\medskip

Proof: Consider $\tilde f$ as in (\ref{eq:ftilde}).
If $x_1, x_2 \in B$ satisfy $|x_1 - x_2|\le |\varpi|\cdot\delta$ (where $\varpi\in\cO_K$ is a uniformizer), then
in the ``$\approx$'' of (\ref{eq:with_x0}), the difference is less than $|f^{(r)}|\cdot |x_1 - x_2|\cdot \delta^{r-1} \le |f^{(r)}|\cdot|\varpi|\cdot\delta^r$, and the right hand side of (\ref{eq:with_x0}) satisfies
\[
|f(x_2) - f(x_1) - d(x_1, x_2)| = |f(x_2) - T^r_{f,x_1}(x_2)| \overset{\text{Case~1}}{<} |f^{(r)}|\cdot |x_1-x_2|^r.
\]
Thus $|\tilde f(x_2) - \tilde f(x_1)| < |f^{(r)}|\cdot|\varpi|\cdot\delta ^r$, i.e., for any open ball $B' \subset B$ of (open) radius $\delta$,
$\tilde f(B')$ is contained in an open ball of (open) radius $|f^{(r)}|\cdot|\varpi|\cdot\delta^{r}$. By Lemma~\ref{lem:balltosmall}, $\tilde f(B)$ is therefore contained in the union of finitely many closed balls of (closed) radius $|f^{(r)}|\cdot|\varpi|\cdot\delta^{r}$. By the characterization of $\sim$ from Step~2.2, each such closed ball corresponds exactly to one equivalence class of $\sim$, so we are done.
\qed2.4

\medskip

These were all the prerequisites needed for Lemma~\ref{lem:equivalence}, so this finishes the proof of Case~2 and hence of the entire Lemma.
\end{proof}

\begin{proof}[Proof of Lemma~\ref{lem:h:sc}]Namely:\\
From $\underbrace{|f -  T^{\le r}_{f,x_0} | < |f^{(r)}\cdot (x-x_0)^{r}|}_{(\ref{eq:c-higher-S})}$  to
$\underbrace{|f -  T^{\le r}_{f,x_0} | \leq |f^{(r)}\cdot (x-x_0)^{r+1}|/\radop(B)}_{(\ref{eq:c-higher})}$.
\medskip

We assume that $K$ is sufficiently saturated.

\medskip

Step 1: The corollary follows if we can verify that for every $a \in K$, Inequality (\ref{eq:c-higher}) holds for every ball $B \subset K \setminus \acl_K(a)$ that is $1$-next to $a$.

\medskip

Proof:
An easy computation shows that given any open ball $B$, Inequality (\ref{eq:c-higher}) holds on $B$ if and only if for every closed subball $B' \subset B$, and every \(x_0, x \in B'\), we have the following corresponding strict inequality:
\begin{equation}\label{eq:c-higher:c}
| f(x) -  T^{\le r}_{f,x_0}(x) | < |f^{(r)}(x_0)\cdot(x - x_0)^{r+1}| / \radcl(B').
\end{equation}
In particular, for $a$ and $B$ as in the assumption of Step~1, every closed subball of $B$ satisfies (\ref{eq:c-higher:c}).
This allows us to apply Lemma~\ref{lem:fullyT}, where we take $\cB$ to be the set of closed balls on which (\ref{eq:c-higher:c}) fails.
The lemma yields a finite $\emptyset$-definable set $C \subset K$
intersecting each ball in $\cB$ and hence also intersecting (as desired) each open ball where (\ref{eq:c-higher}) does not hold.\qed1

\medskip

For the remainder of the proof, let $a \in K$ be given and let $B$ be a ball which is disjoint from $\acl_K(a)$ and $1$-next to $a$. We may assume that $f^{(r)}$ is nowhere $0$ on $B$, since otherwise, it would be $0$ on all of $B$ (since $B \cap \acl_K(\emptyset) = \emptyset$) and Corollary~\ref{cor:high-ord-S} would yield $f(x)  = T^{\le r}_{f,x_0}(x)$ for $x \in B$.

Suppose that (\ref{eq:c-higher}) fails on $B$.
Choose $x, x_0 \in B$ witnessing this failure and fix, for the remainder of the proof,
\begin{equation}\label{eq:counter-d}
\delta := |x - x_0|
\end{equation}
and
\begin{equation}\label{eq:counter-a}
\alpha := \frac{ |f(x) -  T^{\le r}_{f,x_0}(x) | }{ |f^{(r)}(x_0)\cdot (x-x_0)^{r+1}|/\radop(B)};
\end{equation}
in other words, $\alpha > 1$ is the factor by which (\ref{eq:c-higher}) fails for $x_0$, $x$.
Moreover, set
\[
\gamma := \min\{\delta \cdot \alpha, \radop(B)\}
\]
(so that the ball $B_{< \gamma}(x_0)$ is contained in $B$ and is not much bigger than the smallest ball containing $x_0$ and $x$).

\medskip

Step 2: There exists an open ball $B' \subset B$ of radius $\gamma$ containing a ``$(\delta,\alpha)$-counter-example to (\ref{eq:c-higher})''
with the additional properties that $B'$ is $1$-next to some $a' \in K$ and $B' \cap \acl_K(a, a') = \emptyset$. By a $(\delta,\alpha)$-counter-example to (\ref{eq:c-higher}),
we mean a pair $x_0, x \in B'$ with
$|x - x_0| = \delta$ satisfying
\begin{equation}\label{eq:da}
|f(x) -  T^{\le r}_{f,x_0}(x) | \ge \alpha\cdot |f^{(r)}(x_0)\cdot (x-x_0)^{r+1}|/\radop(B).
\end{equation}

\medskip

Proof: If $\gamma = \radop(B)$, we can simply take $B' = B$ and $a' = a$, so now suppose that $\gamma < \radop(B)$.

We claim that the fact that $B$ contains a $(\delta,\alpha)$-counter-example already implies that
each open ball $B_{< \gamma}(z) \subset B$ of radius $\gamma$ contains a $(\delta,\alpha)$-counter-example.
Indeed, the set
\[
Z := \{z \in K \mid \exists x_0, x \in B_{< \gamma}(z)\colon (x_0, x)\text{ is a $(\delta,\alpha)$-counter-example}\}
\]
is definable using only some value group parameters (namely $\alpha$, $\delta$ and $\radop(B)$), so since $B \cap \acl_K(\emptyset) = \emptyset$, we have either $B \subset Z$ or $B \cap Z= \emptyset$. In particular, the existence of a single $(\delta,\alpha)$-counter-example in $B$ already implies $B \subset Z$,
which proves our claim.

Now fix any $a' \in B$ and fix an open ball $B'$ with $\radop(B') = \gamma$, which is $1$-next to $a'$ and such that $B' \cap \acl_K(a, a') = \emptyset$; such a $B'$ exists by saturation of $K$.
Since $\gamma < \radop(B)$, we have $B' \subset B$ and by the previous paragraph, $B'$ contains a $(\delta,\alpha)$-counter-example. Hence $a'$ and $B'$ are as desired.\qed2

\medskip

The remainder of the proof now consists in proving that there is no $(\delta,\alpha)$-counter-example on $B'$ (so that we have a contradiction). More precisely, fix $x_0, x \in B'$ with $|x_0 - x| = \delta$; our goal is to prove that (\ref{eq:da}) does not hold.

\medskip

Step 3: There exists a $d \in \acl_K(a')$ such that
the image $f^{(r)}(B')$ is either a ball $1$-next to $d$, or equal to the singleton $\{d\}$.

\medskip

Proof: Apply Proposition~\ref{prop:range} to compatibly prepare the domain and the image of $f^{(r)}$ using finite $a'$-definable sets $C$ and $D$, where we additionally impose $a' \in C$. Since $B' \cap \acl_K(a') = \emptyset$, $B'$ is a ball $1$-next to $C$ and the claim follows.
\qed3

\medskip

Step~4: We apply Corollary~\ref{cor:high-ord-S} to the function
\[g(x) := f(x) - \frac{d}{r!}\cdot x^r.
\]
Since $g$ is $\{d\}$-definable and $B' \cap \acl_K(d) = \emptyset$, we obtain
\begin{equation}\label{eq:gprep}
|f(x) -  T^{\le r}_{f,x_0}(x) | = |g(x) - T^{\le r}_{g,x_0}(x)| <  | g^{(r)}(x_0)\cdot(x-x_0)^r|
\end{equation}
(where $x_0$ and $x$ are as fixed above Step~3). The first equality follows from the fact that $f$ and $g$ differ by a polynomial of degree $r$, so their $r$-th Taylor approximations differ by the same polynomial.

\medskip

Step~5a: If $f^{(r)}$ is non-constant on $B'$, we obtain
\begin{equation}\label{eq:s6}
|f(x) - T^{\le r}_{f,x_0}(x)| < \radop(f^{(r)}(B')) \cdot |x-x_0|^r.
\end{equation}

\medskip

Proof: By definition of $g$ and by our choice of $d$ (in Step~3), we have
\begin{equation}\label{eq:use-d}
|g^{(r)}(x_0)| = |f^{(r)}(x_0) - d| = \radop(f^{(r)}(B')).
\end{equation}
Combining this with (\ref{eq:gprep}) yields (\ref{eq:s6}).
\qed5a

\medskip

Step~5b: If $f^{(r)}$ is constant on $B'$, we obtain
\[
|f(x) - T^{\le r}_{f,x_0}(x)| = 0.
\]

\medskip

Proof: Instead of (\ref{eq:use-d}), we have $g^{(r)}(x_0) = f^{(r)}(x_0) - d = 0$; apply this in the same way as in Step 5a.\qed5b

\medskip

In the constant case (as in Step 5b), we are already done for the lemma, since the right hand side of (\ref{eq:da}) is non-zero, a contradiction. (Recall that we assume $f^{(r)} \ne 0$.)

The last ingredient for the non-constant case is the following:

\medskip

Step 6: We have $\radop(f^{(r)}(B')) \le \alpha\cdot |f^{(r)}(x_0)|\cdot|x-x_0|/\radop(B)$.

\medskip

Proof: Using that $\rv(f^{(r)})$ is constant on $B$, we obtain $\radop(f^{(r)}(B)) \le |f^{(r)}(x_0)|$.
From this and by applying Lemma~\ref{lem:gammaLin} to $f^{(r)}$, we deduce that $\radop(f^{(r)}(B')) \le |f^{(r)}(x_0)|\cdot \radop(B')/\radop(B)$.
Combining this with $\radop(B') = \gamma \le \alpha\cdot |x - x_0|$ yields the claim.\qed6

\medskip

Now Steps 5a and 6 together imply that (\ref{eq:da}) fails, as desired, so we are done.
\end{proof}

\begin{proof}[Proof of Lemma \ref{lem:h:ct}]Namely:\\
From  $\underbrace{|f -  T^{\le r}_{f,x_0} | \leq |f^{(r)}\cdot (x-x_0)^{r+1}|/\radop(B)}_{(\ref{eq:c-higher})}$  to
$\underbrace{|f -  T^{\le r}_{f,x_0} | \leq |f^{(r+1)}\cdot (x-x_0)^{r+1}|}_{(\ref{eq:t-higher})}$.
\medskip

Let $\cB$ be the set of all closed balls on which (\ref{eq:t-higher}) does not hold.
The strategy is to use Lemma~\ref{lem:fullyT} to find a finite $\emptyset$-definable $C$
meeting every $B \in \cB$. Note that then we are done,
since if (\ref{eq:t-higher}) fails for some $x, x_0 \in B'$, where $B'$ is a ball $1$-next
to $C$, then (\ref{eq:t-higher}) also fails on a ball from $\cB$, namely the smallest (closed) ball containing $x$ and $x_0$.

So as needed for Lemma~\ref{lem:fullyT}, let $a \in K$ be given and let $B$ be an open ball $1$-next to $a$ satisfying $B \cap \acl_K(a) = \emptyset$. We need to verify that (\ref{eq:t-higher}) holds on $B$.

By applying Proposition~\ref{prop:range} to $f^{(r)}$ and to a set $C_0$ containing $a$, we find a $d \in \acl_K(a)$ such that either $f^{(r)}(B) = \{d\}$ or
$f^{(r)}(B)$ is a ball $1$-next to $d$.

Now apply Corollary~\ref{cor:high-ord} to $g(x) := f(x) - \frac{d}{r!}x^r$. Since $B$ is disjoint from the algebraic closure of the parameters used to define $g$,
we obtain
\begin{equation}\label{eq:gT}
|f(x) - T^{\le r}_{f,x_0}(x)| = |g(x) - T^{\le r}_{g,x_0}(x)| \le   |g^{(r)}(x_0)\cdot (x - x_0)^{r+1}| / \radop(B)
\end{equation}
for every $x, x_0 \in B$. As before, the first equality holds because \(f\) and \(g\) differ by a polynomial of degree \(r\). To finish the proof, it now remains to bound the right hand side of (\ref{eq:gT}) by $|f^{(r+1)}(x_0)\cdot (x-x_0)^{r+1}|$.

If $f^{(r)}(B) = \{d\}$, then $g^{(r)}(x_0) = 0$ and we are done.
Otherwise,  $g^{(r)}(B)$ is a ball $1$-next to $0$, so $\radop(g^{(r)}(B)) = |g^{(r)}(x_0)|$.
Moreover, by Corollary~\ref{cor:gammaLin} applied to $g^{(r)}$, we have $\radop(g^{(r)}(B)) = \radop(B)\cdot |g^{(r+1)}(x_0)| = \radop(B)\cdot |f^{(r+1)}(x_0)|$.
Combining these two equations yields the desired bound on the right hand side of (\ref{eq:gT}).
\end{proof}

\section{Resplendency}\label{sec:respl}

The main goal of this section is to show that Hensel minimality behaves well with respect to expansions of the structure by predicates living only on Cartesian powers of $\RV$, one of the main results (Theorem~\ref{thm:resp:h}) being that if the $\cL$-theory of a valued field $K$ is $0$-, $1$- or $\omega$-h-minimal, then so is its $\cLe$-theory, where $\cLe$ is the language of such an $\RV$-expansion of $K$.

The key ingredient to Theorem~\ref{thm:resp:h} is Proposition~\ref{prop:equiv}, which in some sense is even stronger: Any set $X \subset K$ definable in the bigger language $\cLe$ can already be prepared by a finite set definable in the smaller language $\cL$. Using this, it often becomes possible, given a completely arbitrary set $Z \subset \RV^k$ to ``without loss assume that $Z$ is definable''. This turns out to be pretty useful to get rid of some technicalities related to cell decomposition in valued fields; the preparations for this are done in Subsection~\ref{sec:alg:skol}.

Under the assumption of $\omega$-h-minimality,
preparation can also be generalized to more general leading term structures, namely:
One can define $\RV_{I} := (K^\times/(1+I)) \cup \{0\}$ for arbitrary proper definable ideals $I \subset \cO_K$, and for most such $I$, any subset of $K$ which is definable using parameters from $\RV_I$
can be ``$I$-prepared'' (see Theorem~\ref{thm:all:I}). Using this, we deduce that if the theory of a valued field is $\omega$-h-minimal, then so is the
theory of the field with any coarsened (not necessarily definable) valuation (Corollary~\ref{cor:coarse}). In fact, in \cite[Theorem~2.2.7]{CHRV}, it is shown that also $1$-h-minimality is preserved under coarsenings of the valuation, and furthermore, some resplendency results in the mixed characteristic case are developed in \cite{CHRV}.

Note that the proofs in this section need somewhat deeper methods from model theory than the remainder of the paper. In particular, the emphasis shifts from a geometric description of definable sets to questions revolving around the extension of automorphisms.

\subsection{Resplendency for a fixed ideal}
\label{sec:fixI}

As convened earlier, $K$ is a valued field of equi-characteristic zero, considered as a structure in a language $\cL$ containing $\Lval $.
For this entire subsection, we fix a proper definable (with parameters) ideal $I \subset \cO_K$.
We start by defining the $I$-version of preparation.

\begin{defn}[$I$-preparing sets]\label{defn:Inext}
\begin{enumerate}
 \item We define $\RV_I$ to be the disjoint union of the
quotient $K^\times/(1+I)$ with $\{0\}$, and we write $\rv_I$ for the map $K\to \RV_I$ which extends the projection map $K^\times\to K^\times/(1+I)$ by sending $0$ to $0$.
 \item
 We say that a ball $B \subset K$ is \emph{$I$-next} to $c$ for some $c\in K$ if $B = c + \rv_I^{-1}(\xi)$
for some (nonzero) element $\xi$ of $\RV_I$.
We say that $B$ is \emph{$I$-next} to $C$ for some finite (non-empty) set $C\subset K$
if $B$ equals $\bigcap_{c\in C} B_c$
with $B_c$ a ball $I$-next to $c$ for each $c\in C$.
 \item Let $C$ be a finite non-empty subset of $K$.
We say that a set $X\subset K$ is \emph{$I$-prepared} by $C$ if belonging to $X$ depends only on the tuple $(\rv_{I}(x-c))_{c\in C}$,
or, equivalently, if every ball $I$-next to $C$ is either contained in $X$ or disjoint from $X$.
\end{enumerate}
\end{defn}

\begin{remark}
If $I = B_{<\lambda}(0)$ for some $\lambda \le 1$ in $\Gamma_K$, then of course we have $\RV_I = \RV_\lambda$, $I$-next means $\lambda$-next, and $I$-prepared means $\lambda$-prepared.
\end{remark}

By considering $I$ as a member of a $\emptyset$-definable family of proper ideals of $\cO_K$, we see that
$\RV_{I}$ is a definable subset of an imaginary sort. In particular, it makes sense to work using parameters from $\RV_I$.

\begin{defn}[Having $I$-preparation]\label{defn:Iprepared}
We say that $K$ has \emph{$I$-preparation}
if for every set $A\subset K$ and every $(A\cup \RV_I)$-definable subset $X\subset K$,
there exists a finite $A$-definable set $C \subset K$ such that
$X$ is $I$-prepared by $C$.
\end{defn}

\begin{remark}\label{rem:prep:vs:h}
Note how this is related to Hensel minimality: $\Th(K)$ is $0$-h-minimal iff every $K' \equiv K$ has $\cM_{K'}$-preparation, and $\Th(K)$ is $\omega$-h-minimal iff every $K' \equiv K$ has $B_{<\lambda}(0)$-preparation for every $\lambda \le 1$ in $\Gamma^\times_{K'}$.
\end{remark}

By ``resplendent $I$-preparation'', we mean that one can $I$-prepare sets that are definable in arbitrary expansions of $\cL$ by predicates on $\RV_I$. More precisely:

\begin{defn}[Resplendent $I$-preparation]\label{defn:Iresplendent}
\begin{enumerate}
\item
An \emph{$\RV_I$-expansion} of $K$ is an expansion obtained by adding (any)
predicates which live on Cartesian powers of the (imaginary) definable set $\RV_I$.
\item
We say that $K$ has \emph{resplendent $I$-preparation} if for every set $A\subset K$, for every $\RV_I$-expansion of $K$ with language $\cL' \supset \cL$, and
for every $\cL'(A)$-definable subset $X\subset K$,
there exists a finite $\cL(A)$-definable set $C \subset K$ such that $X$ is $I$-prepared by $C$.
\end{enumerate}
\end{defn}

Note that we intentionally require $C$ to be definable in the smaller language $\cL$.

For the remainder of this subsection, we will assume that $I$ is definable without parameters. In this case,
it also makes sense to introduce the notions of $I$-preparation for the theory of $K$:

\begin{defn}[$I$-preparation for theories]\label{defn:I:T}
Suppose that $I$ is $\emptyset$-definable (as convened for the remainder of this subsection).
\begin{enumerate}
\item
Given $K' \equiv K$, we write $I_{K'}$ for the ideal of $\cO_{K'}$ defined by the formula which defines $I$ in $K$.
\item
We say that $\Th(K)$ has \emph{$I$-preparation} if every model $K' \equiv K$ has $I_{K'}$-preparation.
\item
We say that $\Th(K)$ has \emph{resplendent $I$-preparation} if every model $K' \equiv K$ has resplendent $I_{K'}$-preparation.
\end{enumerate}
\end{defn}

(In particular, for $\Th(K)$, having $\cM_K$-preparation is the same as $0$-h-minimality.)

Since adding parameters from $\RV_I$ is a specific kind of $\RV_I$-expansion, resplendent $I$-preparation clearly implies $I$-preparation. The central result of this subsection is the converse given in the following proposition:

\begin{prop}[Resplendency]\label{prop:equiv}
Suppose that $I$ is $\emptyset$-definable. Then the following are equivalent:
\begin{enumerate}[(i)]
\item $\Th(K)$ has $I$-preparation.
\item $\Th(K)$ has resplendent $I$-preparation.
\end{enumerate}
\end{prop}

Note that the proposition in particular implies that $0$-h-minimal theories have resplendent $\cM_K$-preparation;
in Theorem~\ref{thm:all:I}, we will see a strong version of this for $\omega$-h-minimality.

The proof of (i) $\Rightarrow$ (ii) requires a number of lemmas which we will prove now. It will be sufficient to consider models $K$ that are sufficiently saturated, so from now on, we fix such a sufficiently saturated $K$: We assume $K$ to be $\kappa$-saturated for some $\kappa >|\cL|$. As usual, we call a set \emph{small} if its cardinality is less than $\kappa$.

\begin{conv}
For the remainder of Subsection~\ref{sec:fixI}, we consider $\cL$ as a genuine two-sorted language, with sorts $K$ and $\RV_{I}$.
\end{conv}

Let us first rephrase preparation in terms of definability in a certain sublanguage, namely:
\begin{defn}[$\Lbas$]
Let $\Lbas$ be the language $\{0, +, -\} \cup \{s\cdot \mid s\in \QQ\}$ of $\QQ$-vector spaces on the valued field $K$
(where ``$s\cdot$'' denotes multiplication by $s$), together with
the sort $\RV_{I}$ and the map $\rv_{I}$.
\end{defn}

\begin{lem}[Preparation in terms of $\Lbas$]\label{lem:prep eq qf}
For every (not necessarily definable) $X\subset K$ and every $\QQ$-vector space $A \subset K$, the following are equivalent:
\begin{enumerate}
\item There exists a finite set $C\subset A$ such that $X$ is $I$-prepared by $C$.
\item There exists an $\RV_I$-expansion of $K$ with language $\Lebas \supset \Lbas$ such that $X$ is defined by a quantifier free $\Lebas(A)$-formula.
\item There exists an $\RV_I$-expansion of $K$ with language $\Lebas \supset \Lbas$ such that $X$ is defined by a field quantifier free $\Lebas(A)$-formula.
\end{enumerate}
\end{lem}

\begin{proof}
(2) $\Rightarrow$ (3) is clear.

(1) $\Rightarrow$ (2): Let us assume that (1) holds and let $C = \{c_1,\ldots,c_k\}$, $f(x) := (\rv_{I}(x-c_1),\ldots,\rv_{I}(x-c_k))$ and $Y := f(X) \subset \RV_{I}^k$. Then, because $X$ is $I$-prepared by $C$, $X = f^{-1}(Y)$, and this is quantifier free definable in $\Lbas(A)\cup\{Y\}$.

(3) $\Rightarrow$ (1): Every field quantifier free $\Lebas(A)$-formula in a single variable $x$ is equivalent
to one of the form $\phi(\rv_{I}(m_1x + c_1),\ldots,\rv_{I}(m_\ell x + c_\ell))$, where $\phi$ is a formula living entirely in the sort $\RV_I$, $m_i \ne 0$ are rational numbers and $c_i$ are elements of $A$. Since $\rv_{I}(x + c_i/m_i)$ determines
$m_i\rv_{I}(x + c_i/m_i) = \rv_{I}(m_i x+ c_i)$,
$X$ is $I$-prepared by $C = \{- c_1/m_1, \dots, -c_\ell/m_\ell\} \subset A$.
\end{proof}

We now recall some general model theoretic notions.

\begin{notn}
In the remainder of this subsection, we use the following conventions common in model theory:
\begin{itemize}
 \item
Given a set $A$ and a tuple of variables $x$, we write $A^x$ for the Cartesian power of $A$ corresponding to the length of the tuple of variables $x$.
\item
Given sets $A$ and $B$, we sometimes write $AB$ for their union, and we freely interpret tuples as sets.
\end{itemize}
\end{notn}

\begin{defn}[Partial (elementary) isomorphisms]
Suppose that $\cLarb$ is an arbitrary language, $M$ and $N$ are $\cLarb$-structures, $A\subset M$, $B\subset N$ and $f \colon A\to B$ a bijection.
\begin{itemize}
\item We say that $f$ is a \emph{partial $\cLarb$-isomorphism} if for every quantifier free $\cLarb$-formula $\phi(x)$ and $a\in A^{x}$, $M\models\phi(a)$ if and only if $N\models\phi(f(a))$.
\item We say that $f$ is a \emph{partial elementary $\cLarb$-isomorphism} if for every $\cLarb$-formula $\phi(x)$ and $a\in A^{x}$, $M\models\phi(a)$ if and only if $N\models\phi(f(a))$.
\end{itemize}
\end{defn}

Note that any partial $\cLarb$-isomorphism has a unique extension to the $\cLarb$-structure generated by its domain. We implicitly identify $f$ and this extension.

In the following, for a subset $A \subset K$, $\langle A \rangle_\QQ$ denotes the $\QQ$-sub-vector space  of $K$ generated by $A$.

\begin{remark}\label{rem:join iso}
For any set $A \subset K$, the $\Lbas$-substructure of $K$ generated by $A$ consists of $\langle A\rangle_{\QQ}$
together with its image $\rv_{I}(\langle A\rangle_{\QQ})$ in $\RV_{I}$.
In particular, since $\Lbas$ has no language on $\RV_I$ (except for the maps $\rv_{I}$), we have:
To obtain that a sort-preserving map $f\colon A_1 \to A_2$ (for some $A_1, A_2 \subset K \cup \RV_I$)
is a partial $\Lbas$-automorphism, it suffices to verify that the restriction $f|_{A_1 \cap K}$ is a partial $\Lbas$-automorphism and that on $\tilde A_1 := A_1 \cap \rv_I(A_1 \cap K)$, the map induced by $f|_{A_1 \cap K}$ agrees with $f|_{\tilde A_1}$.
\end{remark}

\begin{lem}[Preparation in terms of partial isomorphisms]\label{lem:prep eq iso}
Let $A \subset K$ be a small $\QQ$-sub-vector space.
The following are equivalent:
\begin{enumerate}[(i)]
\item Any $\cL(A\cup \RV_I)$-definable set $X\subset K$ can be $I$-prepared by some finite set $C\subset A$.
\item For all $A_2\subset K$, $c_1, c_2 \in K$ and all (potentially large) $B_1, B_2\subset \RV_I$
with $\rv_I(\langle A,c_1\rangle_\QQ) \subset B_1$,
if $f\colon A B_1 c_1 \to A_2 B_2 c_2$ is a partial $\Lbas$-isomorphism sending $c_1$ to $c_2$
whose restriction $\restr{f}{AB_1}$ is a partial elementary $\cL$-isomorphism, then the entire $f$
is a partial elementary $\cL$-isomorphism.
\item For all $c_1, c_2 \in K$ and all (potentially large) $B\subset \RV_I$ containing $\rv_I(\langle A,c_1\rangle_\QQ)$, any partial $\Lbas(A\cup B)$-isomorphism $f\colon\{c_1\} \to \{c_2\}$, is a partial elementary $\cL(A\cup B)$-isomorphism.
\end{enumerate}
\end{lem}

\begin{remark}\label{rem:wlogRVI}
If (ii) holds for $B_1 = B_2 = \RV_I$, then it also holds in general, since by Remark~\ref{rem:join iso},
the partial $\Lbas$-isomorphism $f\colon A B_1 c_1 \to A_2 B_2 c_2$ extends to $f\colon A c_1 \cup \RV_I \to A_2 c_2 \cup  \RV_I$.
Analogously, we may assume $B = \RV_I$ in (iii).
\end{remark}

\begin{proof}[Proof of Lemma~\ref{lem:prep eq iso}]
(i) $\Rightarrow$ (iii):
Let $f$ be as in (iii). We have to check that for every $\cL(A\cup B)$-definable set $X\subset K$, $c_1 \in X$ if and only if $c_2\in X$. By (i), there exists a finite $C \subset A$ such that $X$ is $I$-prepared by $C$. Since $f$ is an $\Lbas(A\cup B)$-isomorphism and $B$ contains $\rv_I(\langle A,c_1\rangle_\QQ)$,
for all $a\in C$ and all $r\geq 1$, we have
\[
\rv_{I}(c_2-a) = \rv_{I}(f(c_1)-f(a)) = f(\rv_{I}(c_1-a)) = \rv_{I}(c_1-a).
\]
Since $X$ is $I$-prepared by $C$, it follows that $c_1 \in X$ if and only if $c_2\in X$.

\medskip

(iii) $\Rightarrow$ (ii):
Let $f$ be as in (ii). Since $f$ is $\cL$-elementary if and only its restriction to every finite domain is, we may assume $B_i$ small. Using the assumption that $\restr{f}{AB_1}$ is $\cL$-elementary, we can extend $(\restr{f}{AB_1})^{-1}$ $\cL$-elementarily to some $g$ defined at $c_2$. Let $c'_1 := g(c_2)$. Then $g\circ f\colon \{c_1\} \to \{c'_1\}$ is a partial $\Lbas(A\cup B_1)$-isomorphism. Since $\rv_I(\langle A, c_1\rangle_\QQ)\subset B_1$, it follows by (iii) that $g\circ f$ is an elementary $\cL$-isomorphism. As $g$ is also $\cL$-elementary, so is $f$.

\medskip

(ii) $\Rightarrow$ (i):
Let $X$ be as in (i), and let $B \subset \RV_I$ be a finite subset such that $X$ is $\cL(A \cup B)$-definable.

Consider any $c_1, c_2 \in K$ which have the same qf-$\LbasLA$-type over $A \cup B$, where $\LbasLA \supset \Lbas$ is the language of the expansion of $K$ by the full $\cL(A)$-induced structure on $\RV_I$.
Then the map $f\colon c_1 \to c_2$ is an $\LbasLA(A \cup B)$-isomorphism and extends to $f\colon A B_1 c_1 \to A_2 B_2 c_2$, where
$B_i := B \cup \rv_I(\langle A,c_i\rangle_\QQ)$. By definition of $\LbasLA$, the restriction $f|_{AB_1}$ is $\cL$-elementary, so by (ii), also the entire $f$ is $\cL$-elementary.
Since moreover $f$ is the identity on $A \cup B$, this implies that $c_1$ and $c_2$ have the same $\cL$-type over $A \cup B$.

We just proved that the $\cL(A \cup B)$-type of any element $c \in K$ is implied by its qf-$\LbasLA(A \cup B)$-type.
By a classical compactness argument (cf. the proof of \cite[Theorem 3.2.5]{TentZiegler}), it follows that any $\cL(A\cup B)$-formula in one valued field variable is equivalent to a quantifier free $\LbasLA(A\cup B)$-formula. In particular, this applies to our set $X$. Since $\LbasLA(B)$ is the language of an $\RV_I$-expansion of $K$, our claim follows from Lemma~\ref{lem:prep eq qf}.
\end{proof}

\begin{lem}[Back and forth over $\RV_I$]\label{lem:b-a-f over RV}
Suppose that $K$ has $I$-preparation. Then the set of partial elementary $\cL$-isomorphisms $f\colon \RV_I\cup A_1 \to \RV_I\cup A_2$ (where $A_1, A_2$ run over all small subsets of $K$) has the back-and-forth.
\end{lem}

Recall that ``having the back and forth'' means: for any such $f$ and any $c_1 \in K \setminus A_1$,
$f$ can be extended to $c_1$ while staying in that set of maps, and similarly for \(f^{-1}\).

\begin{proof}
Let $f$ be as above. Since partial elementary isomorphisms can always be extended to the algebraic closure of their domain, we may assume that $\acl_K(A_i) = A_i$.
By $I$-preparation, statement (i) of Lemma~\ref{lem:prep eq iso} now holds for $A = A_i$.
Pick any $c_1\in K$, set $B_1 := \rv_I(\langle A_1, c_1\rangle_{\QQ})$
and let $c_2\models f_{*}\mathrm{qftp}_{\Lbas}(c_1/A_1B_1)$. (Such a $c_2$ exists, since $K$ is $|A_1B_1|^+$-saturated.) Let $g$ extend $f$ by sending $c_1$ to $c_2$. By construction, $\restr{g}{A_1B_1c_1}$ is a partial $\Lbas$-isomorphism. This implies that the entire $g$ is a partial $\Lbas$-isomorphism
(cf.\ Remark~\ref{rem:join iso}). Note also that $\restr{g}{A_1\RV_I} = f$ is a partial elementary $\cL$-isomorphism. Hence, by Lemma~\ref{lem:prep eq iso}~(ii), the entire $g$ is $\cL$-elementary.
\end{proof}

We are now ready for the proof of the central result of this subsection:

\begin{proof}[Proof of Proposition~\ref{prop:equiv}]
Set $\cT := \Th(K)$. Recall that the implication (ii) to (i) is trivial, so let us assume (i) (namely, $\cT$ has $I$-preparation) and prove (ii) (namely, $\cT$ has resplendent $I$-preparation):
Suppose that $K$ is an arbitrary model of $\cT$, let $\cLe \supset \cL$ be the language of an $\RV_{I}$-expansion of $K$, let $A$ be a subset of $K$ and let $X \subset K$ be an $\cLe(A)$-definable set;
we need to find a finite $\cL(A)$-definable set $C \subset K$ which $I$-prepares $X$. Since the condition that $C$ $I$-prepares $X$ is first order, we may replace $K$ by a sufficiently saturated elementary extension. For the remainder of the proof, we fix such a $K$.

We may moreover assume that $A = \acl_{\cL,K}(A)$. Indeed, if we find a finite, $\cL(\acl_{\cL,K}(A))$-definable set $C$ which $I$-prepares $X$, then the union of the (finitely many) conjugates of $C$ over $A$ is an $\cL(A)$-definable finite set $I$-preparing $X$.
Therefore, we may apply Lemma~\ref{lem:prep eq iso} (and Remark~\ref{rem:wlogRVI}), meaning that it suffices to show that every partial $\Lbas(A\cup\RV_I)$-isomorphism $f\colon\{c_1\} \to \{c_2\}$ is an elementary $\cLe(A\cup\RV_I)$-isomorphism.

By (i) and Lemma~\ref{lem:prep eq iso} applied in $\cL$, such an $f$ is a partial elementary $\cL(A\cup\RV_I)$-isomorphism, and since $f$ is the identity on $\RV_I$, it is a partial $\cLe(A\cup\RV_I)$-isomorphism. It remains to show that $f$ preserves all $\cLe$-formulas and not just the quantifier free ones.

Let $\cF$ be the class of partial elementary $\cL(\RV_I)$-isomorphisms with small domains $A \subset K$. By Lemma~\ref{lem:b-a-f over RV}, $\cF$ has the back-and-forth. Moreover, any $g\in\cF$ is also a partial $\cLe$-isomorphism, as it fixes $\RV_I$. Recall that if a class of partial $\cLarb$-isomorphisms has the back and forth
(for any given language $\cLarb$), by an easy induction on the structure of formulas, those partial isomorphisms are automatically $\cLarb$-elementary.
Thus any $g\in\cF$, and in particular $f$, is a partial elementary $\cLe$-isomorphism, which is what we had to show.
\end{proof}

We now mention some easy consequences of Proposition~\ref{prop:equiv} and its proof.

\begin{cor}[$\RV$-expansions preserve $\acl$]\label{cor:acl=acl}
Suppose that $\Th(K)$ has $I$-preparation, for some $\emptyset$-definable proper ideal $I \subset \cO_K$.
Then for any $\RV_I$-expansion of $K$ with language $\cLe \supset \cL$ and any $A \subset K$, we have $\acl_{\cL',K}(A) = \acl_{\cL,K}(A)$.
In particular:
\begin{enumerate}
\item If $\Th(K)$ is $0$-h-minimal, then for any $\RV$-expansion of $K$ with language $\cLe \supset \cL$ and any $A \subset K$, we have $\acl_{\cL',K}(A) = \acl_{\cL,K}(A)$.
\item If $\Th(K)$ is $\omega$-h-minimal, then
for any $\lambda\in\Gamma_K^\times$, any $\RV_\lambda$-expansion of $K$ with language $\cLe \supset \cL$ and any $A \subset K$, we have $\acl_{\cL',K}(A) = \acl_{\cL,K}(A)$.
 \end{enumerate}

\end{cor}
\begin{proof}
Any $b \in \acl_{\cLe,K}(A)$ is an element of a finite $\cLe(A)$-definable set $X \subset K$.
By Proposition~\ref{prop:equiv}, $X$ can be $I$-prepared by a finite $\cL(A)$-definable set $C$. This implies $X \subset C$ and hence $b \in \acl_{\cL,K}(A)$.
\end{proof}

\begin{remark}
We already saw some stable embeddedness results in Proposition~\ref{prop:stab} and Remark~\ref{rem:stab}. For similar reasons, $I$-preparation implies stable embeddedness of $\RV_I$. Alternatively, one may deduce the stable embeddedness from Lemma~\ref{lem:b-a-f over RV} and \cite[Appendix, Lemma 1]{ChaHru-ACFA}\footnote{This uses the existence of fully saturated models; getting rid of those is left as an exercise to the reader.}.
\end{remark}

\private{Proof of stab. emb of $\RV_I$ without using fully saturated models:

Let $K$ be a sufficiently saturated model of $\cT$.
It suffices to prove that for every finite $A \subset K$ and for every tuple $\xi$ of elements of $\RV_I$, the type $\tp(\xi/A)$ is already determined by $\tp(\xi/B)$, where $B :=\rv_I(\langle A \rangle_\QQ)$.
Indeed, if $X$ is an $A$-definable subset of a product of some of the sorts of $\RV_I$, then using compactness, the above implies that each $\xi \in X$ is contained in a $B$-definable subset of $X$ and using compactness once more, the union of finitely many of those subsets is equal to $X$.

So now let $A$ and $B = \rv_I(\langle A \rangle_\QQ)$ be as above and suppose that $\xi_1, \xi_2$ are two tuples from $\RV_I$ which have the same $\cL$-type over $B$.
Define $f\colon A \cup B \cup \{\xi_1\} \to A \cup B \cup \{\xi_2\}$ to be the identity on $A \cup B$ and to send $\xi_1$ to $\xi_2$.
This map is an $\Lbas$-isomorphism and its restriction $f|_{B\xi_1}$ is $\cL$-elementary. Applying Lemma~\ref{lem:prep eq iso} repeatedly
(once for each of the finitely many elements of $A$), we deduce that the entire $f$ is $\cL$-elementary, hence showing that $\xi_1$ and $\xi_2$ have the same type over $A$.
}

We conclude this subsection by the result that various notions of Hensel minimality automatically pass to $\RV$-expansions of the structure.

\begin{thm}[$\RV$-expansions preserve Hensel minimality]\label{thm:resp:h}
Let $K$ be a valued field of equi-characteristic $0$ in a language $\cL$ containing the language $\Lval$ of valued fields.
Fix $\ell \in \{0, 1, \omega\}$ and suppose that $\Th_\cL(K)$ is $\ell$-h-minimal.
Moreover, fix an $\RV$-expansion of $K$ with language $\cLe \supset \cL$ (i.e., an expansion by predicates on Cartesian powers of $\RV$).
Then $\Th_{\cLe}(K)$ is also $\ell$-h-minimal.
\end{thm}

Whereas the cases $0$ and $\omega$ follow easily from Proposition~\ref{prop:equiv}, the only proof we have for $\ell = 1$ makes a detour through the criterion given in Theorem~\ref{thm:tame2vf} (which is generalized for $\ell \ge 2$ in \cite{Verm:h-min}). In \cite{Verm:h-min}, \Cref{thm:resp:h} is generalized to hold for all $\ell$.

\begin{proof}[Proof of Theorem~\ref{thm:resp:h}]
The case $\ell = 0$ follows directly from Proposition~\ref{prop:equiv}: Since $\Th_{\cL}(K)$ has $\cM_K$-preparation, it even has resplendent $\cM_K$-preparation, so every $\cLe(A)$-definable set $X \subset K$ (for $A \subset K$) can be prepared by a finite $\cL(A)$-definable (and hence in particular $\cLe(A)$-definable) subset of $K$.

For the case $\ell = \omega$, we use a similar argument, where $\cM_K$ is replaced by $I = B_{< \lambda}(0)$ for arbitrary $\lambda \le 1$ in $\Gamma^\times_K$. Note that indeed, the given $\RV$-expansion of $K$ can also be considered as an $\RV_\lambda$-expansion, by pulling back all the new predicates along the canonical map $\RV_\lambda \to \RV$.

Finally, in the case $\ell = 1$, we use the $1$-h-minimality criterion given in Theorem~\ref{thm:tame2vf}: Let $f\colon K \to K$ be $\cLe(A)$-definable, for $A \subset K \cup \RV$. We need to prove the two conditions (T1) and (T2) stated in the theorem.

First of all, note that we already know that $\Th_{\cLe}(K)$ is $0$-h-minimal.

By Corollary~\ref{cor:acl=acl} (applied with $I = \cM_K$), for every $x \in K$, we have $f(x) \in \acl_{\cL,K}(A \cup \{x\})$.
Let $C_x$ be an $\cL(A)$-definable family of finite sets such that $f(x) \in C_x$ for all $x \in K$. Lemma~\ref{lem:average} yields an $\cL(A)$-definable family of injective maps $g_x\colon C_x \to \RV^k$. Using this, we define maps $h\colon K \to \RV^k, x \mapsto g_x(f(x))$ and $\tilde f\colon K \times \RV^k \to K$ with $\tilde f(x, \xi) = g_x^{-1}(\xi)$
if $\xi \in g_x(C_x)$ and $\tilde f(x, \xi) = 0$ otherwise.
Note that we obtain $f(x) = \tilde f(x, h(x))$ for all $x \in K$.
Also note that $h$ is $\cLe(A)$-definable and $\tilde f$ is $\cL(A)$-definable.

Applying Corollary~\ref{cor:prep} (in $\cLe$) to the graph of $h$ yields a finite $\cLe(A)$-definable set $C'$ such that $h$ is constant on every ball $B$ $1$-next to $C'$. Moreover, using Lemma~\ref{lem:gammaLin} in $\cL$, we find an $\cL(A)$-definable family of sets $C_\xi$ preparing $x \mapsto \tilde f(x, \xi)$ in the sense of that lemma for every $\xi \in \RV^k$. Now let
$C$ be the union of $C'$ and all those $C_\xi$. That union is still finite (using Corollary~\ref{cor:finiterange:set} in $\cL$), and it prepares $f$ in the way Condition (T1) requires it. Indeed, on each fixed ball $B$ $1$-next to $C$, we have $f(x) = \tilde f(x, \xi)$ for one fixed $\xi$, our application of Lemma~\ref{lem:gammaLin} ensured that $x \mapsto \tilde f(x, \xi)$ is prepared in the desired sense.

Condition (T2) now is also clear: For each $\xi$, $\tilde f(\cdot, x)$ has only finitely many infinite fibers (by Lemma~\ref{lem:fin-inf} in $\cL$). Taking the union of all those sets (for all $\xi$) still yields a finite set.
\end{proof}

\subsection{Changing the ideal}

We now switch back to considering $\cL$ as a single-sorted language.
In this subsection, we prove that $\omega$-h-minimality implies $I$-preparation for most ideals $I$, where ``most'' means ``$\Gamma_K$-open'' in the following sense:

\begin{defn}[Notions of openness]\label{defn.G-open}
By an \emph{open ball ideal} in $\cO_K$, we mean an ideal of the form $B_{<\lambda}(0)$ for some
$\lambda \le 1$ in $\Gamma_K$. By a \emph{$\Gamma_K$-open ideal} in $\cO_K$, we mean an ideal which is equal to the union of all open balls ideals it contains.
\end{defn}

In other words, an ideal is $\Gamma_K$-open if its image in $\Gamma_K$ is open with respect to the interval topology on $\Gamma_K$. Note that none of the two above openness notions coincides with being topologically open in the valued field topology.

\begin{remark}
If the value group is discrete, then every proper ideal is $\Gamma_K$-open; otherwise, an ideal is $\Gamma_K$-open if and only if it is not a closed ball. In \cite{Verm:h-min}, the condition of $\Gamma_K$-openness on $I$ in \Cref{thm:all:I} is removed.
\end{remark}

\begin{thm}[$I$-preparation]\label{thm:all:I}
Let $K$ be a valued field of equi-characteristic $0$, considered as a structure in a language $\cL \supset \Lval$. Suppose that $\Th(K)$ is $\omega$-h-minimal.
Then $K$ has resplendent $I$-preparation (see Definition~\ref{defn:Iresplendent}) for every proper $\Gamma_K$-open definable (with parameters) ideal $I$ of $\cO_K$.
\end{thm}

\begin{proof}[Proof of Theorem~\ref{thm:all:I}]
Let $K$ and $I$ be as in the theorem, let
$A \subset K$, fix an $\RV_I$-expansion of $K$ in some language $\cLe \supset\cL$, and let $X \subset K$ be $\cLe(A)$-definable.
We need to find a finite $\cL(A)$-definable set $C\subset K$ such that $X$ is $I$-prepared by $C$. By passing to an elementary extension, we may assume that $K$ is sufficiently saturated as an $\cLe$-structure.

Fix $\lambda \in \Gamma^\times_K$ such that $B_{<\lambda}(0)$ is contained in $I$.
Then up to interdefinability, the $\cL'$-structure on $K$ is an $\RV_{\lambda}$-expansion of the $\cL$-structure on $K$,
namely by the preimage of each of the $\cL'$-definable sets in
$\RV_I^k$ under the map $\RV_{\lambda}^k \to \RV_{I}^k$.

By Lemma~\ref{lem:addconst}, adding a $\lambda$ as a constant to the language preserves $\omega$-h-minimality, i.e., $\cT := \Th_{\cL(\lambda)}(K)$ is still $\omega$-h-minimal. In that language, $B_{<\lambda}(0)$ is $\emptyset$-definable, so Proposition~\ref{prop:equiv} implies that $\cT$ has resplendent $B_{<\lambda}(0)$-preparation. In particular, there exists a finite $\cL(A, \lambda)$-definable set $C$ such that no ball $\lambda$-next to $C$ ``intersects $X$ properly'', i.e., every such ball is either contained in $X$ or disjoint from $X$.
Using Corollary~\ref{cor:finiterange:set}, we may even assume that $C$ is $\cL(A)$-definable.

The above set $C$ might depend on the choice of $\lambda$, so let us denote it by $C_\lambda$ instead. (Contrary to what the notation might suggest, this dependence is not definable, in general.)
Using a similar compactness argument as in Proposition~\ref{prop:uniform},
we now make $C_\lambda$ independent of $\lambda$, as follows:
The condition that $X$ intersects no ball $\lambda$-next to $C_\lambda$ properly can be expressed by a first order formula in $\lambda$, and there
are only a bounded number of choices for $C_\lambda$. It follows that if no finite $\cL(A)$-definable $C$ works for all $\lambda$, then by our saturation assumption on $K$, we find a $\lambda$ such that for every $C$, some ball $\lambda$-next to $C$ intersects $X$ properly, which is a contradiction. Let now $C$ be a set that works for every $\lambda$.

To finish the proof, it remains to prove that every ball $I$-next to $C$ is either contained in $X$ or disjoint from $X$. Let $a, a'\in K$ be such that $\rv_{I}(a - c) = \rv_{I}(a' - c)$ for all $c \in C$. Since $I$ is the union of the open ball ideals $B_{<\lambda}(0)$ it contains, $\frac{a - c}{a'-c} \in 1 + I$ implies $\frac{a - c}{a'-c} \in B_{<\lambda}(1)$ for one such $\lambda$. We may moreover choose a single $\lambda$ such that this holds for all $c \in C$. Then $\bigcap_{c\in C}\rv_{\lambda}^{-1}(\rv_{\lambda}(a-c))$ is a ball $\lambda$-next to $C$ containing both $a$ and $a'$. It follows that $a\in X$ if and only if $a'\in X$.
\end{proof}

From this, we can now deduce that $\omega$-h-minimality is preserved under coarsening of the valuation:

\begin{cor}[Coarsening the valuation]\label{cor:coarse}
Suppose that $\Th_\cL(K)$ is $\omega$-h-minimal and that $|\cdot|_c\colon K \to \Gamma_{K,c}$ is a non-trivial coarsening of the valuation $|\cdot|$ on $K$ with valuation ring $\cO_{K,c}$ (non-trivial meaning $\cO_{K,c} \ne K$).
Let $K_c$ be the expansion of $K$ by a predicate for $\cO_{K,c}$, considered as a  valued field with the valuation $|\cdot|_c$, and denote by
$\cL_c \supset \cL$ the corresponding language.
Then $\Th_{\cL_c}(K_c)$ is also $\omega$-h-minimal.
\end{cor}

\begin{proof}
Since $\cO_{K,c}$ is an $\rv$-pullback (i.e., a preimage in $K$ of a subset of $\RV$), Theorem~\ref{thm:resp:h} implies that $\Th_{\cL_c}(K)$ is $\omega$-h-minimal for the valuation $|\cdot|$.

Let $K'$ be $\cL_c$-elementary equivalent to $K$; we need to verify that $K'$ has $B_{<\lambda}(0)$-preparation for every $\lambda \le 1$ in $\Gamma^\times_{K',c}$.
This follows from Theorem~\ref{thm:all:I} since $B_{<\lambda}(0)$ is $\Gamma_{K'}$-open, which in turn holds because $B_{<\lambda}(0)$ is the preimage under $K' \to \Gamma_{K'}$ of the set of those $\mu \in \Gamma_{K'}$
which get sent to an element less than $\lambda$ by the map
$\Gamma_{K'} \to \Gamma_{K', c}$.
\end{proof}

\begin{remark}\label{coarsening1-h}
In \cite{CHRV}, it is shown that also $1$-h-minimality is preserved under coarsening of the valuation, and, in \cite{Verm:h-min}, this is further generalized to all $\ell\ge 1$.
\end{remark}

\subsection{Algebraic Skolem functions}
\label{sec:alg:skol}

As usual, $K$ is a valued field of equi-characteristic $0$, in a language $\cL \supset \Lval$.

The mere statement of cell decomposition in valued fields is usually very technical, one reason being that one cannot definably pick centers of the cells.
This problem is easier to deal with if one expands the structure by certain Skolem functions. Usually, one would not want to modify the structure in such a way. However, in this subsection, we provide tools which make it possible, in many situations, to assume the existence of such Skolem functions (without losing power nor generality).

The first thing to note is that properly adding the required Skolem functions preserves Hensel minimality; this follows from Theorem~\ref{thm:resp:h} and the observation that adding ``algebraic'' Skolem functions on $\RV$ is enough (Lemma~\ref{lem:alg:skol:K}). This by itself is useful whenever one only wants to prove that every definable set (or function) has some good properties. A similar approach has been recently followed in \cite{CFL}.

Sometimes, one wants to use cell decomposition to prove that every definable object yields some other kind of definable object. If one proves such a result in a structure expanded by Skolem functions, one needs to be able to get rid of the Skolem functions again afterwards. There is probably no general recipe for this, but Lemma~\ref{lem:undo-K-to-RV} is a useful tool; we apply it e.g.\ in the proof of Theorem~\ref{thm:T3/2.mv}.

\begin{lem}\label{lem:alg:skol:K}
Suppose that $\Th(K)$ is $0$-h-minimal and
that for every set $A' \subset \RV$, we have $\acl_{\RV}(A') = \dcl_{\RV}(A')$.
Then for every set $A \subset K$, we have $\acl_K(A) = \dcl_K(A)$.
\end{lem}

\begin{proof}
Let $C \subset K$ be a finite $A$-definable set; we need to prove that $C$ is contained in $\dcl_K(A)$.
By Lemma~\ref{lem:average}, there exists an $A$-definable bijection $f\colon C \to C' \subset \RV^k$.
By Proposition~\ref{prop:stab}, $C'$ is definable with parameters from $\dcl_{\RV}(A)$, so
$C' \subset \acl_{\RV}(\dcl_{\RV}(A))$. By assumption, this implies $C' \subset \dcl_{\RV}(A)$, which,
using $f$, implies $C \subset \dcl_K(A)$.
\end{proof}

\begin{remark}\label{rem:alg:skol}
The above property ``$\acl_K(A) = \dcl_K(A)$'' holds in all models if and only if algebraic Skolem functions exist, i.e.:
For all integers $n\geq 0$ and $m > 0$ and every $\emptyset$-definable set $X\subset K^{n+1}$ with the property that the coordinate projection map $p\colon X \to K^{n}$ has fibers of cardinality precisely $m$, there is a $\emptyset$-definable function $f:  K^{n}\to K$ whose graph is a subset of $X$.
\end{remark}

\begin{prop}[Obtaining $\acl = \dcl$]\label{prop:exist:alg:skol}
Given a valued field $K$ (in a language $\cL$) whose theory is $\ell$-h-minimal, for $\ell \in \{0, 1, \omega\}$, there exists an
$\RV$-expansion of $K$ in a language $\cLas \supset \cL$ such that $\Th_{\cLas}(K)$ is $\ell$-h-minimal and such that for every model $K'\models \Th_{\cLas}(K)$ and for every subset $A \subset K'$, we have $\acl_{\cLas,K'}(A) = \dcl_{\cLas,K'}(A)$.
\end{prop}

\begin{proof}
Let us be lazy and simply expand $K$ by a predicate for each subset of $\RV_K^n$ (for every $n$); let $\cLas$ be the corresponding language. Then $\Th_{\cLas}(K)$ is still $\ell$-h-minimal by Theorem~\ref{thm:resp:h}. Since in this expansion, we have algebraic Skolem functions in $\RV_K^n$, we in particular have $\acl_{\RV}(A') = \dcl_{\RV}(A')$ for every $A' \subset \RV_{K'}$ and for every model $K' \equiv_{\cLas} K$. Now Lemma~\ref{lem:alg:skol:K} implies
$\acl_K(A) = \dcl_K(A)$ for every $A \subset K'$.
\end{proof}

One may use Proposition \ref{prop:exist:alg:skol} to ensure the condition that $\acl$ equals $\dcl$; the next lemma can be useful to get back to the original language.

\begin{lem}[Undoing algebraic skolemization]\label{lem:undo-K-to-RV}
Suppose that $\Th_{\cL}(K)$ is $0$-h-minimal and
that we have an $\RV$-expansion of $K$ in some language $\cL' \supset\cL$. Let $\chi'\colon K^n \to \RV^{k'}$ be an $\cL'$-definable map (for some $k' \ge 0$).
Then there exists an $\cL$-definable map $\chi\colon K^n \to \RV^k$ (for some $k \ge 0$) such that $\chi'$ factors over $\chi$, i.e.
$\chi' = g\circ \chi$, for some function $g\colon \RV^k \to \RV^{k'}$ (which is automatically $\cL'$-definable).
\end{lem}

\begin{proof}
We do an induction over $n$. The case $n = 0$ is trivial, so assume $n \ge 1$.

Fix $a \in K$ and set $\chi'_a(\tilde a) := \chi'(a, \tilde a)$, for $\tilde a \in K^{n-1}$.
By induction, we find an $\cL(a)$-definable map $\chi_a\colon K^{n-1} \to \RV^k$
such that $\chi'_a = g_a \circ \chi_a$ for some $\cL'(a)$-definable $g_a$.
By compactness, we may assume that $\chi_a$ and $g_a$ are definable uniformly in $a$,
so that in particular, we obtain an $\cL'$-definable set $W = \{(a, \zeta', g_a(\zeta')) \mid a \in K, \zeta' \in \RV^k\}$.
Using Corollary~\ref{cor:prep}, we find a finite $\cL'$-definable set $1$-preparing $W$, and by
Corollary~\ref{cor:acl=acl}, that set is contained in a finite $\cL$-definable set $C$.

Finally, let $f\colon K \to \RV^\ell$ be an $\cL$-definable map as provided by Lemma~\ref{lem:InextFam}, namely such that each nonempty fiber of $f$
is either a singleton in $C$ or a ball $1$-next to $C$. We claim that the map
$\chi(a, \tilde a) := (f(a), \chi_a(\tilde a))$ is as desired.
Indeed, it is clearly $\cL$-definable, and since $C$ $1$-prepares $W$, the function $g_a$ is determined by $f(a)$,
so that $f(a)$ and $\chi_a(\tilde a)$ together determine $\chi'(a, \tilde a) = \chi'_a(\tilde a) = g_a(\chi_a(\tilde a))$.
\end{proof}

\section{Geometry in the Hensel minimal setting}
\label{sec:compactn}

In this section, we deduce geometric results in $K^n$ under the assumption of $1$-h-minimality. Some of those are already known under stronger assumptions,  like $C^k$ properties of definable functions (Subsection \ref{sec:cont}),
cell decomposition (Subsection \ref{sec:cd}), and dimension theory (Subsection \ref{sec:dim}). Since some proofs are very similar to those in many other papers, we will be somewhat succinct.

As highlights, we then prove a version of the Jacobian Property in many variables (Subsection \ref{sec:sjp}) which allows us to get general t-stratifications for definable sets in $1$-h-minimal structures (Subsection \ref{sec:t-strat}).
We also provide a higher-dimension version of the Taylor-approximation result (Subsection \ref{sec:taylor-box}), the aim being to lay the
ground for the first axiomatic approach in the non-archimedean context of analogues of results by Yomdin--Gromov and Pila--Wilkie on parameterizations and point counting.

The version of cell decomposition presented in Subsection \ref{sec:cd} uses a simplified notion of cells, by temporarily adding certain Skolem functions to the language, as detailed in Subsection~\ref{sec:alg:skol}. Since that approach might not work for all potential applications,
we conclude Section~\ref{sec:compactn} by linking with more classical viewpoints on cells (Subsection \ref{sec:cd:classical}).

\subsection{Continuity and Differentiability}\label{sec:cont}

We start by proving that $1$-h-minimality implies that definable
functions are almost everywhere continuous and even $C^k$. Note that all this follows directly from Theorem~\ref{thm:t-high-high}; moreover, the notion of ``almost everywhere'' makes more sense once dimension theory is developed (see Subsection~\ref{sec:dim}). For those reasons, for now the statements will be a bit weaker and the proofs a bit sketchy.

\begin{thm}[Continuity]\label{thm:C0}
Assume that $\Th(K)$ is $1$-h-minimal. For every definable function $f\colon X \subset K^n \to K$, the set $U$ of those $u \in X$ such that $f$ is continuous on a neighborhood of $u$ is dense in $X$.
\end{thm}

For the moment, we only give the proof in the case where $X$ is open. The general case is proved
in a joint induction with cell decomposition (more precisely with Addenda~\ref{add:cd:cont:f} and \ref{add:cd:cont:c} of Theorem~\ref{thm:cd:alg:skol}),
so we postpone that proof until Subsection~\ref{sec:cd}.

The proof given here uses a different approach than the one used for similar results in \cite{CCL-PW}, building on the following
two lemmas which might be of independent interest.

\begin{lem}[Continuity on small boxes]\label{lem:const:type:cont}
Assume that $\Th(K)$ is $1$-h-minimal. If $B \subset K^n$ is a product of balls such that all elements of $B$ have
the same type over $\RV$ and $f\colon B \to K$ is a $\emptyset$-definable function, then $f$ is uniformly continuous on $B$.
\end{lem}

\begin{proof}
Set $B =: \hat B \times B' \subset K^{n-1} \times K$.
Using Lemma~\ref{lem:gammaLin}, we deduce that for every fixed $\hat a \in \hat B$, the map $a' \mapsto f(\hat a, a')$ is continuous on $B'$.
Moreover, this continuity is uniform in both, $\hat a$ and $a'$, since all elements $(\hat a, a') \in B$ have the same type over $\RV$.
(The type $\tp(\hat a a'/\RV)$ knows which $\delta$ works for which $\epsilon$, where $\epsilon, \delta>0$ are from the definition of continuity at $(\hat a, a')$.)
After applying the same argument to all other coordinates, we deduce that $f$ as a whole is uniformly continuous in $B$:
For every $\epsilon > 0$, there exists a $\delta > 0$ such that if $a_1, a_2 \in B$ with $|a_1 - a_2| < \delta$ differ only in one coordinate,
then $|f(a_1) - f(a_2)| < \epsilon$. If $a_1$ and $a_2$ differ in several coordinates, apply this repeatedly.
\end{proof}

\begin{lem}[Balls with constant type]\label{lem:ex:const:type}
Assume that $K$ is $|\cL|^+$-saturated and that $\acl_{K}$ equals $\dcl_{K}$.
If $X \subset K^n$ is a $\emptyset$-definable set with non-empty interior,
then $X$ contains a ball $B$ such that all elements of $B$ have
the same type over $\RV$ (i.e., $\tp(a/\RV) = \tp(a'/\RV)$ for all $a, a' \in B$).
\end{lem}

\begin{proof}[Proof of Lemma~\ref{lem:ex:const:type}]
We may assume that $X$ is of the form $\hat B \times B'$ for some balls $\hat B \subset K^{n-1}, B' \subset K$, and by induction, we may assume that all elements of $\hat B$ have the same type over $\RV$.

Since any $\emptyset$-definable function $f\colon K^{n-1} \to K$ is continuous on $\hat B$ (by Lemma~\ref{lem:const:type:cont}), given such an $f$, we can first shrink $\hat B$
so that $B' \setminus f(\hat B)$ still contains a ball and then shrink $B'$ so that
$\hat B \times B'$ is disjoint from the graph of $f$. After possibly further shrinking $\hat B$, we then obtain
\begin{equation}\label{eq:rveq}
\rv(a'_1 - f(\hat a_1)) = \rv(a'_2 - f(\hat a_2)).
\end{equation}
for any $(\hat a_i, a'_i) \in \hat B \times B'$. By saturation, this shrinking of $\hat B$ can be done for all $\emptyset$-definable $f$ simultaneously.
In particular, $B'$ is now disjoint from $\acl_K(\hat a_1)$ for every $\hat a_1 \in \hat B$, and
one deduces that all elements of $\hat B \times B'$ have
the same type over $\RV$, namely as follows:
Given $(\hat a_i, a'_i) \in \hat B \times B'$ for $i=1,2$, we find some $a''_1 \in B'$ such that $\tp(\hat a_1a''_1 / \RV) = \tp(\hat a_2a'_2 / \RV)$. This implies $\rv(a''_1 - f(\hat a_1)) = \rv(a'_2 - f(\hat a_2))$ for every $\emptyset$-definable $f$ and hence (using (\ref{eq:rveq}) and Lemma~\ref{lem:type-0-h-min}) $\tp(\hat a_1a''_1 / \RV) = \tp(\hat a_1a'_1 / \RV)$.
\end{proof}

\begin{proof}[Proof of Theorem~\ref{thm:C0} when $X$ is open]
We may suppose that $f$ is $\emptyset$-definable, that $K$ is $|\cL|^+$-saturated and that $\acl_K$ equals $\dcl_K$ (by Proposition~\ref{prop:exist:alg:skol}).
Suppose that the set $U$ from the theorem is not dense in $X$. Choose a ball $B \subset X \setminus U$ as provided by Lemma~\ref{lem:ex:const:type}. By Lemma~\ref{lem:const:type:cont}, $f$ is continuous on $B$, which is a contradiction to $B \not\subset U$.
\end{proof}

We now come to differentiability of definable functions.
The definition of $C^k$ is the usual one:

\begin{defn}[$C^k$]\label{defn:deriv}
Let $U \subset K^m$ be open and let $f\colon U \to K^n$ be a map.
We say that $f$ is $C^0$ if it is continuous, and that it is $C^k$ for $k \ge 1$ if there exists a $C^{k-1}$-map $\Jac f\colon U \to K^{n\times m}$
(the Jacobian of $f$) such that
for every $u \in U$, we have $\lim_{x \to u} \frac{|f(x) - f(u)-((\Jac f)(u))(x - u)|}{|x - u|} = 0$.
In the case $n = 1$, we also write $\grad f$ instead of $\Jac f$ and call it the gradient of $f$.
\end{defn}

\begin{thm}[$C^k$]\label{thm:Ck}
Assume that $\Th(K)$ is $1$-h-minimal and fix $k \ge 0$.
For every definable function $f\colon K^n \to K$, the set $U$ of those $u \in K^n$ such that $f$ is $C^k$ on a neighborhood of $u$ is dense in $K^n$.
\end{thm}

\begin{remark}\label{rem:Ck}
Clearly, $U$ is open and definable over the same parameters as $f$,
so $f$ is $C^k$ on a definable dense open subset of $K^n$.
\end{remark}

\begin{remark}
In Subsection~\ref{sec:dim}, we will introduce a notion of dimension, and we will see that $U$ being dense implies that $K^n \setminus U$ has dimension less than $n$. Thus,  definable functions are almost everywhere $C^k$ in this rather strong sense.
\end{remark}

\begin{proof}[Proof of Theorem~\ref{thm:Ck}]\footnote{We would like to thank Yatir Halevi for pointing out that a previous argument for \Cref{thm:Ck} was wrong.}
The case $k = 0$ is Theorem~\ref{thm:C0}. The cases $k \ge 2$ follow from the case $k = 1$ by induction, so now assume $k = 1$. By using the same strategy as in the proof of Theorem~\ref{thm:C0}, it suffices to prove the following: If $B \subset K^n$ is a product of balls such that all elements of $B$ have the same type over $\RV$ and $f\colon B \to K$ is a $\emptyset$-definable function, then $f$ is $C^1$ on $B$.

So let such a $B$ be given. To prove that $f$ is $C^1$ on $B$, we use three ingredients:

Firstly, we use that all second partial derivatives of $f$ have constant valuation on $B$, say, bounded by $\lambda$.
Secondly, for each $i \le n$, we apply Theorem~\ref{thm:high-ord} with $r = 1$ to the function $f$ considered as a function in the $i$-th coordinate and with all other coordinates fixed. Using that all elements of $B$ have the same type over $\RV$, we deduce that for $x, x' \in B$ differing only in one coordinate, we have
\begin{equation}\label{eq:proof:Ck}
|f(x) - f(x') - \nabla f(x')\cdot (x - x')| \le \lambda|x - x'|^2
\end{equation}
Thirdly, we apply Theorem~\ref{thm:high-ord} with $r = 0$ in a similar way to each of the first partial derivatives of $f$.

Now the same computation as in the proof of Theorem~\ref{thm:t-high-high} implies that $f$ is $C^1$ on $B$. The idea of this computation is the following: Given arbitrary $x$ and $x'$ in $B$, we change one coordinate at a time, apply \eqref{eq:proof:Ck} in each step, and use
the third ingredient to bound the difference between the value of $\nabla f$ at the starting point and its value at the intermediate points. In the end, we obtain the same bound as in \eqref{eq:proof:Ck} for arbitrary $x, x' \in B$.
\end{proof}

\subsection{Cell decomposition}\label{sec:cd}

In this subsection we present a version of cell decomposition results that are simpler than usual, by imposing the following condition.

\begin{assumption}\label{ass:acl=dcl}
We assume in this subsection that we have algebraic Skolem functions in $K$; or, equivalently, that for every $K' \equiv K$ and every $A \subset K'$, we have $\acl_{K'}(A) = \dcl_{K'}(A)$. We will abbreviate this by ``$\acl$ equals $\dcl$ (in $\Th(K)$)''.
\end{assumption}

The tools provided in Subsection~\ref{sec:alg:skol} often make it possible to reduce to the case of languages where this assumption holds, as we will see
in later subsections of Section~\ref{sec:compactn}.
In fact, Assumption~\ref{ass:acl=dcl} does more than just simplify the arguments: it also allows one to formulate stronger results like piecewise Lipschitz continuity results as in Theorem \ref{thm:cd:alg:piece:Lipschitz}.
This is similar to \cite{CFL}, where this condition on $\acl$ is furthermore used to obtain parametrization results with finitely many maps (in the Pila-Wilkie sense \cite{PW}), but not under axiomatic assumptions.

Close variants of the results of this subsection appear in \cite{CFL}; indeed, the tameness condition of \cite{CFL} is very similar to the
$1$-h-minimality criteria from Theorem \ref{thm:tame2vf}; see Remark \ref{rem::tame2vf}

In the following, for $m \le n$, we denote the projection $K^n \to K^{m}$ to the first $m$ coordinates by $\pi_{\le m}$, or also by $\pi_{<m+1}$.

\begin{defn}[Cells, twisted boxes]\label{defn:cell}
Fix any parameter set $A \subset K\eq$ and consider a non-empty $A$-definable set $X\subset  K^n$ for some $n$, and, for $i=1,\ldots, n$, values $j_i$ in $\{0,1\}$ and $A$-definable functions $c_i:\pi_{<i}(X)\to K$.
Then $X$ is called an \emph{$A$-definable cell} with \emph{center tuple $c=(c_i)_{i=1}^n$} and of \emph{cell-type $j= (j_i)_{i=1}^n$} if it is of the form
$$
X = \{x\in K^n\mid  (\rv(x_i-c_i(x_{<i})))_{i=1}^n
   \in R  \},
$$
for a (necessarily $A$-definable) set
$$
R \subset  \prod_{i=1}^n (j_i\cdot \RV^\times),
$$
where $x_{<i}=\pi_{<i}(x)$ and where
$0 \cdot \RV^\times = \{0\} \subset \RV$, and $1 \cdot \RV^\times = \RV^\times \subset \RV$.

If $X$ is such a cell, then, for any $r\in R$, the subset
$$
\{x \in K^n\mid  \rv(x_i-c_i(x_{<i})))_{i=1}^n   = r  \},
$$
of $X$ is called a \emph{twisted box} of the cell $X$.
We also call $X$ itself a twisted box if it is a cell consisting
of a single twisted box (i.e., if $R$ is a singleton).
\end{defn}

\begin{remark}\label{rem:cell:homeo}
Given a cell $X$ as above, let $\pi\colon K^n \to K^d$ be the projection to those coordinates $i$ for which $j_i = 1$.
Then $\pi(X)$ is a cell of cell-type $(1, \dots, 1)$, and the restriction of $\pi$ to $X$ is injective. If the components $c_i$ of the center tuple of $X$ are continuous, then $\pi|_X$ is even a homeomorphism.
\end{remark}

There are many variants of cell decomposition results preparing different kinds of data. We first state the simplest version
and then formulate other versions as addenda:

\begin{thm}[Cell decomposition]\label{thm:cd:alg:skol}
Suppose that $\acl$ equals $\dcl$ in $\Th(K)$ (see Assumption~\ref{ass:acl=dcl}), and that $\Th(K)$ is $1$-h-minimal.
Consider a $\emptyset$-definable set $X\subset  K^n$ for some $n$.
Then there exists a cell decomposition of $X$, more precisely, a partition of $X$ into finitely many $\emptyset$-definable cells $A_\ell$.
\end{thm}

\begin{remark}\label{rem:cd:alg:skol}
Given finitely many $\emptyset$-definable sets $X_i \subset K^n$, we can find a cell decomposition of $K^n$ such that
each $X_i$ is a union of cells, namely by applying Theorem~\ref{thm:cd:alg:skol} to suitable intersections of the sets $X_i$ and $K^n \setminus X_i$.
\end{remark}

\begin{addendum}[Preparation of $\RV$-sets]\label{add:cd:alg:prep}
On top of the assumptions from Theorem \ref{thm:cd:alg:skol}, let also a $\emptyset$-definable set $P\subset X\times \RV^k$ be given for some $k$. We consider $P$ as the function sending $x\in X$ to the fiber $P_x := \{\xi \in \RV^k \mid (x,\xi) \in P\}$.

Then the cells $A_\ell$ from Theorem \ref{thm:cd:alg:skol} can be taken such that moreover $P$ (seen as function)
is constant on each twisted box of each cell $A_\ell$.
\end{addendum}

\begin{addendum}[Continuous functions]\label{add:cd:cont:f}
On top of the assumptions from Theorem \ref{thm:cd:alg:skol} (with Addendum~\ref{add:cd:alg:prep}, if desired), suppose that finitely many $\emptyset$-definable functions $f_j\colon X \to K$ are given. Then the $A_\ell$ can be taken such that moreover, the restriction $f_j|_{A_\ell}$ of each function $f_j$ to each cell $A_\ell$ is continuous.
\end{addendum}

In the one-dimensional case, we can prepare the domain and the image of the functions in a compatible way:

\begin{addendum}[Compatible preparation of domain and image]\label{add:cd:alg:range}
Under the assumptions of Addendum~\ref{add:cd:cont:f}, if the ambient dimension $n$ is equal to $1$,
we may moreover impose that for each $\ell$ and each $j$, $f_j|_{A_\ell}$ is either constant or injective, $f_j(A_\ell)$ is a $\emptyset$-definable cell and $f_{j}$ maps each twisted box of $A_\ell$ onto a twisted box of $f_j(A_\ell)$.
\end{addendum}

\begin{addendum}[Continuous centers]\label{add:cd:cont:c}
In Theorem~\ref{thm:cd:alg:skol} (with any of Addenda~\ref{add:cd:alg:prep}, \ref{add:cd:cont:f}, if desired), we may assume that for each cell $A_\ell$, each component $c_i\colon \pi_{<i}(A_\ell) \to K$ of its center tuple is continuous. In particular, each twisted box of each cell $A_\ell$ of cell-type $(1, \dots, 1)$ is an open subset of $K^n$.
\end{addendum}

\begin{remark}
In Addenda~\ref{add:cd:cont:f} and \ref{add:cd:cont:c}, one can replace continuous by $C^k$, provided that one fixes a reasonable definition of a function being $C^k$ on non-interior points of its domain.
\end{remark}

To formulate the next addendum, we need the notion of Lipschitz continuity, which is defined in the usual way:

\begin{defn}[Lipschitz continuity]\label{defn:Lipschitz}
For a valued field $K$ and an element $\lambda$ in its value group $\Gamma_K^\times$, a function $f:X\subset K^n\to K^m$ is called Lipschitz continuous with Lipschitz constant $\lambda$ if for all $x$ and $x'$ in $X$ one has
$$
|f(x)-f(x')| \leq \lambda |x-x'|,
$$
where the norm of tuples is, as usual, the sup-norm.
We call such $f$ shortly $\lambda$-Lipschitz.

Call $f$ locally $\lambda$-Lipschitz, if for each $x\in X$ there is an open neighborhood $U$ of $x$  such that the restriction of $f$ to $U$ is $\lambda$-Lipschitz.
\end{defn}

\begin{addendum}[$1$-Lipschitz centers]
\label{add:cd:Lip:comp}
Theorem~\ref{thm:cd:alg:skol} (with any of Addenda~\ref{add:cd:alg:prep}, \ref{add:cd:cont:f}, \ref{add:cd:cont:c}, if desired),
is also valid in the following variant:
Instead of imposing that $A_\ell$ itself is a cell in the sense of Definition~\ref{defn:cell}, we only impose that $\sigma_{\ell}(A_\ell)$ is a cell, for some coordinate permutation $\sigma_\ell\colon K^n \to K^n$.
In that version, we may moreover impose that $\sigma_{\ell}(A_\ell)$ is of cell-type $(1, \dots, 1, 0, \dots, 0)$ and that each component $c_i$ of the center tuple $(c_i)_i$ of $\sigma_{\ell}(A_\ell)$ is $1$-Lipschitz.
\end{addendum}

Closely related to that addendum, we also have the following reformulation of the piecewise continuity result of \cite{CFL} in the Hensel minimal setting.

\begin{thm}[Piecewise Lipschitz continuity]\label{thm:cd:alg:piece:Lipschitz}
Suppose that $\acl$ equals $\dcl$ in $\Th(K)$ (see Assumption~\ref{ass:acl=dcl}), and that $\Th(K)$ is $1$-h-minimal.
Consider a $\emptyset$-definable set  $X\subset  K^n$ for some $n$ and a $\emptyset$-definable function $f:X\to K$. Suppose that $f$ is locally $1$-Lipschitz.
Then there exists a finite partition of $X$ into $\emptyset$-definable sets $A_\ell$ such that the restriction of $f$ to $A_\ell$ is $1$-Lipschitz, for each $\ell$.
\end{thm}

All of the above results in this subsection have already been proved under related assumptions. For the convenience of the reader, we nevertheless provide short, self-contained versions of the proofs, with exception of Addendum~\ref{add:cd:Lip:comp} and Theorem~\ref{thm:cd:alg:piece:Lipschitz} which will not be used in this paper and where we will only explain how to adapt the proof from \cite{CFL}. However, we do provide a proof
of the weak version of Addendum~\ref{add:cd:Lip:comp} stated as Proposition~\ref{prop:twibox:1Lip}, since that proposition is used to deduce the existence of t-stratifications.

\begin{proof}[Proof of Theorem \ref{thm:cd:alg:skol} with Addendum \ref{add:cd:alg:prep}]
It suffices to find a cell decomposition adapted to a set $P \subset K^n \times \RV^k$.

Case $n = 1$:
In this case,
Corollary~\ref{cor:prep} provides
a cell decomposition, namely: Let $C \subset K$ $1$-prepare $P$. Using the assumption that $\acl$ equals $\dcl$, each element $c \in C$ is $\emptyset$-definable, so we obtain a cell decomposition as follows: For each ball $B$ $1$-next to $C$,
choose (in a definable way) a $c(B) \in C$ such that $B$ is $1$-next to $c(B)$.
Now the cell decomposition consists of two cells with center $c$ for each $c \in C$, namely (i) $\{c\}$ and (ii) the union of all those $B$ $1$-next to $C$ for which $c(B)$ equals $c$.

Case $n > 1$:
We apply the case $n = 1$ to each fiber
$P_a = \{(b,\xi) \in K \times \RV^k \mid (a,b,\xi) \in P\}$, where $a$ runs over $K^{n-1}$. (Note that by compactness this works uniformly, so that in particular we get definable cell centers $K^{n-1} \to K$.)
Then we finish by applying induction to a set $P' \subset K^{n-1} \times \RV^{k'}$ ``describing'' the fibers:
For each $0$-cell $\{c\} \subset K$ of the fiber at $a \in K^{n-1}$,
$P'_a$ encodes the set $P_{a,c} \subset \RV^k$;
for each $1$-cell $X \subset K$ of the fiber at $a \in K^{n-1}$,
$P'_a$ encodes (a) the set denoted by $R$ in Definition~\ref{defn:cell} and (b), for each $\xi \in R$, the fiber $P_{a,b} \subset \RV^k$,
where $b \in K$ is an arbitrary element of the twisted box corresponding to $\xi$ (i.e., $\rv(b - c) = \xi$, where $c$ is the center of $X$).
\end{proof}

\begin{proof}[Proof of Theorem~\ref{thm:C0} and Addenda~\ref{add:cd:cont:f}, \ref{add:cd:cont:c}]
For $n = 0$, all three results are trivial. We now assume that all three results are already known for $n - 1$ and we deduce them for $n$. Concerning Theorem~\ref{thm:C0}, note that we may as well
assume that $\acl$ equals $\dcl$ (by Proposition~\ref{prop:exist:alg:skol}).

Addendum~\ref{add:cd:cont:c}: First, find a cell decomposition with possibly non-continuous centers. By inductively applying Addendum~\ref{add:cd:cont:f} to each component of the center tuple of each cell, we may refine the cell decomposition to get continuous centers.

Theorem~\ref{thm:C0}: We may suppose that $X$ has empty interior, since the proof given at the end of Subsection~\ref{sec:cont} applies to the interior of $X$.
Choose a cell decomposition of $X$ with continuous centers. No cell $A_\ell \subset X$ is of cell-type $(1, \dots, 1)$ (since such cells have non-empty interior,
by the ``in particular'' part of Addendum~\ref{add:cd:cont:c}). Thus the homeomorphism from Remark~\ref{rem:cell:homeo} allows us to reduce the problem on $A_\ell$
to one of lower ambient dimension. Apply induction.

Addendum~\ref{add:cd:cont:f}: The above proof of Theorem~\ref{thm:C0} also yields a finite partition of $X$ such that $f$ is continuous on each piece.
Apply this to each of the given functions $f_j$ and then choose a cell decomposition respecting all pieces from all those partitions.
\end{proof}

\begin{proof}[Proof of Addendum~\ref{add:cd:alg:range}]
Using Lemma~\ref{lem:fin-inf} and our assumption $\acl=\dcl$,
we find a partition of $X$ such that on each piece, $f_j$ is continuous and either constant or injective for each $j$; assume without loss that $X$ is a single such piece. In a similar way, assume without loss that $X$ is a cell and that each $f_j$ has the Jacobian Property (Definition~\ref{defn:JP}) on each twisted box of $X$. Since constant functions pose no problem, we assume that all $f_j$ are injective.

Next, choose a finite $\emptyset$-definable set $\tilde C \subset K$ $1$-preparing $f_j(B)$ for every $j$ and every twisted box $B$ of $X$ (using Corollary~\ref{cor:prep}) and set $C := \bigcup_j f^{-1}_j(\tilde C)$. After a further finite partition of $X$, we may assume that either (i) $X$ consists of a single twisted box or that (ii) $C$ is empty.

In Case (i), choose any $c \in C \subset X$ and decompose $X$ as \(\{c\}\cup (X\setminus \{c\})\), both of which are \(\emptyset\)-definable cells. Then the desired properties follow from the Jacobian Property, namely the image of both cells are cells with center $f_j(c)$.

In Case (ii), definably choose, for each twisted box $B$ of $X$ and each $j$, an element $\tilde c_{j,B} \in \tilde C$ in such a way that $f_j(B)$ is $1$-next to $\tilde c_{j,B}$. By a further finite partition of $X$, we may assume that $\tilde c_{j,B}$ does not depend on $B$. Then
$f_j(X)$ is a cell with center $\tilde c_{j,B}$ and we are done.
\end{proof}

As explained above, we will not give detailed proofs of Theorem \ref{thm:cd:alg:piece:Lipschitz} and its related result from Addendum \ref{add:cd:Lip:comp} to Theorem \ref{thm:cd:alg:skol}, and we do not use these results in this paper. We nevertheless specify where this is worked out, under very closely related assumptions.
\begin{proof}[Proof of Theorem \ref{thm:cd:alg:piece:Lipschitz} and Addendum \ref{add:cd:Lip:comp} to Theorem \ref{thm:cd:alg:skol}]
Under our Assumption \ref{ass:acl=dcl}  that $\acl$ equals $\dcl$, but assuming a notion of tameness with an angular component map $\ac$ (instead of $1$-h-minimality with $\rv$), both results are proved in \cite{CFL}, and the proof readily adapts. (By Theorem \ref{thm:tame2vf} and Remark \ref{rem::tame2vf}, $1$-h-minimality and the tameness notion are very closely related.)
\end{proof}

(The proof of Addendum~\ref{add:cd:Lip:comp} and its variant in \cite{CFL} essentially comes from \cite{CCL-PW}, apart from the improvement made possible by the assumption $\acl_K=\dcl_K$.)

\subsection{Dimension theory}\label{sec:dim}

Under the assumption of $1$-h-minimality, there is a good notion of dimension of definable subsets of $K^n$. It can be defined in various equivalent ways; here is one possible definition.

\begin{defn}[Dimension]\label{defn:dim}
We define the dimension of a non-empty definable set $X \subset K^n$ as the maximal integer $m$ such that there is a $K$-linear function $\ell : K^n\to K^m$ such that $\ell(X)$ has non-empty interior in $K^m$. If $X$ is empty, we set $\dim X := -\infty$.
\end{defn}

\begin{remark}
In Proposition~\ref{prop:dim:basic} (\ref{dim:1}), we will see that one could equivalently only consider coordinate projections $\ell\colon K^n \to K^m$ (instead of arbitrary linear maps).
\end{remark}

In \cite{CLb}, there already exists a proof that dimension behaves well if one imposes a suitable condition on $\Th(K)$ called $b$-minimality. By the following proposition, $1$-h-minimality implies $b$-minimality, so for several results about dimension, we could simply refer to \cite{CLb}. For the convenience of the reader, we will nevertheless give a self-contained proof of those results.

\begin{prop}[b-minimality]\label{prop:b-min}
Assume that $\Th(K)$ is $1$-h-minimal.
Then the two sorted structure on $(K,\RV)$ obtained from $K$ by adding the sort $\RV$ and the map $\rv$ is $b$-minimal in the sense of Definition 2.1 of \cite{CLb}, with $K$ as main sort. More specifically, the structure $(K,\RV)$ is $b$-minimal with centers and preserves all balls in the sense of Definitions 5.1 and 6.2 of \cite{CLb}.
\end{prop}
\begin{proof}
The axioms of Definition 2.1 of \cite{CLb} clearly hold, and Definition 5.1 of \cite{CLb} follows from the Jacobian Property as formulated in Corollary \ref{cor:JP}.
\end{proof}

The definition of dimension given in \cite[Definition~4.1]{CLb} is different than ours, but the results from \cite[Section 4]{CLb} imply that the definitions are equivalent:
If $X \subset K^n$ is a finite union of cells, then the dimension of $X$ in our sense equals the dimension of $X$ in the sense of \cite{CLb}, namely the maximum of the dimensions of the cells, where the
dimension of a cell of cell-type $(j_i)_{i=1}^n$
is $\sum_i j_i$.

The following proposition summarizes the good properties of dimension; in particular, we have definability of dimension, as in o-minimal structures. Property (\ref{prop:dim:frontier}) is new.

\begin{prop}[Dimension theory]\label{prop:dim:basic}
Assume that $\Th(K)$ is $1$-h-minimal. Let $X\subset K^n$, $Y\subset K^n$ and $Z\subset K^{m}$ be non-empty definable sets, and let $f:X\to Z$ be a definable function. Then the following properties hold.
\begin{enumerate}
 \item\label{dim:1}
 For any $d \le n$, we have $\dim X \ge d$ if and only if there exists a projection $\pi\colon K^n \to K^d$ to a subset of the coordinates
 such that $\pi(X)$ has non-empty interior. In particular,
 $\dim X = 0$ if and only if $X$ is finite.
 \item\label{dim:2} $\dim(X \cup Y) = \max\{\dim X, \dim Y\}$.
 \item\label{dim:3}\label{dim:defble} For any $d \le n$, the set of $z\in Z$ such that $\dim f^{-1}(z) = d$ is definable over the same parameters as $f$.
 \item\label{dim:4} If all fibers of $f$ have dimension $d$, then $\dim X = d + \dim Z$.
 \item\label{dim:5}\label{prop:dim:local} There exists an $x \in X$ such that the local dimension of $X$ at $x$ is equal to the dimension of $X$, i.e.,
   such that for every open ball $B \subset K^n$ around $x$, we have $\dim (X \cap B) = \dim X$.
 \item\label{dim:6}\label{prop:dim:frontier} One has $\dim (\overline X\setminus X) < \dim X$, where $\overline X$ is the topological closure of $X$, for the valuation topology.
\end{enumerate}
\end{prop}


Although most likely, property (\ref{prop:dim:frontier}) can be proved in a similar way as Theorem (1.8) of \cite{vdD}, we postpone that proof until the end of Subsection~\ref{sec:t-strat}, where we have t-stratifications at our disposal, which make the proof much simpler.
We do however right away prove the ``easy'' case of Property (\ref{prop:dim:frontier}), namely when $\dim X$ is equal to the ambient dimension $n$.

As announced, we present a proof of the proposition which does not rely on $b$-minimality. Readers willing to take $b$-minimality as a blackbox can jump directly to the paragraph ``Property (\ref{prop:dim:local})\dots'' (on p.~\pageref{jump:here}) in the proof of the proposition. Indeed, (\ref{dim:1})--(\ref{dim:4}) are proved in \cite[Section 4]{CLb}, except for the ``in particular'' part of (\ref{dim:1}), which is easy to deduce from the fact that infinite definable subsets of $K$ have non-empty interior.

Our own proof of (\ref{dim:1})--(\ref{dim:4}) consists in showing that dimension agrees with the acl-dimension (which we quickly recall below). In particular, the key is to show that $\acl$ has the exchange property (so that acl-dimension makes sense).

\begin{lem}
If $\Th(K)$ is $1$-h-minimal, then $\acl_K$ has the exchange property, i.e., for any $a \in K^n$ and any $b,c \in K$ satisfying $c \in \acl_K(a, b) \setminus \acl_K(a)$, we have $b \in \acl_K(a, c)$.
\end{lem}

\begin{proof}
Let $a,b,c$ be given as in the lemma; we need to show that $b \in \acl_K(a, c)$. We may assume that $a$ is the empty tuple, otherwise adding it to the language. Since $\RV$-expansions do not change $\acl_K$ (by Corollary~\ref{cor:acl=acl}), we may moreover assume that $\acl$ equals $\dcl$ (using Proposition~\ref{prop:exist:alg:skol}). Thus our assumptions yield $c \in \dcl_K(b)$, which implies that there exists a $\emptyset$-definable function $f\colon K \to K$ satisfying $f(b) = c$.
If $f^{-1}(c)$ is finite, we deduce that $b \in \acl_K(c)$, and we are done. Otherwise, $c$ lies in the ($\emptyset$-definable) set $Y$ of those $y \in K$ for which $f^{-1}(y)$ is infinite. By Lemma~\ref{lem:fin-inf}, $Y$ is finite, so this contradicts $c \notin \acl_K(\emptyset)$.
\end{proof}

Since $\acl$ has the exchange property, it defines a pre-geometry and yields a notion $\dim^\acl(X)$ of dimension of a definable set $X \subset K^n$; we quickly recall how it is defined.
Firstly, note that we may assume that $K$ is $\aleph_0$-saturated (otherwise replacing it by an elementary extension).
One calls a tuple $b = (b_1, \dots, b_m) \in K^m$ algebraically independent over a set $A \subset K$, if $b_i \notin \acl_K(A, b_1, \dots, b_{i-1})$ for every $i \le m$. Now, given an $A$-definable set $X \subset K^n$, for some finite $A \subset K$, we define $\dim^\acl(X)$ as the maximal $m$ for which there exists a coordinate projection $\pi\colon K^n \to K^m$ such that $\pi(X)$ contains a tuple which is algebraically independent over $A$. The theory of pre-geometries shows that $\dim^\acl$ is well-defined (i.e., does not depend on the specific finite set $A$) and well-behaved; in particular, it satisfies Properties (\ref{dim:2}) and (\ref{dim:4}) of Proposition~\ref{prop:dim:basic}, and using
$\exists^\infty$-elimination (Lemma~\ref{lem:finite}), it also satisfies (\ref{dim:3}).

To see that our dimension agrees with $\dim^\acl$, we need the following lemma:

\begin{lem}\label{lem:int}
Suppose that $K$ is $\aleph_0$-saturated and that $A \subset K$ is finite. An $A$-definable set $X \subset K^m$ has non-empty interior if and only if $\dim^\acl X = m$, i.e., if and only if it contains a tuple which is algebraically independent over $A$. More precisely, if $a \in X$ is algebraically independent over $A$, then $a$ lies in the interior of $X$.
\end{lem}

\begin{proof}
For the implication from left to right, suppose without loss that $X$ is a ball. Then we find the desired algebraically independent tuple $a \in X$ by choosing one coordinate after the other, each one outside of the algebraic closure of $A$ and the previously chosen coordinates; this is possible by $\aleph_0$-saturation.

For the other direction (including the ``more precisely'' part), suppose that $a \in X$ is algebraically independent over $A$ and write it as $a = (\hat a, a') \in K^{m-1} \times K$.
On the one hand,
$\hat a$ is algebraically independent over $A \cup \{a'\}$, so by induction, it lies in the interior of the fiber $X_{a'} \subset K^{m-1}$. Fix $\lambda \in \Gamma_K^\times$ such that $B_{<\lambda}(\hat a) \subset X_{a'}$.
On the other hand hand, we have $a' \notin
\acl_K(A \cup \{\hat a\})$. By saturation, we find a ball $B \subset K \setminus \acl(A \cup \{\hat a\})$ containing $a'$, and by $0$-h-minimality in the form of Lemma~\ref{lem:type-0-h-min}, all elements of $B$ have the same type over $A\cup \{\hat a\} \cup \RV\eq$ as $a'$. In particular, we have $B_{<\lambda}(\hat a) \subset X_{a''}$ for every $a'' \in B$, which shows that $X$ contains $B \times B_{<\lambda}(\hat a)$.
\end{proof}

\begin{proof}[Proof of Proposition~\ref{prop:dim:basic}, except for (\ref{prop:dim:frontier}) when $\dim X < n$]
We claim that for every definable set $X \subset K^n$, we have $\dim^\acl(X) = \dim(X)$.

To prove this claim, we suppose that $K$ is $\aleph_0$-saturated.
Set $m := \dim^\acl(X)$.
Then there exists a coordinate projection $\ell\colon K^n \to K^m$ such that $\ell(X)$ contains an algebraically independent tuple (over the parameters needed to define $X$). By Lemma~\ref{lem:int}, this
tuple lies in the interior of $\ell(X)$, which shows that $\dim(X) \ge m$. For the other inequality, suppose for contradiction that we have a linear map $\ell\colon K^n \to K^{m'}$ for some $m' > m$ such that $\ell(X)$ has non-empty interior. Then by Lemma~\ref{lem:int}, $\ell(X)$ contains an algebraically independent tuple (over the parameters of $\ell(X)$), so $\dim^\acl(\ell(X)) = m'$, but this contradicts that taking the image under a definable map cannot increase $\dim^\acl$.

This already proves (\ref{dim:2})--(\ref{dim:4}), and we also get
(\ref{dim:1}) for free, using once more that a definable set in $K^d$ has non-empty interior if and only if its dimension is $d$. (The in particular part of (\ref{dim:1}) holds trivially for $\dim^\acl$.)

\medskip

\refstepcounter{dummy}\label{jump:here}
Property (\ref{prop:dim:local}) is proved in \cite{FornHal} in a much more general context; here is a much shorter proof in the present setting:
We may assume that $\acl$ equals $\dcl$ in $\Th(K)$ (using Proposition~\ref{prop:exist:alg:skol}), so that we can apply cell decomposition (Theorem~\ref{thm:cd:alg:skol}) to $X$; we also use Addendum~\ref{add:cd:cont:c} to get continuous centers.

Choose a cell $A_\ell \subset X$ of maximal dimension (of cell-type $(j_i)_{i=1}^n$, with $\sum_i j_i = \dim X =: d$). Then for any $x \in A_\ell$, the local dimension of $X$ at $x$ is $d$. Indeed,
the projection $\pi(A_\ell) \subset K^d$ to the coordinates $\{i \le n \mid j_i = 1\}$ is a cell of cell-type $(1, \dots, 1)$ with continuous center, and hence open. Thus for every sufficiently small ball $B \subset K^n$ around $x$, we have $\pi(X \cap B) = \pi(B)$, witnessing $\dim (X \cap B) \ge d$.

To prove Property (\ref{prop:dim:frontier}) in the case $\dim X = n$, we again
first expand the structure so that $\acl$ equals $\dcl$ and then find a cell decomposition of $\overline X \setminus X$. Since every $n$-dimensional cell has non-empty interior, no such cells can be contained in $\overline X \setminus X$. This implies
$\dim (\overline X\setminus X) < n$.
\end{proof}

\subsection{Jacobian Properties in many variables}
\label{sec:sjp}

There are different ways to generalize the Jacobian Property (Definition~\ref{defn:JP})
to functions $f$ in several variables. The one presented in this subsection (which we now call the Supremum Jacobian Property) has been introduced in \cite{Halup} and is used to obtain t-stratifications. (To be precise, \cite[Definition~2.19]{Halup} is a bit weaker than Definition~\ref{defn:sup-prep} below).

First of all, we need a specific higher-dimensional version of the $\rv$-map.

\begin{notn}[$\ltz$]
Given $\lambda, \mu \in \Gamma_K$, we define $\lambda \ltz \mu$ as $\lambda < \mu  \vee \lambda = \mu = 0$.
\end{notn}

\begin{defn}[Higher-dimensional $\RV$]
For every $n \ge 1$, we define $\RV^{(n)}$ as the quotient $K^n/\mathord{\sim}$, where $x \sim x' \iff |x - x'| \ltz |x|$.
We write $\rv^{(n)}$ for the canonical map $K^n \to\RV^{(n)}$. (For matrices $M \in K^{n \times m}$, we will use the more suggestive notation
$\rv^{(n \times m)}(M)$ instead of $\rv^{(n \cdot m)}(M)$.)
\end{defn}

(Recall that for $x \in K^n$, $|x|$ denotes the maximum norm of $x$.)

\begin{remark}
Note that $\RV^{(1)}$ is just the usual $\RV$.
For $n \ge 2$, $\RV^{(n)}$ is not the same as $\RV^n$, but $\rv^{(n)}$ factors over coordinate-wise $\rv$, so that we have a natural surjection $\RV^n \to \RV^{(n)}$. Moreover, the maximum norm on $K^n$ factors over $\RV^{(n)}$.
\end{remark}

\begin{remark}\label{rem:rvn:GLn}
As explained in \cite[Section~2.2]{Halup}, $\rv^{(n)}$ interacts well with $\GL_n(\cO_K)$; in particular, given $M \in \GL_n(\cO_K)$ and $x \in K^n$, $\rv^{(n)}(Mx)$ is determined by $\rv^{(n \times n)}(M)$ and $\rv^{(n)}(x)$.
\end{remark}

\begin{defn}[Sup-Jac-prop, sup-preparation]\label{defn:sup-prep}
For $X \subset K^n$ open and $f\colon X \to K$, we say that $f$ has the
\emph{Supremum Jacobian Property} (\emph{sup-Jac-prop} for short) on $X$ if
$f$ is $C^1$ on $X$, $\rv^{(n)}(\grad f)$ is constant on $X$, and for every $x_0$ and $x$ in $X$ we have:
\begin{equation}\label{eq:T3/2,mv}
|f(x) - f(x_0) -  ((\grad f)(x_0))\cdot(x - x_0) | \ltz  |\grad f |\cdot |x-x_0|.
\end{equation}
As usual, we consider $(\grad f)(x_0)$ as a matrix with a single row, which we multiply with the column vector $x - x_0$ in the usual way.
We say that a map $\chi\colon X \to \RV^k$ \emph{sup-prepares} $f$ (for some $k \ge 0$) if each $n$-dimensional fiber $F \subset K^n$ of $\chi$ is open and $f$ has the sup-Jac-prop on each such $F$.
\end{defn}

\begin{remark}\label{rem:move:x_0}
One easily checks that the validity of (\ref{eq:T3/2,mv}) does not depend on the precise value of $(\grad f)(x_0)$, but only on $\rv^{(n)}((\grad f)(x_0))$, so it does not play a role whether we evaluate $\grad f$ at $x_0$ or at any other point of $X$.
\end{remark}

\begin{remark}
In the case $n = 1$, (\ref{eq:T3/2,mv}) is equivalent to $\rv(f(x) - f(x_0)) = \rv(f'(x_0))\cdot \rv(x - x_0)$ (which is exactly the main condition of the one-dimensional Jacobian Property; see Definition~\ref{defn:JP}).
For $n\ge 2$ however, (\ref{eq:T3/2,mv}) does not always determine
$\rv(f(x) - f(x_0))$.
Indeed, if e.g.\ $x - x_0$ is orthogonal to $(\grad f)(x_0)$, then (\ref{eq:T3/2,mv}) only imposes an upper bound on $|f(x) - f(x_0)|$.
\end{remark}

\begin{remark}
One cannot expect to be able to sup-prepare definable functions in the stronger sense that $\rv(f(x) - f(x_0))$ is equal to
$\rv(((\grad f)(x_0))\cdot(x - x_0))$ within fibers of $\chi$ (which would correspond to replacing the right hand side of (\ref{eq:T3/2,mv}) by $|((\grad f)(x_0))\cdot(x - x_0)|$). Indeed, consider for example $f(x, y) = y - x^2$, fix any $(x_0, y_0) \in K^2$ and any $\epsilon \in K^\times$, and set $(x,y) := (x_0 + \epsilon, y_0 + 2x_0\epsilon)$. Then
$((\grad f)(x_0,y_0))\cdot((x,y) - (x_0,y_0)) = 0$ but
$f(x,y) - f(x_0,y_0) \ne 0$. For any $\chi$ potentially preparing $f$, we can make such choices such that $(x_0,y_0)$ and $(x, y)$ lie in the same fiber.
\end{remark}

The following lemma states that the sup-Jac-prop is preserved by certain transformations.

\begin{lem}[Preservation of sup-Jac-prop]\label{lem:liptrans}
Let $X \subset K^m$ and $Y \subset K^n$ be open subsets
and let $\alpha\colon X \to Y$ be a $C^1$-map.
Suppose that $\rv^{(n \times m)}(\Jac \alpha)$ is constant on $X$ and that
for every $x_1, x_2 \in X$ we have
\begin{equation}\label{eq:scaling}
|\alpha(x_2) - \alpha(x_1)| = |\Jac \alpha|\cdot |x_2 - x_1|
\end{equation}
and
\begin{equation}\label{eq:alpha}
\rv^{(n)}(\alpha(x_2) - \alpha(x_1)) = \rv^{(n)}((\Jac \alpha)(x_1)\cdot (x_2 - x_1))
\end{equation}
Finally, suppose that $f\colon Y \to K$ is a $C^1$-map such that $f \circ \alpha$ has the sup-Jac-prop (on $X$).
Then $f$ satisfies (\ref{eq:T3/2,mv}) for all $x_0, x \in \alpha(X)$.
\end{lem}

\begin{proof}
Let $x_1, x_2 \in X$ be given and set $y_i := \alpha(x_i)$ and $z_i := f(y_i)$.
In the following, gradients and Jacobians will always be computed at $x_1$ or $y_1$; we
will omit those points from the notation.

What we need to show is:
\begin{equation}\label{eq:lt-want}
|z_2 - z_1 - (\grad f)\cdot (y_2 - y_1)| \ltz  |\grad f| \cdot |y_2 - y_1|.
\end{equation}
By assumption, we have
\begin{equation}\label{eq:lt-have}
|z_2 - z_1 - (\grad (f\circ \alpha))\cdot (x_2 - x_1)| \ltz  |\grad (f\circ \alpha)| \cdot |x_2 - x_1|.
\end{equation}
Applying $\grad (f\circ \alpha) = (\grad f)\cdot (\Jac \alpha)$
to the right hand side of (\ref{eq:lt-have}) gives
\[
|\grad (f\circ \alpha)| \cdot |x_2 - x_1| \le |\grad f| \cdot |\Jac \alpha| \cdot |x_2 - x_1|
\overset{(\ref{eq:scaling})}{=} |\grad f| \cdot |y_2 - y_1|.
\]
On the left hand side of (\ref{eq:lt-have}), we do the following:
\[
(\grad (f\circ \alpha))\cdot (x_2 - x_1) = (\grad f)\cdot (\Jac \alpha)\cdot (x_2 - x_1) \approx
(\grad f)\cdot (y_2 - y_1),
\]
where in the ``$\approx$'', we use (\ref{eq:alpha}) to get an error $e$ with  $e \ltz |\grad f| \cdot |y_1 - y_2|$ (and this is what ``$\approx$'' means here).
Putting things together yields (\ref{eq:lt-want}), as desired.
\end{proof}

The main result of this subsection is that every definable function on $K^n$ can be sup-prepared:

\begin{thm}[Sup-preparation]\label{thm:T3/2.mv}
Suppose that $\Th(K)$ is $1$-h-minimal.
For every $\emptyset$-definable function $f\colon K^n \to K$, there exists a $\emptyset$-definable map $\chi\colon K^n\to \RV^k$ (for some $k\geq 0$) sup-preparing $f$
(in the sense of Definition~\ref{defn:sup-prep}).
\end{thm}

The proof needs some kind of cell decomposition with $1$-Lipschitz centers, as e.g.\ provided by Theorem~\ref{thm:cd:alg:skol}, Addendum~\ref{add:cd:Lip:comp}. To keep this paper more self-contained (since we did not give the proof of Addendum~\ref{add:cd:Lip:comp} in full detail), we will instead prove and use the following weaker version of the addendum; more precisely, this proposition is proved in a joint induction with Theorem~\ref{thm:T3/2.mv}.

\begin{prop}[Twisted boxes with $1$-Lipschitz centers]\label{prop:twibox:1Lip}
Assume $1$-h-minimality and that $\acl$ equals $\dcl$ (in the sense of Assumption~\ref{ass:acl=dcl}).
Then, for every $\emptyset$-definable set $X \subset K^n$, there exists a
$\emptyset$-definable map $\chi\colon X \to \RV^{k'}$ such that each nonempty fiber $F$ of $\chi$ is, up to permutation of coordinates,
an $\RV$-definable twisted box of cell-type $(1, \dots, 1, 0, \dots, 0)$ with $1$-Lipschitz center (i.e., each component $c_i\colon \pi_{<i}(F) \to K$ of the center tuple is $1$-Lipschitz).
\end{prop}

Here, by ``twisted box of cell-type etc.'', we mean: cell of cell-type etc., consisting of a single twisted box; see Definition~\ref{defn:cell}.

A key step in the proof of the proposition consists in swapping
two coordinates, to make the derivative of some center smaller.
This uses the following lemma, which has a similar but simpler proof than Cases 1 and 2 of the proof of Proposition 2.4 of \cite{CCLLip}.
By a ``genuine box'' in $K^n$, we mean a Cartesian product of balls in $K$ (as opposed to the more general notion of ``twisted box'').

\begin{lem}\label{lem:swap:twibox}
Assuming only that $K$ is a valued field in a language containing $\Lval$:
Suppose that \[X = \{(x_1,x_2) \in K^2 \mid \rv(x_1 - c_1) = \xi_1,\ \rv(x_2 - c_2(x_1)) = \xi_2\}\]
is an $A$-definable twisted box, for some set of parameters $A \subset K$ and that $c_2$ has the Jacobian Property. Then either $X$ is a genuine box (i.e., a Cartesian product of two balls), or,  we have $X = \{(x_1,x_2) \in K \times Y  \mid  \rv(x_1 - c_2^{-1}(x_2)) = \xi_3\}$ for some suitable constant $\xi_3 \in \dcl_\RV(A)$ and where $Y$ is the projection of $X$ to the second coordinate.
\end{lem}

\begin{proof}
If $c_2$ is constant there is nothing to do. So, suppose that $c_2$ is not constant.
Then the image $c_2(\pi_{\le 1}(X))$ is an open ball of radius $\rho := |\xi_1|\cdot |c_2'|$ (where $c_2'$ is the derivative of $c_2$). If $\rho \le |\xi_2|$, then whether \(\rv(x_2 - c_2(x_1)) = \xi_2\) holds does not depend on the choice on \(x_1\in \pi_{\le 1}(X)\) and \(X\) is a genuine box.
Otherwise, if $\rho >  |\xi_2|$, then \(c_2^{-1}\) is defined on the whole of \(Y\) and the Jacobian property implies that \(\rv(x_2 -c_2(x_1)) = \xi_2\) if and only if \(\rv(x_1 - c_2^{-1}(x_2)) = -\rv(c_2')^{-1}\cdot\xi_2\).
\end{proof}

\begin{proof}[Proof of Proposition~\ref{prop:twibox:1Lip} in the case $n = 1$]
Prepare $X$ by a finite set $C$ and let $\chi$ be the map given by Lemma~\ref{lem:InextFam}.
\end{proof}

\begin{proof}[Proof of Proposition~\ref{prop:twibox:1Lip},
assuming Proposition~\ref{prop:twibox:1Lip} and Theorem~\ref{thm:T3/2.mv} for $n - 1$]

By an ``$\RV$-partition'' of $X$, we mean a partition into fibers of an $\RV$-definable
map $X \to \RV^k$. Note that if we are already given an $\RV$-partition of $X$, it suffices to prove the proposition for each fiber individually. (Then put everything together using compactness.)

We say that $X$ is a ``thick graph'' (of the map $c_n$) if it is of the form $\{(y, x_n) \in Y \times K\mid \rv(x_n - c_n(y)) = \xi\}$ for some $Y \subset K^{n-1}$, some $c_n \colon Y \to K$, and some $\xi \in \RV$. (Note that we allow $\xi = 0$, which means that $X$ is just the graph of $c_n$.)

Note that it suffices to obtain the claim in the last coordinate, i.e., to $\RV$-partition $X$ into sets that are,
up to permutation of coordinates, thick graphs of $1$-Lipschitz functions $c_n\colon Y \to K$.
After that, the proposition follows by applying induction to $Y$.

Using cell decomposition, we reduce to the case where $X$ is a thick graph of a function $c_n\colon Y \to K$ and $Y$ is a twisted box. In particular, $X$ is a twisted box. We may assume that $Y$ is either open or has empty interior (by treating the interior separately).

Step 1: If $Y$ has empty interior, we apply induction to $Y$ to reduce to the case that $Y$ is a twisted box with $1$-Lipschitz centers and we translate the centers away so that $X$ lives in a subspace of $K^n$ where some of the coordinates are $0$; then apply induction once more to finish. Note that, translating the variables of a $1$-Lipschitz function by $1$-Lipschitz functions yields a  $1$-Lipschitz function and hence the original $X$ does have the required properties.

So suppose from now on that $Y$ is open. By partitioning $Y$, we may assume that $c_n$ is $C^1$ and that $|\partial c_n/\partial x_{i}|$ is constant on $Y$ for each $i$. (Lower-dimensional pieces are treated as in Step 1.)

Step 2: Assume $|\grad c_n| \le 1$.  Using the $n-1$ case of Theorem~\ref{thm:T3/2.mv}, we may assume that $c_n$ has the sup-Jac-prop. This, together with
$|\grad c_n| \le 1$, implies that $c_n$ is $1$-Lipschitz, and hence we are done.

Step 3: So now suppose $|\grad c_n| > 1$. We do an induction on the number of partial derivatives of $c_n$ satisfying $|\partial c_n/\partial x_{i}| > 1$. We suppose without loss that
$|\partial c_n/\partial x_{n-1}| = |\grad c_n|$. Let $Z$ be the projection of $Y$ to the first $n-2$ coordinates.
By further partitioning, we reduce to the case where, for each individual $a \in Z$,
the function $c_n(a, \cdot)$ has the Jacobian Property (using Corollary~\ref{cor:JP} and compactness) and
has an open ball as domain. (Again, lower-dimensional pieces are treated as in Step 1.)

Note that for each $a \in Z$, the fiber $X_{a} \subset K^2$ is a twisted box. By further partitioning $Z$, we may assume that either all of them or none of them are genuine boxes.

Step 3.a: If all fibers are genuine boxes: by induction on $n$, we may assume that the projection $\tilde X$ of $X$ to the coordinates $1, \dots, n-2, n$ is a thick graph of a $1$-Lipschitz function $\tilde c_n\colon Z \to K$.
(This involves permuting the coordinates $1, \dots, n-2, n$.) Then  $X$ is a thick graph of $c_n(z, x_{n-1}) := \tilde c_n(z)$, which is $1$-Lipschitz, so we are done.

Step 3.b: If no fiber is a genuine box, we apply the map $\sigma\colon K^n \to K^n$ swapping the coordinates $n-1$ and $n$. By compactness and Lemma~\ref{lem:swap:twibox},
$\sigma(X)$ is the thick graph of the function $c_{n,\mathrm{new}}$ sending $(x_1, \dots, x_{n-2}, c_n(x_{n-1}))$ to $x_{n-1}$. Since $|\partial c_{n,\mathrm{new}}/\partial x_i|\le |\partial c_{n}/\partial x_i|$ for $i \le n - 2$ and
$|\partial c_{n,\mathrm{new}}/\partial x_n|= 1/ |\partial c_{n}/\partial x_{n-1}| < 1$, $c_{n,\mathrm{new}}$ has fewer partial derivatives bigger than $1$, so we can finish by the induction from Step~3.
\end{proof}

\begin{proof}[Proof of Theorem~\ref{thm:T3/2.mv} in the case $n = 1$]
This follows directly from Corollary~\ref{cor:JP}
and Lemma~\ref{lem:InextFam}: The corollary yields a finite $\emptyset$-definable set $C$ such that (\ref{eq:T3/2,mv}) holds on every ball $1$-next to $C$, and Lemma~\ref{lem:InextFam} then yields a map $\chi\colon K \to \RV^k$ whose $1$-dimensional fibers are exactly those balls.
\end{proof}

\begin{proof}[Proof of Theorem~\ref{thm:T3/2.mv}, assuming Proposition~\ref{prop:twibox:1Lip} for $n$ and Theorem~\ref{thm:T3/2.mv} for $n - 1$]
The proof consists of three parts.

\medskip

Part 1: Some preliminaries:

\medskip

\begin{claim}\label{cl:3/2-acldcl}
It suffices to prove the theorem under the assumption that $\acl$ equals $\dcl$.
\end{claim}

\begin{proof}
Let $\cLas \supset \cL$ be as given by Proposition~\ref{prop:exist:alg:skol}, i.e., $\Th_{\cLas}(K)$ is still $1$-h-minimal, and in $\Th_{\cLas}(K)$, we have $\acl$ equals $\dcl$. Assuming that Theorem~\ref{thm:T3/2.mv} holds for this language, we find an $\cLas$-definable $\chi'\colon K^n\to \RV^{k'}$ sup-preparing $f$.
Since the $\cLas$-structure on $K$ is an $\RV$-expansion of the $\cL$-structure, Lemma~\ref{lem:undo-K-to-RV} provides an $\cL$-definable $\chi\colon K^n\to \RV^k$ such that each fiber $F$ of $\chi$
is contained in a fiber of $\chi'$; in particular, (\ref{eq:T3/2,mv}) holds whenever $(\grad (f|_F))(x_0)$ is defined. It remains to refine $\chi$ in such a way that each of its $n$-dimensional fibers is open.
We do this by splitting each $n$-dimensional fiber $F$ into its interior $\mathring F$ and the remainder. Then indeed $\mathring F$ is open, and
$F \setminus \mathring F$ has dimension less than $n$,
by Proposition~\ref{prop:dim:basic} (\ref{prop:dim:frontier}) applied to
$K^n \setminus F$. (Note that we only use the case of Proposition~\ref{prop:dim:basic} (\ref{prop:dim:frontier}) which we already proved right after the proposition.)
\qedhere(\ref{cl:3/2-acldcl})
\end{proof}

So for the remainder of the proof we assume that $\acl$ equals $\dcl$ (so that we can apply Cell Decomposition and Proposition~\ref{prop:twibox:1Lip}).

Recall that we inductively assume that the theorem holds up to dimension $n-1$. From this, we deduce the following for functions defined on some $n$-dimensional neighborhoods of lower-dimensional subsets of $K^n$.

\begin{claim}\label{cl:on-graphs}
Given any $(n-1)$-dimensional $\emptyset$-definable $Z \subset K^{n}$ and any $\emptyset$-definable $C^1$-function $f$ to $K$ defined on an open neighborhood of $Z$, there exists a $\emptyset$-definable map $\chi\colon Z \to \RV^k$ (for some $k \ge 0$) such that if $F \subset Z$ is an $(n-1)$-dimensional fiber of $\chi$, then
(\ref{eq:T3/2,mv}) holds for every pair $x_0, x \in F$.
\end{claim}

\begin{proof}
We can (and will repeatedly) partition $Z$ into fibers of a $\emptyset$-definable map $Z \to \RV^k$. (If the claim holds for each fiber of such a partition, we then obtain the desired map $Z \to \RV^{k'}$ using compactness.)
By partitioning $Z$ into the twisted boxes of a cell decomposition, we may assume that $\rv^{(n)}(\grad f)$ is constant. By Proposition~\ref{prop:twibox:1Lip}, we may assume that $Z$ is a twisted box of cell-type $(1, \dots, 1, 0)$ with
$1$-Lipschitz center.

Let $\hat Z$ be the projection of $Z$ to the first $n-1$ coordinates, so that $Z$ is the graph of a $1$-Lipschitz function $c \colon \hat Z \to K$. Apply Theorem~\ref{thm:T3/2.mv} to $c$ and to $f \circ \alpha$, where $\alpha\colon \hat Z \to Z, x \mapsto (x, c(x))$, and partition $\hat Z$ and $Z$ accordingly, i.e., so that after the partition, $c$ and $f \circ \alpha$ have the sup-Jac-prop on $\hat Z$. Using that $c$ is $1$-Lipschitz, one obtains that the map $\alpha\colon \hat Z \to Z, x \mapsto (x, c(x))$ satisfies the assumptions of Lemma~\ref{lem:liptrans}. Therefore, the fact that $f \circ \alpha$ satisfies (\ref{eq:T3/2,mv}) on $\hat Z$ implies that $f$ satisfies (\ref{eq:T3/2,mv}) on $Z$, as desired.
\qedhere(\ref{cl:on-graphs})
\end{proof}

Let now a $\emptyset$-definable function $f\colon K^n \to K$ be given (with $n \ge 2$); we need to find a $\emptyset$-definable map $K^n \to \RV^k$ sup-preparing $f$. We more generally allow the domain of $f$ to be any $\emptyset$-definable set $X \subset K^n$.
As in the proof of Claim~\ref{cl:on-graphs}, if we have a $\emptyset$-definable map $\chi\colon X \to \RV^k$, it suffices to sup-prepare the restrictions of $f$ to each fiber of $\chi$. Moreover, fibers of dimension less than $n$ can always be neglected. This argument will be applied repeatedly.

We write elements of $K^n$ as $(x, y)$, with $x\in K^{n-1}$ and $y \in K$.

\medskip

Part 2: Reducing to the case where
$X$ is of the form $\bar X\times B$ for some $\bar X \subset K^{n-1}$ and some ball $B \subset K$, $f$ is $C^1$, $\rv^{(n)}(\grad f)$ is constant, and
\begin{enumerate}
\item[(SJP1)] for each fixed $x \in \bar X$, the function $f(x, \cdot)$ has the sup-Jac-prop.
\end{enumerate}

\medskip

We start by partitioning $K^n$ as follows:

\begin{enumerate}[(a)]
 \item By repeatedly applying the case $n = 1$ (and using compactness), we may assume that $f$ has the sup-Jac-prop fiberwise: for every fixed $x \in K^{n-1}$ and every coordinate permutation $\sigma\colon K^n \to K^n$, the map $y \mapsto f(\sigma(x,y))$ has the sup-Jac-prop.
 \item We moreover assume that $f$ is $C^1$ (using Theorem~\ref{thm:Ck}) and that
 $\rv^{(n)}(\grad f)$ is constant.
\end{enumerate}

By applying Proposition~\ref{prop:twibox:1Lip} and permuting coordinates, we may assume that $X$ is a twisted box with $1$-Lipschitz center:
\begin{equation}\label{eq:F1}
X = \{(x, y) \in K^{n-1} \times K \mid x \in \bar X,\ \rv(y -c(x)) = \rho\}
\end{equation}
for some $\rho \in \RV^\times$, some definable $1$-Lipschitz $c\colon \bar X \to K$ and where $\bar X := \pi_{\le n-1}(X) \subset K^{n-1}$ is the projection to the first $n-1$ coordinates. We may assume that $c$ is $C^1$, and applying the Theorem inductively to $c$ allows us to moreover assume that $c$ has the sup-Jac-prop.

Set $X' := \bar X \times \rv^{-1}(\rho) \subset K^{n-1} \times K$. The bijection $\alpha\colon X' \to X, (x,y) \mapsto (x, y + c(x))$
satisfies the assumptions of Lemma~\ref{lem:liptrans}, so to prove that $f$ has the sup-Jac-prop on some set $F \subset X$, it suffices
to verify that $f' := f \circ \alpha$ has the sup-Jac-prop on $\alpha^{-1}(F)$.
In other words, it remains to sup-prepare $f'\colon X' \to K$.

By (b), $f'$ is $C^1$ and $\rv^{(n)}(\grad f') = \rv^{(n)}((\grad f)(\Jac \alpha))$ is constant (using Remark~\ref{rem:rvn:GLn});
by (a),
$f'(x,\cdot)$ has the sup-Jac-prop for each $x \in \bar X$; thus we are done with Part~2.

\medskip

Part 3: Finishing under the assumptions obtained in Part~2.

\medskip

Recall that $X = \bar X \times B \subset K^{n-1} \times K$.

\begin{enumerate}
 \item[(SJP2)] By induction, we find a map $\chi\colon X \to \RV^k$ such that for each fixed $y \in B$, $g(\cdot, y)$ is sup-prepared by $\chi(\cdot, y)$.
\end{enumerate}
We choose a cell decomposition of $X$ with continuous centers and we refine $\chi$ in such a way that the fibers of $\chi$ are exactly the twisted boxes of the cells. Given such a cell $A_\ell$, let $\bar A_\ell := \pi_{\le n-1}(A_\ell)$ be its projection and let
$c_\ell\colon \bar A_\ell \to K$ be the last component of its center tuple.

Let $Z_\ell$ be the intersection of the graph of $c_\ell$ with $X$.
We apply Claim~\ref{cl:on-graphs} to $Z_\ell$ and $f$, yielding a map $\chi_\ell\colon Z_\ell \to \RV^k$,
we extend $\chi_\ell$ by $0$ to the whole graph of $c_\ell$,
and we replace $\chi$ by the refinement $(x, y) \mapsto (\chi(x, y), (\chi_\ell(x, c_\ell(x)))_\ell)$.
In this way, we achieved the following:
\begin{enumerate}
\item[(SJP3)]
Given any $(x_1, y_1), (x_2, y_2)$ in the same $n$-dimensional fiber of $\chi$, and given any $\ell$, the pair of points $(x_i, c_\ell(x_i))$ ($i =1,2$) satisfies (\ref{eq:T3/2,mv}), provided that both of those points $(x_i, c_\ell(x_i))$ lie in $X$.
\end{enumerate}
Note that since this refinement of $\chi$ only depends on $x$,
each fiber $F$ of $\chi$ is still of the form
\begin{equation}\label{eq:fiber}
F = \{(x, y) \in K^{n-1} \times K \mid x \in \bar F,\rv(y -c_\ell(x)) = \xi\} \subset \bar X \times B,
\end{equation}
for some $\ell$, some $\xi \in \RV$ and where $\bar F = \pi_{\le n-1}(F)$.
Using one last refinement of $\chi$ (also depending only on $x$), we may assume that $\bar F$ is either open or has dimension less than $n-1$,
so that if $F$ is $n$-dimensional, it is open.
To finish the proof of the theorem, we will prove that $f$ has the sup-Jac-prop on each such $n$-dimensional fiber $F$.

We already know that $f$ is $C^1$ on $F$ and that $\rv^{(n)}(\grad f)$ is constant on $F$, so it remains to verify (\ref{eq:T3/2,mv}); thus
let $(x_1,y_1), (x_2, y_2) \in F$ be given.

Recall (Remark~\ref{rem:move:x_0}) that
in (\ref{eq:T3/2,mv}), it does not matter at which point of $F$ we evaluate the gradient $\grad f$. Using this, an easy computation
shows that (\ref{eq:T3/2,mv}) can be verified
in several steps, jumping to certain intermediate points $(x_3,y_3) \in F$ first, namely:
If (\ref{eq:T3/2,mv}) holds for $(x_1,y_1), (x_3, y_3)$ and also for $(x_3,y_3), (x_2, y_2)$, and if moreover
$|(x_1,y_1) - (x_3, y_3)| \le |(x_1,y_1) - (x_2, y_2)|$, then (\ref{eq:T3/2,mv}) follows for $(x_1,y_1), (x_2, y_2)$.
In a similar way, we can also jump through several intermediate points.
Note also that the intermediate points can be arbitrary points of $X$ (and do not need to lie in $F$), since $\rv^{(n)}(\grad f)$ is constant on all of $X$.

We use the notation from (\ref{eq:fiber}) and distinguish three cases:

\medskip

Case 1: $|c_\ell(x_1) - y_1| > |y_2 - y_1|$.
Then we have $(x_1, y_2) \in F$, so we can jump from  $(x_1, y_1)$ to $(x_1, y_2)$ by (SJP1) and from $(x_1, y_2)$ to $(x_2, y_2)$ by (SJP2).

\medskip

Case 2: $|c_\ell(x_2) - y_2| > |y_2 - y_1|$: analogous to Case 1.

\medskip

Case 3: $|c_\ell(x_i) - y_i| \le |y_2 - y_1|$ for $i = 1,2$:
From $y_1, y_2 \in B$, we deduce $c_\ell(x_i) \in B$,
so that $(x_i, c_\ell(x_i))$ lies in $X$. Moreover, we have $|(x_i,y_i) - (x_i, c_\ell(x_i))| \le |(x_1, y_1) - (x_2, y_2)|$.
This means that we can jump from $(x_1,y_1)$ to $(x_1, c_\ell(x_1))$ by (SJP1), then to $(x_2, c_\ell(x_2))$ by (SJP3), and then to $(x_2,y_2)$ by (SJP1) again.
\end{proof}

\subsection{t-stratifications}
\label{sec:t-strat}

In \cite{Halup}, a notion of stratifications in valued fields has been introduced, called ``t-stratifications''. Intuitively, given a definable set $X \subset K^n$, a t-stratification captures, for every ball $B \subset K^n$, the dimension of the space of directions in which $X \cap B$ is ``roughly translation invariant''. This strengthens classical notions of stratifications (like Whitney or Verdier stratifications), which capture rough translation invariance only locally.

The existence proof of t-stratifications given in \cite{Halup} is carried out under some axiomatic assumptions, namely \cite[Hypothesis~2.21]{Halup}. Those assumptions hold in valued fields with analytic structure (in the sense of \cite{CLip}) by
\cite[Proposition~5.12]{Halup} and in power-bounded $T$-convex structures by \cite{Gar.powbd}. We will now show that the assumptions hold in any $1$-h-minimal theory of equi-characteristic $0$, hence implying that t-stratifications exist
in this generality. By the examples of $1$-h-minimal theories given in Section~\ref{sec:examples}, this generalizes
both of the above results.

In this entire section, let $K$ be an equi-characteristic $0$ valued field with $1$-h-minimal theory.

We quickly recall the necessary definitions related to t-stratifications. First, here is the precise notion of ``roughly translation invariant'':

\begin{defn}[Risometries, translatability]\label{defn:trble}
Let $B \subset K^n$ be a ball.
\begin{enumerate}
\item A bijection $f\colon B \to B$ is a \emph{risometry}
 if for every $x_1, x_2 \in B$, we have $\rv^{(n)}(f(x_1) - f(x_2)) = \rv^{(n)}(x_1 - x_2)$.
\item A map $\chi\colon K^n \to Q$ (for some arbitrary set $Q$)
is \emph{$d$-translatable} on $B$, for some $0 \le d \le n$, if there exists a definable (with parameters) risometry $f\colon B \to B$ and a $d$-dimensional sub-vector space $V \subset K^n$ such that for every $x, x' \in B$ satisfying $x - x' \in V$, we have $\chi(f(x)) = \chi(f(x'))$.
\item
A subset $X \subset K^n$ is called $d$-translatable if its characteristic function $1_X\colon K^n \to \{0,1\}$ is $d$-translatable.
\end{enumerate}
\end{defn}

\begin{defn}[t-stratifications]
Let $\chi\colon K^n \to Q$ be a map (for some arbitrary set $Q$) and let $A$ be a set of parameters. An \emph{$A$-definable t-stratification reflecting} $\chi$ is a partition of $K^n$ into $A$-definable sets $S_0, \dots, S_n$ with the following properties:
\begin{enumerate}
 \item $\dim S_d \le d$.
 \item
 Set $\chi'(x) := (\chi(x), d(x)) \in Q \times \{0,\dots, n\}$, where $d(x)$ is defined such that $x \in S_d$ for every $x \in K^n$.
 For each $d \le n$ and each open or closed ball $B \subset S_d \cup \dots \cup S_n$,
this map $\chi'$ is $d$-translatable on $B$.
\end{enumerate}
\end{defn}

\begin{thm}[t-stratifications]\label{thm:t-strats}
Let $K$ be a valued field of equi-characteristic $0$ with $1$-h-minimal theory, and
let $\chi \colon K^n \to Q$ be a $\emptyset$-definable map, where $Q$ is a sort of $\RV\eq$.
Then there exists a $\emptyset$-definable t-stratification $(S_i)_{i \le n}$ reflecting $\chi$.
\end{thm}

\begin{proof}
According to \cite[Theorem~4.12]{Halup}, the existence of t-stratifications follows from \cite[Hypothesis~2.21]{Halup}, so the theorem follows from the following lemma.
\end{proof}

\begin{lem}
Hypothesis~2.21 of \cite{Halup} holds in $1$-h-minimal theories.
\end{lem}

\begin{proof}
The hypothesis consists of the following four parts.
\begin{enumerate}
 \item $\RV$ is stably embedded:

 This is Proposition~\ref{prop:stab}.
\item  Definable maps from $\RV$ to $K$ have finite image:

This is (a special case of) Corollary~\ref{cor:finiterange}.

\item
For every $A \subset K \cup \RV\eq$, every $A$-definable set $X \subset K$ can be $1$-prepared by a finite $A$-definable set $C \subset K$:

This is clear from the definition of $1$-h-minimality and
Lemma~\ref{lem:addconst}.

\item
The theory has the Jacobian Property in the sense of \cite[Theorem~2.19]{Halup}, namely: For every $A \subset K \cup \RV\eq$, every $A$-definable function $f\colon K^n \to K$ can be sup-prepared (in the sense of Definition~\ref{defn:sup-prep}) by an $A$-definable map $\xi\colon K^n \to Q$, where $Q$ is a sort of $\RV\eq$:

Add $A$ as constants to the language. Then (4) is just Theorem~\ref{thm:T3/2.mv}.\qedhere
\end{enumerate}
\end{proof}

Note that the proof we gave here also simplifies the proofs from \cite{Halup} (in the case of fields with analytic structure) and \cite{Gar.powbd} (in the case of $T$-convex structures): In those papers, the proof of (4) was done using a complicated inductive argument using the existence
of t-stratifications in lower dimension. This has been replaced by the more direct proof of our Theorem~\ref{thm:T3/2.mv}.

\medskip

We end this subsection with the promised proof of the missing part of Proposition \ref{prop:dim:basic}, namely that for definable sets $X \subset K^n$, the frontier
 $\overline X\setminus X$ has lower dimension than $X$:

\begin{proof}[Proof of Proposition \ref{prop:dim:basic} (\ref{prop:dim:frontier})]
Choose a t-stratification reflecting
the Cartesian product $\chi(x) := (1_X(x), 1_Y(x))$ of the characteristic functions of $X$ and of the frontier $Y := \bar X \setminus X$, and set $d := \dim Y$.
For dimension reasons, $Y$ contains at least one point $y \in S_d$.
(Note that the definition of t-stratification implies $Y \cap S_j = \emptyset$ for $j > \dim Y$; see \cite[Lemma~3.10]{Halup}.) Assuming $\dim X \le d$, we will show that $y$ cannot be contained in the topological closure of $X$.

Since $S_{\le d-1} := S_0 \cup \dots \cup S_{d-1}$ is closed (by \cite[Lemma~3.17 (a)]{Halup}), there exists a ball $B \subset S_d \cup \dots \cup S_n$ containing $y$.
Let $f\colon B \to B$ be a risometry and $V \subset K^n$ be a vector space witnessing $d$-translatability of $\chi$ on this $B$, as in Definition~\ref{defn:trble}. Since $f$ is a homeomorphism (and $y \notin X$), to obtain $y \notin \bar X$, it suffices to show that $X' := f^{-1}(X \cap B)$ is closed in $B$.
Let $\pi \colon K^n \to K^n/V$ be the canonical map. The definition of $d$-translatability implies that $X' = B \cap \pi^{-1}(\pi(X'))$. Now the assumption $\dim X \le d$ implies $\dim \pi(X') = 0$, so indeed $X'$ is closed in $B$.
\end{proof}

\subsection{Taylor approximation on boxes disjoint from a lower dimensional set}
\label{sec:taylor-box}

We prove a higher-dimensional version of the Taylor approximation Theorem~\ref{thm:high-ord}.
By a \emph{box} in $K^n$, we mean a Cartesian product of balls in $K$.

\begin{thm}[Taylor approximations on boxes]\label{thm:t-high-high}
Given a $\emptyset$-definable function $f\colon K^n\to K$,
there exists a $\emptyset$-definable set $C\subset K^n$ of dimension at most $n-1$ such that for any box $B \subset K^n \setminus C$, $f$ is $(r+1)$-fold differentiable on $B$, for each $i\in\NN^n$ with $|i|= r+1$ one has that  $|f^{(i)}|$ is constant on $B$, and we have
\begin{equation}\label{eq:t-high-high}
|f(x) -  T^{\le r}_{f,x'}(x)| \le \max_{|i| = r+1} |f^{(i)}(x')(x - x')^i|
\end{equation}
for every $x, x' \in B$. (Here, $i$ runs over $n$-tuples and we use the usual multi-index notation.)
\end{thm}

One may search for a variant of (\ref{eq:t-high-high}) that not only holds on boxes disjoint from $C$ but that even holds globally on each open fiber of a map $K^n \to \RV^k$, as in Theorem~\ref{thm:T3/2.mv}; however, we do not yet know how to formulate  such a variant in general.

\begin{proof}[Proof of Theorem~\ref{thm:t-high-high}]
We do an induction over $n$, the case $n = 1$ being Theorem~\ref{thm:high-ord}. Applying the induction hypothesis fiberwise (and using compactness) allows us to find a $C \subset K^n$ ($\emptyset$-definable, of dimension less than $n$) such that for every box $B = \hat B \times B_n \subset K^{n-1} \times K$ disjoint from $C$, for $x_n \in B_n$ and for every $\hat x, \hat x' \in \hat B$, we have:
\begin{equation}\label{eq:hh:1}
 |f(x) -\sum_{|\hat i| \le r} \frac{f^{(\hat i, 0)}(\hat x', x_n)}{\hat i!}(\hat x - \hat x')^{\hat i}| \le \max_{|\hat i| = r+1} |f^{(\hat i, 0)}(\hat x', x_n)(\hat x - \hat x')^{\hat i}|
\end{equation}
For each $\hat i \in \NN^{n-1}$ with $|\hat i| \le r$,  we moreover apply the $n = 1$ case of the theorem to $f^{(\hat i, 0)}(\hat x, \cdot)$ (for each fixed $\hat x$) and $r' := r - |\hat i|$, so that for $\hat x' \in \hat B$ and $x_n, x'_n \in B'$, we have:
\begin{equation}\label{eq:hh:2}
|f^{(\hat i, 0)}(\hat x', x_n) - \sum_{i_n = 0}^{r'} \frac{f^{(\hat i, i_n)}(\hat x', x'_n)}{i_n!}(x_n - x_n')^{i_n}| \le | f^{(\hat i, r' + 1)}(\hat x', x'_n)(x_n - x'_n)^{r' + 1}|.
\end{equation}
Using (\ref{eq:hh:2}) to estimate the $f^{(\hat i, 0)}(\hat x', x_n)$ from the left hand side of
(\ref{eq:hh:1}) yields (\ref{eq:t-high-high}), as desired,
where we also use that we may assume that $|f^{(\hat i, 0)}(\hat x', \cdot)|$ is constant for each $\hat x'$ and each $\hat i$, so that
in the right hand side of (\ref{eq:hh:1}), we can replace
$|f^{(\hat i, 0)}(\hat x', x_n)|$ by $|f^{(\hat i, 0)}(\hat x', x'_n)|$.
\end{proof}

\subsection{Classical cells}\label{sec:cd:classical}

We end this section by recalling the older, more classical notion of cells which has to be used in the absence of the condition that $\acl$ equals $\dcl$ in $\Th(K)$, and by stating the corresponding classical cell decomposition results under the assumption of $1$-h-minimality.

\begin{defn}[Reparameterized cells]\label{defn:reparamcell}
Consider integers $n,k\geq 0$, a $\emptyset$-definable set $X\subset  K^n$, and a $\emptyset$-definable function
$$
\sigma : X\to \RV^k.
$$
Then $(X,\sigma)$ is called a $\emptyset$-definable reparameterized cell (reparameterized by $\sigma$),
if, for each $\xi\in \RV^k$, the set $\sigma^{-1}(\xi)$, when non-empty, is a $\xi$-definable cell with some center tuple $c_\xi$ (see Definition \ref{defn:cell}), such that moreover $c_\xi$ depends definably on $\xi$ and such that the cell-type of $\sigma^{-1}(\xi)$ is independent of $\xi$.  If $(X,\sigma)$ is such a reparameterized cell
then, by a twisted box of $X$ we mean a twisted box of $\sigma^{-1}(\xi)$ for some $\xi$ as in Definition \ref{defn:cell}, and similarly, by the center tuple and the cell-type of $(X,\sigma)$, we mean the definable family of the center tuples of the cells $\sigma^{-1}(\xi)$ (with family parameter $\xi$), and the cell-type of the $\sigma^{-1}(\xi)$, respectively.
\end{defn}

\begin{remark}
In the above definition, one can always modify $\sigma$ in such a way that afterwards,
each $\sigma^{-1}(\xi)$ is either empty or a single twisted box. In a different direction, if the language $\cL$ has an angular component map $\ac$ sending $K$ to its residue field,
then one can take $\sigma$ from Definition \ref{defn:reparamcell} to be residue field valued (instead of $\RV$-valued).
Either of those additional assumptions on $\sigma$ can in particular be imposed on the cells appearing in the following theorem.
\end{remark}

\begin{thm}[Reparameterized cell decomposition]\label{thm:rcd}
Suppose that $\Th(K)$ is $1$-h-minimal.
Consider $n,k\geq 0$ and $\emptyset$-definable sets  $X\subset  K^n$ and $P\subset X\times \RV^k$. We consider $P$ as the function sending $x\in X$ to the fiber $P_x := \{\xi \in \RV^k \mid (x,\xi) \in P\}$.
Then there exists a finite decomposition of $X$ into $\emptyset$-definable reparameterized cells $(A_i,\sigma_i)$ such that moreover $P$ (considered as a function) is constant on each twisted box of each $A_i$.

The other addenda of Theorem \ref{thm:cd:alg:skol} can be adapted in a similar way.
In particular, for the analogue of Addendum \ref{add:cd:alg:range} of \ref{thm:cd:alg:skol}, given finitely many $\emptyset$-definable  $f_j:X\to K$ and assuming $n=1$, we can moreover assume that there are $\emptyset$-definable reparameterized cells $(B_{ij},\tau_ {ij})$  such that $f_j(A_i)=B_{ij}$ and such that any twisted box of $A_i$ is mapped by $f_j$ onto a twisted box of $B_{ij}$.

For the analogue of Addendum \ref{add:cd:Lip:comp} of \ref{thm:cd:alg:skol}, up to allowing a well-chosen coordinate permutation for each for each $i$, we can moreover for each $\xi$ obtain that the center of $\sigma_i^{-1}(\xi)$ is $1$-Lipschitz.
\end{thm}
\begin{proof}
The proof is analogous to the proof of Theorem \ref{thm:cd:alg:skol} and its addenda.
\end{proof}
Note that the above version of Addendum \ref{add:cd:Lip:comp} is weaker than the original Addendum \ref{add:cd:Lip:comp} of Theorem \ref{thm:cd:alg:skol}, since instead of obtaining finitely many $1$-Lipschitz centers,
we now only obtain that for each $\xi$ separately, the center of
$\sigma_i^{-1}(\xi)$ is $1$-Lipschitz;
this corresponds to an infinite partition.

A similar phenomenon also arises with Theorem \ref{thm:cd:alg:piece:Lipschitz} in the absence of the condition that $\acl=\dcl$, as in Theorem ~2.1.7 of \cite{CCL-PW}: $f$ being locally $1$-Lipschitz implies globally $1$-Lipschitz only after some infinite partition of the domain.

\subsection{Motivic integration}

In this subsection, we show that the Cluckers--Loeser style of motivic integration from \cite{CLoes} can be generalized to the $1$-h-minimal setting under some natural extra assumptions about the induced structure on $\VG$ and $\RF$. It is beyond the scope of this paper to recall what motivic integration is. Instead, we rely on \cite{CLbounded}, which generalizes an essential part of \cite{CLoes} to more characteristics. In \cite{CLbounded}, some axiomatic assumptions are given under which motivic integration works; below, we just verify those assumptions. More precisely, we check the conditions of  \cite[Definition~3.9]{CLbounded} (in the variant from \Cref{rem:T-fields} below) for the theory $\Th_{\cL}(K)$ where $K$ is a Henselian, discretely valued field of equi-characteristic zero with an $\cL$-structure which is $1$-h-minimal and satisfies some natural conditions;
motivic integration is developed uniformly for such $K$ in \cite{CLbounded} (as well as in mixed characteristic assuming bounds on the ramification). 
Note that the framework with motivic additive characters and motivic Fourier transforms from \cite{CLexp} can be generalized similarly to the $1$-h-minimal setting under the same extra assumptions; as this would take more time to explain, we leave such adaptation to some future work.

\begin{remark}\label{rem:T-fields}
We indicate  a small adaptation of the framework of \cite{CLbounded} in the equi-characteristic zero case. If one only considers equi-characteristic zero Henselian valued fields, then the language $\cL_{\rm high}$ used in \cite{CLbounded} can be replaced by the Denef-Pas language (i.e., one omits the higher order angular component maps $\ac_n$ for $n>1$ from $\cL_{\rm high}$, one omits the symbol $\pi$ for the uniformizer, and, one keeps only the sorts $\VF$, $\RF$, and $\VG$). One can adapt all definitions and results from \cite{CLbounded} correspondingly, as long as one only considers equi-characteristic zero Henselian valued fields.
Let us call this variant (with the Denef-Pas language) the non-high variant of \cite{CLbounded}.
\end{remark}

\begin{thm}\label{thm:mot.int}
Let $K$ be an equi-characteristic zero, Henselian, discretely valued field, equipped with $\cL$-structure for some language $\cL \supset \Lval$.
Suppose that $\Th(K)$ is $1$-h-minimal and that we additionally have the following:
\begin{enumerate}
 \item The language $\cL$ contains an angular component map $\overline{\mathrm{ac}}\colon K \to k := \cO_K/\cM_K$.
 \item The induced structure on $\Gamma_K^\times$ is that of the pure ordered abelian group, i.e., every definable subset of $(\Gamma_K^\times)^\ell$ (for every $\ell$) is definable in the language $\{+,<\}$ on $\Gamma_K^\times$.
 \item The value group $\Gamma_K^\times$ and the residue field $k$ are orthogonal, i.e., every definable subset of $k^m \times (\Gamma_K^\times)^\ell$ (for every $m$ and $\ell$) is a finite union of sets of the form $Y \times Z$, with $Y \subset k^m$ and $Z \subset (\Gamma_K^\times)^\ell$.
\end{enumerate}
Then the $\cL$-theory $\cT := \Th(K)$  satisfies the assumptions from the non-high variant of \cite[Definition~3.9]{CLbounded} (see \Cref{rem:T-fields}).
\end{thm}

\begin{remark}
In Items (1) and (3) of \Cref{thm:mot.int}, the residue field $k$ is naturally considered as a subset of $\RV_K$, consisting of the elements of norm $1$ inside $\RV_K$, and the zero element.
\end{remark}

\begin{proof}[Proof of \Cref{thm:mot.int}]
According to \cite[Definition~3.9]{CLbounded} and with the terminology of \cite[Section 3.2]{CLbounded}, we have to verify that $\cT$ is split, finitely $b$-minimal, has the Jacobian property, and has at least one model which is a $(0,0,0)$-field, all in their non-high variants from \Cref{rem:T-fields} above.

By Proposition~\ref{prop:b-min}, $\cT$ is $b$-minimal, and the additional condition ``finitely'' follows from Lemma~\ref{lem:fin-inf}. By Corollary~\ref{cor:JP} (with $\lambda = 1$), $\cT$ has the Jacobian property, and the given field $K$ is model which is a $(0,0,0)$-field.

That $\cT$ is split (\cite[Definition~3.6]{CLbounded}) means that
the induced structure on $\Gamma_K^\times$ is the pure ordered abelian group (which is one of our assumptions) and that $k$ and $\Gamma_K^\times$ are ``orthogonal with control of parameters'' in the following sense:
For every set $A \subset K \cup \RV\eq$ of parameters, every $A$-definable subset of $k^m \times (\Gamma_K^\times)^\ell$ is a finite union of sets of the form $Y \times Z$, where $Y \subset k^m$ and $Z \subset (\Gamma_K^\times)^\ell$ are $A$-definable.

It turns out that using the order on $\Gamma_K^\times$, this strong form of orthogonality follows from the normal one, namely as follows:
Let an $A$-definable set $X\subset k^m \times (\Gamma_K^\times)^\ell$ be given, and let $(\Gamma_K^\times)^\ell = \bigcup_i Z_i$ be the partition defined by the condition that $z, z' \in (\Gamma_K^\times)^r$ lie in the same set $Z_i$ if and only if the fibers $X_z$ and $X_{z'}$ (in $k^m$) are equal. By
(normal) orthogonality of $k$ and $\Gamma_K^\times$, this partition is finite. Using the order on $\Gamma_K^\times$, we can definably distinguish between the sets $Z_i$, so each $Z_i$ is $A$-definable. Denote by $Y_i \subset k^m$ the fiber corresponding to $Z_i$ (so that $X_z = Y_i$ for every $z  \in Z_i$). Those are also $A$-definable, and as desired, $X$ is the union of the sets $Y_i \times Z_i$.
\end{proof}

\section{Examples of Hensel minimal structures}\label{sec:examples}

In this section, we provide many examples of Hensel minimal valued fields of equi-characteristic $0$, with various kinds of languages. Note that the structures from \cite{Rid,Scanlon} for which a proper subfield of $K$ is definable (e.g.~as the fixed field of an automorphism) are not covered by our framework.  In some cases, we prove $\omega$-h-minimality, in others, we only obtain $1$-h-minimality (namely for power-bounded $\cT$-convex structures). We also provide examples of valued fields of mixed characteristic.

So far, we introduced Hensel minimality in equi-characteristic $0$ only. A mixed characteristic analogue is developed in the sequel \cite{CHRV}.
For the sake of the examples below, we quickly give a definition of mixed characteristic Hensel minimality already here.

\subsection{Mixed characteristic}\label{ssec:mixed}
As usual, we fix a theory $\cT$ of valued fields in a language $\cL$ containing the language $\Lval$ of valued fields. We do require the characteristic of models to be $0$, but we allow the models to have arbitrary residue field characteristic.

\begin{notn}[Equi-characteristic $0$ coarsening]
Given a model $K \models \cT$,
we write $\cO_{K,\eqc}$ for the smallest subring of $K$ containing $\cO_K$ and $\QQ$ and we let $|\cdot|_{\eqc}\colon K \to \Gamma_{K,\eqc}$ be the corresponding valuation. (Thus, $|\cdot|_{\eqc}$ is the finest coarsening of $|\cdot|$ which has equi-characteristic $0$; note that $|\cdot|_{\eqc}$ can be a trivial valuation on $K$.) If $|\cdot|_{\eqc}$ is a nontrivial valuation (i.e., $\cO_{K,\eqc} \ne K$), then we
also use the following notation: $\rv_{\eqc}\colon K \to \RV_\eqc$ is the leading term map with respect to $|\cdot|_{\eqc}$;
given $\lambda \in \Gamma_{K,\eqc}$, $\rv_{\lambda}\colon K \to \RV_\lambda$ is the leading term map with respect to $\lambda$;
and $\cL_{\eqc}$ is the extension of $\cL$ by a predicate symbol for $\cO_{K,\eqc}$.
\end{notn}

\begin{defn}[Hensel minimality in mixed characteristic]\label{defn:mixed}
Let $\cT$ be a theory
of valued fields of characteristic $0$ (and arbitrary residue field characteristic) in a language $\cL$ containing the language $\Lval$ of valued fields. Let $\ell\geq 0$ be either an integer or $\omega$. We say that $\cT$ is \emph{$\ell$-\heqc-minimal}
if for every model $K \models \cT$ the following holds:
If the valuation $|\cdot|_{\eqc}$ on $K$
is non-trivial, then
the $\cL_{\eqc}$-theory of $K$, when considered
as a valued field with the valuation $|\cdot|_{\eqc}$, is $\ell$-h-minimal in the sense of Definition~\ref{defn:hmin}.
\end{defn}

\begin{remark}\label{rem:1heqc}
For theories of valued fields of equi-characteristic $0$, $\ell$-\heqc-minimality is trivially equivalent to $\ell$-h-minimality, for each $\ell$ (since $\cO_{K,\eqc} = \cO_K$).
\end{remark}

A natural adaptation of the definition of h-minimality to mixed characteristic would consist in (suitably) replacing $\lambda$-preparation by $|m|\lambda$-preparation, for some nonzero integer $m$. In \cite{CHRV}, we show in the case of $1$-h-minimality that this adaptation is equivalent to
the above notion of $1$-\heqc-minimality. Moreover, we prove that many preparation results from this paper can be adapted to mixed characteristic in a similar manner, replacing $\lambda$-preparation by $|m|\lambda$-preparation.

\subsection{Valued fields with or without analytic structure}
\label{sec:analyt}

In this subsection, we prove $\omega$-h-minimality of arbitrary Henselian valued fields $K$  of equi-characteristic $0$ with analytic structure in the sense of \cite{CLip}. If $K$ has mixed characteristic and analytic structure, we will show that it is $\omega$-\heqc-minimal, which means that its equi-characteristic $0$ coarsenings are $\omega$-h-minimal.
As so often, obtaining results in the positive equi-characteristic case seems completely out of reach at present.

The pure valued field language is a special case of an analytic structure on $K$. Nevertheless, we will treat this case separately in this section, avoiding the machinery of analytic structures and instead building only on classical quantifier elimination results. Note that before, the only known proofs of the Jacobian Property (Corollary~\ref{cor:JP}) either went via analytic structures (as in \cite{CLoes, CLip}) or were restricted to algebraically closed valued fields (as in \cite{HK}).

The proofs of $\omega$-h-minimality (i) in the pure field language and (ii) in fields with analytic structure are very similar, the main differences being that, in case (ii) we require additional input from \cite{CLip}. We therefore formulate both proofs simultaneously, tagging differences with (i) and (ii).

We fix the following language and structure for the remainder of this subsection.
Note that in the mixed characteristic case, the Definition~\ref{defn:mixed} of  $\omega$-\heqc-minimality requires us to work with two different valuations simultaneously (though in Case~(i), just understanding the coarser valuation is enough for the proof of $\omega$-h-minimality).
\begin{enumerate}[(i)]
 \item Let $\cL := \Lval \cup \{0,1,-\}  = \{+,-,0,1,\cdot,\cO_K\}$ be the pure valued field language together with $0,1,-$ with their natural meaning, and let $K$ be a Henselian valued field of equi-characteristic $0$, considered as an $\cL$-structure.
\item
 Let $\cA=(A_{m,n})_{m,n}$ be a separated Weierstrass system in the sense of \cite{CLip}. Let $K$ be a characteristic zero valued field (possibly of positive residue field characteristic) with a separated analytic $\cA$-structure, as in \cite[Definitions 4.1.5 and 4.1.6]{CLip}. Note that by \cite[Proposition 4.5.10 (i)]{CLip}, any such $K$ is Henselian.

 We denote the valuation ring of $K$ by $\cO_{K,\fine}$, and we fix an equi-characteristic $0$ coarsening
 $\cO_K \supset \cO_{K,\fine}$. (If $\cO_{K,\fine}$ itself is of equi-characteristic $0$, one can as well
  choose $\cO_K = \cO_{K,\fine}$.) We write $\cM_{K,\fine}$, $\cM_{K}$ for the corresponding maximal ideals and $|\cdot|_{\fine}\colon K \to \Gamma_{K,\fine}$, $|\cdot|\colon K \to \Gamma_K$
 for the corresponding valuations.

We let, still in this case (ii), $\cL$ be the extension of the language from (i) by a
function symbol for field division (extended by zero on zero), by a predicate symbol for $\cO_{K,\fine}$ and
by one function symbol for each $f$ in $A_{m,n}$, interpreted as a function $\cO_{K,\fine}^m\times\cM_{K,\fine}^n\to K$ via the analytic $\cA$-structure on $K$ (and extended by $0$ outside of its domain).
\end{enumerate}

\begin{thm}[Fields with analytic structure]\label{thm:hmin analytic}
Let $(K, |\cdot|)$ be an equi-characteristic $0$ valued field in a language $\cL$ as above in (i) or (ii). (Recall that in particular, in Case~(ii) there is a valuation ring $\cO_{K,\fine}$ of $K$ which may differ from $\cO_K$.)
Then the $\cL$-theory of $K$ with the valuation $|\cdot|$ is $\omega$-h-minimal.
\end{thm}

In one word, the idea of the proof of Theorem \ref{thm:hmin analytic} goes via quantifier elimination of valued field quantifiers, in a language to which
the sort $\RV_\lambda$ and the map $\rv_\lambda$ has been added, for some given $\lambda \le 1$ in $\Gamma_K^\times$. The conditions of $\omega$-h-minimality are then easily checked. Similar quantifier elimination results have been proved but not yet with such an $\RV_\lambda$, so we take care to give sufficient details.

The first technical ingredient, inspired by \cite{vdDHM}, is the following:

\begin{lem}\label{lem:lin an}
Let \(F\leq K\) be a subfield which is moreover an $\cL$-substructure of $K$, let $\lambda \le 1$ be an element of $\Gamma_K^\times$ and let \(\tau(x)\) be an \(\cL(F)\)-term, with \(x\) a single variable. Then there exists a finite set $C \in K$ consisting only of algebraic elements over $F$ such that
\(\rv_\lambda(\tau(x))\) only depends on $(\rv_{\lambda}(x-c))_{c \in C}$.

Moreover, if we restrict (in our desired property for $C$) to the \(x\) that are solutions of a given degree \(d\) polynomial equation $P \in F[X]$, (i.e., if
we want to obtain the implication $(\rv_{\lambda}(x-c))_{c \in C} = (\rv_{\lambda}(x'-c))_{c \in C} \Rightarrow \rv_\lambda(\tau(x)) = \rv_\lambda(\tau(x'))$ only under the assumption $P(x) = P(x') = 0$), then the elements of \(C\) can all be assumed to have degree strictly less than \(d\) over \(F\).
\end{lem}

\begin{proof}[Proof in Case (i)]
Note that our $\cL(F)$-term $\tau$ is simply a polynomial in $F[x]$, so the first part of the lemma is immediate from \cite[Proposition 3.6]{Flen}, and its proof, namely using for $C$ the set of all roots of all derivatives of $\tau$ (including the roots of $\tau$ itself). Note that the proof of \cite[Proposition 3.6]{Flen} yields that the ``Swiss cheeses'' $U_i$ appearing in the statement of the proposition are
$1$-prepared by $C$.

For the second part, choose a polynomial $Q$ of degree less than $\deg (P)$ such that the polynomial $\tau$ is congruent to $Q$ modulo $P$. Then for every $x \in K$ which is a zero of $P$, we have $\tau(x) = Q(x)$, so we may as well replace $\tau$ by $Q$. Then choosing $C$ as in the proof of the first part yields the claim.
\end{proof}

\begin{proof}[Proof in Case (ii)]
We apply \cite{CLip} to $K$ with the valuation $|\cdot|_{\fine}$.
The main idea is to use \cite[Theorem 5.5.3]{CLip} to reduce to Case (i).
As a preparation, note that it is sufficient to obtain the conclusion of the lemma for $x \in \cO_K$. Indeed,
the $x$ outside of $\cO_K$ can be treated by applying the lemma separately to the $\cL(F)$-term $\tau(1/x)$.

By \cite[Theorem 5.5.3 and Remark 5.5.4]{CLip}, there is a cover of the valuation ring $\cO_{K_{\mathrm{alg}}}$ of the algebraic closure of $K$ by finitely many $F$-annuli $\cU_i$ (cf.\ \cite[Definition~5.1.1]{CLip}) and there is a finite \(F\)-definable set \(S \subset K\), such that, on each \(\cU_i \setminus S\), we have \[\tau(x) = G_i(x)/H_i(x)\cdot E_i(x),\] where \(G_i, H_i\in F[x]\) are polynomials and where \(E_i\) is an \(\cL(F)\)-term which is a strong unit on $\cU_i$ (cf.\ \cite[Definition~5.1.4]{CLip}).
Indeed, that these data are defined over $F$ follows from \cite[Remark 5.5.4]{CLip} with $K'=F$.

Let \(P_S \in F[x]\) be the polynomial whose set of roots is \(S\) and let \(P_j \in F[x]\) be the collection of polynomials appearing in the definition of the $F$-annuli \(\cU_i\). (In particular, each $\cU_i$ is
defined by a boolean combination of inequalities between the valuations of the \(P_j\).) Since $|\cdot|_{\fine}$ factors over $\rv_\lambda$, whether an element $x \in Z$ lies in $\cU_i$ is determined by $(\rv_\lambda(P_{j}(x)))_j$.
Moreover, by \cite[Lemma 6.3.12]{CLip} and \cite[Remark A.1.12]{CLips}, $\rv_\lambda(E_i)$ only depends on $(\rv_\lambda(P_{j}(x)))_j$.
Thus to prove the lemma for $\tau$, it suffices to prove it for the polynomials $G_i$, $H_i$, $P_j$ and \(P_S\) (considered as functions on $K$). But this has already been done in Case (i).
\end{proof}

The second ingredient is a quantifier elimination result. We
fix $\lambda \le 1$ in $\Gamma_K^\times$ and consider the following expansion of $K$ in a language $\Lqe \supset \cL$:
We add $\RV_\lambda$ as a new sort, together with the map $\rv_\lambda$ and with the $\emptyset$-induced structure on $\RV_\lambda$,
 i.e., one predicate for each $\emptyset$-definable subset of $\RV_\lambda^n$, for every $n$.

\begin{prop}\label{prop:QE analytic}
The $\Lqe$-theory of $K$ eliminates field quantifiers.
\end{prop}

The particularity of Proposition \ref{prop:QE analytic} is not only that it has $\rv_\lambda$ and $\RV_\lambda$ on top of the analytic structure, but also that there are two valuation rings at play, namely $\cO_{K,\fine}$ (which may be of mixed characteristic) and $\cO_K$ (which is of equi-characteristic zero).

\begin{proof}[Proof of Proposition \ref{prop:QE analytic}]
First note that in the case $\lambda = 1$, this result is known:
in Case (i), this is \cite[Proposition~4.3]{Flen} (or also \cite[Theorem~B]{Basarab})
and in Case (ii), it is \cite[Theorem 3.10]{Rid}. This already implies a partial result for arbitrary $\lambda \le 1$, namely:
\begin{condition}\label{cond:field-vars}
Every $\Lqe$-formula $\phi(x)$ having only $K$-variables is equivalent to a field quantifier free formula.
\end{condition}
 Indeed, $\phi(x)$ is equivalent to an $\Lqe[1](\lambda)$-formula $\psi(x)$ without $K$-quantifiers, and each $\RV$-quantifier of $\psi(x)$ can easily be replaced by some $\RV_\lambda$-quantifiers.

To prove the general case, we need the following variants of results from \cite{Flen}.

Fix a non-zero polynomial $P \in K[x]$ and an element $a_0 \in K$.
Let $c_i \in K$ be the coefficients of $P$ when developed around $a_0$, i.e., $P(x) = \sum_i c_i (x-a_0)^i$.

Given $b \in K$, we say, cf. \cite[Definition 3.1]{Flen}, that $P$ has a \(\lambda\)-collision at $b$ around \(a_0\) if \(|P(b)| < \lambda \max_i |c_i(b-a_0)^i|\). Note that whether \(P\) has a \(\lambda\)-collision at \(b\) around \(a_0\) only depends on $\rv_\lambda(c_i)$ (for all $i$) and $\rv_\lambda(b - a_0)$. (This will be useful later.)

\begin{claim}\label{cl:col root}
Suppose that the above $P \in K[x]$ has a \(\lambda\)-collision at \(b\) around \(a_0\). Then there exists an integer $n\geq 0$ with $n < \deg P$ and an element $b'\in K$ ``close to $b$'' such that \(P^{(n)}(b') = 0\). Here,  ``close to $b$'' means $\rv_\lambda(b'-a_0) = \rv_\lambda(b-a_0)$ if $n = 0$ and
$\rv(b'-a_0) = \rv(b-a_0)$ if $n \ge 1$, and $P^{(n)}$ stands for the $n$-th derivative, with $P^{(0)}=P$.
\end{claim}

\begin{proof}
Without loss, we may assume that \(a_0 = 0\), \(b=1\) and \(\max_i |c_i| = 1\). Let \(Q = \res(P)\), the reduction of $P$ modulo $\cM_{K}$. We have \(Q(1) = 0\). Let \(n\) be such that $1$ is a root of multiplicity $1$ of \(Q^{(n)}\). Then Hensel's Lemma yields a root \(b'\) of \(P^{(n)}\) such that \(|b'-1| \leq |P^{(n)}(1)|\). Since $|P^{(n)}(1)| \le 1$, this implies $\rv(b') = \rv(b)$. In the case $n = 0$, that $1$ is not a root of $Q'$ implies $|P'(1)| = 1$.
Together with $|P(1)| \le \lambda$, Hensel's Lemma implies $|b' - b| < \lambda$ and hence $\rv_\lambda(b') = \rv_\lambda(b)$.
\qedhere(\ref{cl:col root})
\end{proof}

\begin{claim}\label{cl:exists root}
Suppose that the above $P \in K[x]$ has no common zero with any of its proper derivatives and fix $\xi\in\RV^\times_\lambda$. The following are equivalent:
\begin{enumerate}
\item There exists a root $b \in K$ of \(P\) with \(\rv_\lambda(b-a_0) = \xi\);
\item there exists \(b\in K\) such that\\
(2a) \(\rv_\lambda(b-a_0) = \xi\), \(P\) has a \(\lambda\)-collision at \(b\) around \(a_0\), and\\
(2b) for every root $a \in K$ of every proper derivative of $P$,
\(P\) has a \(\lambda\)-collision at \(b\) around $a$.
\end{enumerate}
\end{claim}

\begin{proof}
If \(b\) is a root of \(P\), then \(P\) has a \(\lambda\)-collision at \(b\) around any \(a\neq b\). So (1) implies (2). Let us now assume that we have \(b\) such that (2) holds. If (1) does not hold, then by Claim~\ref{cl:col root} (and (2a)) there exists a root \(a\) of some proper derivative of $P$ with \(\rv(a-a_0) = \rv(b-a_0)\). Pick the closest such \(a\) to \(b\). By Claim~\ref{cl:col root}, around \(a\) this time (and using (2b)), there exists a root \(c\) of some \(P^{(m)}\) with \(\rv(c-a) = \rv(b-a)\); in particular, $|b-c| < |b-a|$, a contradiction to our choice of $a$.
\qedhere(\ref{cl:exists root})
\end{proof}

We now come back to the actual proof of field quantifier elimination. We abbreviate ``field quantifier free'' by ``fqf''.
It suffices to prove the following: Suppose
that $K' \equiv K$ is $|K|^+$-saturated, that $A \subset K \cup \RV_{K,\lambda}$ and $A' \subset K' \cup \RV_{K',\lambda}$ are substructures and that
$\alpha\colon A \to A'$ is an fqf-elementary bijection, i.e., that it preserves the validity of fqf formulas.
Then for any $a \in K$, there exists an $a'$ such that $\alpha$ extends to an fqf-elementary map sending $a$ to $a'$.

For $\alpha$ to be fqf-elementary, it suffices that it is an isomorphism of substructures and that $\alpha|_{A \cap \RV_{K,\lambda}}$ is fqf-elementary. Indeed, suppose that $\phi$ is an
fqf $\Lqe(A)$-sentence that holds in $K$. Then without loss,
$\phi = \psi((\rv_\lambda(\tau_i))_i)$ for some $\Lqe(A\cap K)$-terms $\tau_i$ and for $\psi(y)$ an fqf $\Lqe(A \cap \RV_{K,\lambda})$-formula.
Let $\xi_i \in A \cap \RV_{K,\lambda}$ be the interpretation of $\rv_\lambda(\tau_i)$. Then $\phi$ follows from $\bigwedge_i \tau_i = \xi_i$ (which is quantifier free) and
$\psi((\xi_i)_i)$ (which is an fqf $\Lqe(A \cap \RV_{K,\lambda})$-sentence).

We may assume $\RV_{K,\lambda} \subset A$, since $\alpha_{A \cap \RV_{K,\lambda}}$ extends to an fqf-elementary map on $\RV_{K,\lambda}$ and the union of this extension with the original $\alpha$ is an isomorphism of substructures. In particular, when further extending $\alpha$, we now only need to make sure that it remains an isomorphism of substructures.
In terms of formulas, this means (by a usual compactness argument) that given a quantifier free $\Lqe(A)$-formula $\phi(x)$ with $x$ a valued field variable, we need to check that $K \models \exists x\,\phi(x)$ implies $K' \models \exists x\,\phi^\alpha(x)$ (where ``$\phi^\alpha$'' is the $\Lqe(A')$-formula obtained from $\phi$ by applying $\alpha$ to the parameters from $A$).

We may assume that $F := K \cap A$ is a subfield. In Case (ii), this is automatic, since $\Lqe$ contains field division. In Case (i), the ring homomorphism $\alpha_{K \cap A}$ uniquely extends to the fraction field $F$ of $K \cap A$, and extending $\alpha$ in this way yields an isomorphism of substructures, since
for $\frac{b}{b'} \in F$ ($b, b' \in K \cap A$), we have
$\rv_\lambda(\frac{b}{b'}) = \frac{\rv_\lambda(b)}{\rv_\lambda(b')}$.

From now on, we identify $A$ with its image $\alpha(A)$.
Fix $a \in K$ and fix a quantifier free $\Lqe(A)$-formula $\phi(x)$ such that $K \models \phi(a)$ holds.
We need to show that $K' \models \exists x\colon \phi(x)$. To do so, we will successively reduce to simpler formulas, until we can get rid of the $K$-quantifier $\exists x$.

Let $P \in F[x]$ be the minimal polynomial of $a$ over $F$ and set $d := \deg P$. If \(a\) is transcendental over $F$, we set \(P := 0\) and \(d := \infty\).
By induction on $d$, we may assume:
\begin{condition}\label{cond:F}
$F$ contains all roots $b$ in \(K\) of polynomials over $F$ of degree strictly less than \(d\).
\end{condition}

As before, we can assume that the above formula $\phi$ is of the form $\phi(x) = \psi((\rv_\lambda(\tau_i(x)))_i)$ for some $\cL(F)$-terms $\tau_i$.
By Lemma~\ref{lem:lin an} (and (\ref{cond:F}))
$\phi(x)$ is equivalent, in the structure $K$, to
a formula of the form \[\phi'(x) = \psi'((\rv_\lambda(x - c_j)_j) \wedge P(x) = 0\] for some $c_j \in F$, where
$\psi'(y) = \psi((\eta_i(y))_i)$ for suitable $\Lqe(F)$-definable functions $\eta_i$.
We claim that this equivalence also holds in $K'$, so that we can without loss replace $\phi$ by $\phi'$.

To prove the claim, note that the
equivalence $\phi \leftrightarrow \phi'$ follows from an $\Lqe(F)$-sentence $\chi$, namely
$\chi = \bigwedge_i \forall x \in K\colon \tau_i(x) = \eta_i((\rv_\lambda(x - c_j))_j)$.
Since $\chi$ only uses $K$-parameters, we already know (by (\ref{cond:field-vars})) that it is equivalent to a fqf $\Lqe(F)$-sentence $\chi'$ (modulo only the $\Lqe$-theory of $K$, i.e., without using the specific embedding of $F$ into $K$). Thus the truth of $\chi'$ is preserved by $\alpha$, so that we obtain the desired equivalence in $K'$.\footnote{In Case (i), we do not need to invoke another field quantifier elimination result. As in Lemma~\ref{lem:lin an}, we assume that each \(\tau_i\) has degree smaller than \(P\) and we choose, as $c_j$, all the roots of derivatives of all \(\tau_i\), including \(\tau_i\) itself. Then, by \cite[Proposition 3.6]{Flen}, for every \(x\) and $i$, there exists a $j$ such that if we write \(\tau_i(x) = \sum_k a_{k} (x-c_{j})^k\), the sum
\(\sum_k \rv_\lambda(a_{k})\rv_\lambda(x-c_{j})^k\) is well-defined and hence equal to \(\rv_\lambda(\tau_i(x))\). Using that the well-definedness of the sum is an fqf condition, we can define $\eta_i \colon (\rv_\lambda(x-c_{j}))_j \mapsto \rv_\lambda(\tau_i(x))$ without field quantifiers, so that the equality $\eta_i(\rv_\lambda(x-c_{j}))_j = \rv_\lambda(\tau_i(x))$ is preserved by $f$.}

Next, note that we can get rid of all the $c_j$ appearing in $\phi'$ except for the one closest to $a$, i.e., denoting that closest $c_j$ by $c$, we can replace $\phi'$ by
\[\phi''(x) = \psi''(\rv_\lambda(x - c)) \wedge P(x) = 0.\] Indeed, one can easily choose $\psi''$ in such a way that $K \models \phi''(a)$ and that
$\Th_{\Lqe}(K)$ implies $\phi'' \to \phi'$.

Now $\exists x \colon \phi''(x)$ is equivalent to
\(\exists\xi \in \RV_\lambda\colon \psi''(\xi)\wedge(\exists x\colon  \rv_\lambda(x-c) = \xi \wedge P(x) = 0)\),
and it remains to get rid of the $\exists x$ in that formula.
If $P$ is the zero polynomial, then the $\exists x$ part is trivially true and we are done.
Otherwise, by Claim \ref{cl:exists root}, the existence of such an $x$ is equivalent to the existence of an $x$ with $\rv_\lambda(x-c) = \xi$ such that $P$ has a $\lambda$-collision at $x$ around certain points $b_j$ from $F$ (namely around $c$ and around the roots of the derivatives of $P$, which are in $F$ by (\ref{cond:F})). For fixed $P$ and $b_i$, the existence of such a collision is determined by $\rv_\lambda(x - b_j)$, so
it remains to eliminate the $\exists x$ from an $\Lqe(B)$-formula of the form
$\exists x\colon \psi'''((\rv_\lambda(x - b_j))_j)$ (with $\psi'''$ fqf, expressing that the collisions exist). This can then be further simplified to a formula of the form $\exists x\colon \bigwedge_j \rv_\lambda(x - b_j) = \xi_j$ (where we take $\xi_j := \rv_\lambda(a - b_j))$).
This formula now expresses that the intersection of certain balls is non-empty, a condition which can easily be seen to only depend on $\rv_\lambda(b_j - b_{j'})$ (see \cite[Proposition~4.1]{Flen} for details). Thus we are done.
\end{proof}

\begin{proof}[Proof of Theorem \ref{thm:hmin analytic}]
We need to show that for every $K' \equiv K$ and every $A \subset K'$, every $(A \cup \RV_{K',\lambda})$-definable set $X = \phi(K') \subset K'$ can be $\lambda$-prepared by a finite $A$-definable set $C$.

By Proposition \ref{prop:QE analytic}, we may assume that $\phi$ contains no field quantifiers, so that it
suffices to $\lambda$-prepare the graph of functions of the form \(\rv_\lambda(\tau(x))\), where \(\tau\) is an \(\Lqe(A)\)-term. (Here, ``preparing a graph'' is in the sense of Definition \ref{defn:uniform}.)

By the first half of Lemma \ref{lem:lin an}, such a graph can be prepared by a finite set $C$ of elements that are algebraic over
the field $F \le K'$ generated by $A$. Since such a set $C$ is contained in a finite $A$-definable set $C'$, we are done.
\end{proof}

\begin{cor}[Equi-characteristic $0$ examples]\label{cor:ex:equi}
Let $K$ be a Henselian valued field of equi-characteristic $0$ in a language $\cL$
containing the pure valued field language $\Lval = \{+,\cdot,\cO_K\}$. Then in each of the following cases, $\Th_{\cL}(K)$ is $\omega$-h-minimal.
\begin{enumerate}
 \item $\cL$ is equal to $\Lval$.
 \item\label{ex:analytic:00} The $\cL$-structure $K$ is an expansion of the $\Lval$-structure by an analytic $\cA$-structure in the sense of \cite{CLip}, for some separated Weierstrass system $\cA$.
\end{enumerate}
\end{cor}

\begin{proof}
These are just examples of Theorem~\ref{thm:hmin analytic}, namely with $\cO_K = \cO_{K,\fine}$.
\end{proof}

\begin{cor}[Mixed characteristic examples]\label{cor:ex:mixed}
Let $K$ be a Henselian valued field of mixed characteristic in a language $\cL$. Then in each of the following cases, $\Th(K)$ is $\omega$-\heqc-minimal (as in Definition~\ref{defn:mixed}).
\begin{enumerate}
 \item $\cL$ is the pure valued field language $\Lval = \{+,\cdot,\cO_K\}$.
 \item $K$ is a finite field extension of $\QQ_p$ and $\cL$ is the sub-analytic language from \cite{vdDHM} (which is a variant on the language from \cite{DvdD}).
\item The $\cL$-structure $K$ is an expansion of the $\Lval$-structure by an analytic $\cA$-structure in the sense of \cite{CLip}, for some separated Weierstrass system $\cA$.
\end{enumerate}
\end{cor}

\begin{proof}
Given $K' \equiv_{\cL} K$, let $|\cdot|_{\eqc}$ be the finest equi-characteristic $0$ coarsening
of the valuation $|\cdot|$, and let $\cL_{\eqc}$ be the extension of $\cL$ by a predicate symbol for the valuation ring $\cO_{K,\eqc}$ corresponding to $|\cdot|_{\eqc}$.
Suppose that $|\cdot|_{\eqc}$ is non-trivial.
Under those assumptions, we need to show that $\Th_{\cL_{\eqc}}(K'_{\eqc})$ is $\omega$-h-minimal, where $K'_{\eqc}$ is the field $K'$ considered as a valued field
with the valuation $|\cdot|_{\eqc}$.

(1)
If $\cL$ is the pure valued field language, we consider $\cL_{\eqc}$ as an extension of $\cL_{\val,\eqc} := \{+,\cdot,\cO_{K',\eqc}\}$ (which is also the pure valued field language) by a predicate symbol
for $\cO_{K'}$. By Theorem~\ref{thm:hmin analytic}, $\Th_{\cL_{\val,\eqc}}(K')$ is $\omega$-h-minimal.
Since the map $K' \to  \Gamma_{K'}$ factors over $\RV_\eqc$ (where $\RV_\eqc$ denotes the leading term structure with respect to $|\cdot|_\eqc$),
the $\cL_\eqc$-structure on $K'$ is an $\RV_\eqc$-expansion of the $\cL_{\val,\eqc}$-structure, so Theorem~\ref{thm:resp:h} implies that $\Th_{\cL_{\eqc}}(K')$ is also $\omega$-h-minimal.

(2)
This language $\cL$ defines an analytic structure on $K$ (by \cite[Section~4.4]{CLip} (2)) and hence also on $K'$. Hence, it suffices to prove (3).

(3)
The language $\cL_{\eqc}$ has the shape of the language called $\cL$ in the above Case~(ii) of Section \ref{sec:analyt}, so
by Theorem~\ref{thm:hmin analytic}, $\Th_{\cL_{\eqc}}(K')$ is $\omega$-h-minimal.
\end{proof}

\subsection{$\Tomin$-convex valued fields}
\label{sec:Tcon}

Fix a language $\Lomin$ containing the language of ordered rings
and fix a complete o-minimal $\Lomin$-theory $\Tomin$ containing the theory RCF of real closed fields.
Given a pair of models $K_0 \prec K$ of $\Tomin$, we can turn $K$
into a valued field by using the convex closure of $K_0$ in $K$ as the valuation ring $\cO_K$.
We suppose that $\cO_K \ne K$ and we let $\cL$ be the extension of $\Lomin$ by a predicate symbol for $\cO_K$.
In \cite{DL.Tcon1,Dri.Tcon2} van den Dries--Lewenberg obtained various results about the model theory of such valued fields $K$ as $\cL$-structures.
In particular, the theory $\cT := \Th_{\cL}(K)$ only depends on $\Tomin$ (and not on the choice of $K$ and $K_0$, provided that $\cO_K \ne K$) \cite[Corollary~3.13]{DL.Tcon1}. (Van den Dries--Lewenberg call such a ring $\cO_K$ a ``$\Tomin$-convex subring of $K$''. Accordingly, and following other subsequent literature, we call $K$ a ``$\Tomin$-convex valued field''.)

We will prove that this theory $\cT$ is
$1$-h-minimal, under the assumption that no fast-growing functions are definable in $\Tomin$. In $\RR$, ``no fast-growing functions'' means that every definable function is eventually bounded by a function of the form $x \mapsto x^n$.
In arbitrary real closed fields, the right generalization is
power-boundedness; see \cite{Mil.powBd}:

\begin{defn}[Power-bounded]
A \emph{power function} in $K$ is an $\Lomin$-definable function $g\colon K^\times \to K^\times$ which is an endomorphism of the multiplicative group $K^\times$. We call the $\Lomin$-structure $K$ (and its theory $\Tomin$) \emph{power-pounded}, if for every $\cL(K)$-definable function $f\colon K \to K$, there exists a
power function $g$ such that $|f(x)| \le g(x)$ for all sufficiently large $x \in K$.
\end{defn}

From now on, we will assume that $\Tomin$ is power-bounded.

The proof that $\Th(K)$ is $1$-h-minimal is essentially contained in the existing literature: Using the criteria given in Theorem~\ref{thm:tame2vf}, this can be deduced from \cite[Theorems~2.1 and 2.9]{Gar.powbd}. However, Theorem~2.9 is a lot deeper than what we really need, so we give a more direct proof below (mainly following the ideas from \cite{Gar.powbd}).

\begin{lem}\label{lem:tcon-0-h}
The theory of $K$ (as an $\cL$-structure) is $0$-h-minimal.
\end{lem}

\begin{proof}
We assume that $K$ is sufficiently saturated and use the criterion given by Lemma~\ref{lem:type-0-h-min}: Given a parameter set $A \subset K$ and a ball $B \subset K \setminus \acl_K(A)$, we need to verify that all elements of $B$ have the same type over $A\cup \RV$. We may assume $A = \acl_K(A)$.

Since $B \cap A = \emptyset$ and since $\Lomin$-types over $A$ correspond to cuts in $A$,
all elements of $B$ have the same $\Lomin$-type over $A$.
Thus \cite[Lemma~3.15]{Yin.tcon} applies and tells us that for any $x, x' \in B$, there exists an automorphism of $K$ fixing $A$ and $\RV$ but sending $x$ to $x'$. This shows that $x$ and $x'$ have the same type over $A \cup \RV$.
\end{proof}

The following lemma states that $\cL$-definable functions are piecewise $\Lomin$-definable. This is already stated in \cite[Lemma~2.6]{Dri.Tcon2}, but we shall use a variant from \cite{Yin.tcon}:

\begin{lem}[{\cite[Lemma~3.3]{Yin.tcon}}]\label{lem:piecewise-omin}
Let $f\colon K \to K$ be an $\cL(A)$-definable function, for some
$A \subset K \cup \RV$. Then there exists a partition of $K$ into finitely many $\cL(A)$-definable sets $X_i$ and $\Lomin(A \cap K)$-definable functions $g_i\colon K \to K$ such that $f|_{X_i} = g_i|_{X_i}$ for each $i$.
\end{lem}

It might be a bit unclear from the formulation of that lemma in \cite{Yin.tcon} whether it is intended that parameters from $\RV$ are allowed. In any case, the proof given in \cite{Yin.tcon} goes through with parameters from $\RV$.

\begin{thm}[$\Tomin$-convex examples]\label{thm:Tcon}
Let $\Tomin$ be a power-bounded o-minimal theory containing the theory of real closed fields, in a language $\Lomin$ containing the language of rings.
Let $\cT$ be the theory of $\Tomin$-convex valued fields, in the language $\cL$ which is obtained by extending $\Lomin$ by a predicate symbol for the valuation ring (as explaind at the beginning of this subsection).
Then $\cT$ is $1$-h-minimal.
\end{thm}

\begin{proof}
We use the criteria from Theorem~\ref{thm:tame2vf} to prove $1$-h-minimality, so let an $A$-definable map $f\colon K \to K$ be given, for some $A \subset K \cup \RV$.
Let $X_i \subset K$ and $g_i\colon K \to K$ be as obtained from Lemma~\ref{lem:piecewise-omin}.

Condition~(T2) of Theorem~\ref{thm:tame2vf} holds for each $g_i$ by o-minimality (namely, the set $\{d \in K \mid g_i^{-1}(d)$ is infinite$\}$ is finite), so it also holds for $f$.

Condition~(T1), too, follows for $f$ if we can prove it for each $g_i$. Indeed, take the union of the sets $C$ obtained for all the $g_i$, and further enlarge $C$ so that it $1$-prepares each $X_i$. (This is possible, since by Lemma~\ref{lem:tcon-0-h}, we already have $0$-h-minimality.) So it remains to prove Condition~(T1) for an $\Lomin(A \cap K)$-definable function $g_i$.

By o-minimality, we find a finite $(A \cap K)$-definable $C \subset K$ such that $g_i$ is continuously differentiable on $K \setminus C$. Further enlarge $C$ (using $0$-h-minimality and Corollary~\ref{cor:prep}) so that it $1$-prepares the map $K \to \Gamma_K, x \mapsto |g_i'(x)|$.

Let $B \subset K$ be a ball $1$-next to $C$. We claim that Condition~(T1) is satisfied with $\mu_B := |g_i'(x)|$ for any $x \in B$. Indeed, let $x_1, x_2 \in B$ be given, with $x_1 \ne x_2$. By the Mean Value Theorem for o-minimal fields,
there exists an $x_3$ in-between (and hence also in $B$) such that $g_i(x_1) - g_i(x_2) = g_i'(x_3)\cdot (x_1 - x_2)$. Taking valuations on both sides implies $|g_i(x_1) - g_i(x_2)| = \mu_B\cdot |x_1 - x_2|$, as desired.
\end{proof}

\begin{remark}
The assumption that $\Tomin$ is power-bounded is necessary to obtain $1$-h-minimality of $\cT$. Indeed, in the presence of an exponential map, we can define $K \to \RV, x \mapsto \rv(e^x)$, whose fibers are exactly the translates of the maximal ideal $B_{<1}(0)$ and which hence cannot be $1$-prepared in the sense of Corollary~\ref{cor:prep}.
\end{remark}

\begin{remark}
We were not able to prove that power-bounded $\Tomin$-convex valued fields are $\omega$-h-minimal. This is one of the main reasons why we only assume $1$-h-minimality in most of the paper.
\end{remark}

Using methods from non-standard analysis, results in a $\Tomin$-convex valued field $K$ can often be translated into results about $K$ as an $\Lomin$-structure.
We finish this subsection by stating what our Taylor approximation result (Theorem~\ref{thm:high-ord}) becomes under such a translation, namely a version of Taylor approximation which has some uniformity even when one approaches a bad point.

\begin{cor}\label{cor:arch}
Let $\Tomin$ be a power-bounded o-minimal theory containing the theory of real closed fields, in a language $\Lomin$ containing the language of rings.
Let $K \models \Tomin$ be a model, let $f\colon K \to K$ be an $\Lomin$-definable function and let $r \in\ \NN$ be given.
Then there exists a finite $\Lomin$-definable set $C \subset K$ and a constant $c \in K_{>0}$ such that for every pair $x_0, x \in K$ satisfying
\begin{equation}\label{eq:arch.ass}
c\cdot |x - x_0| < \min_{a \in C}|x_0 -  a|,
\end{equation}
we have
\begin{equation}\label{eq:arch.imp}
|f(x) -  T^{\le r}_{f,x_0}(x) | \leq c\cdot |f^{(r+1)}(x_0)\cdot (x-x_0)^{r+1}|
\end{equation}
(where $|\cdot|$ denotes the usual absolute value and $T^{\le r}_{f,x_0}$ is the Taylor polynomial of $f$ around $x_0$ of degree $r$; see Definition~\ref{defn:taylor}).

More generally, if $(f_q)_{q \in K^m}$ is an $\Lomin$-definable family of functions $K \to K$, then
we obtain an $\Lomin$-definable family of sets $(C_q)_{q \in K^m}$ and a constant $c \in K_{>0}$ which is independent of $q$ such
that the above holds for every $q \in K^m$.
\end{cor}

\begin{remark}
Note that this result would be false without the assumption on power-boundedness. Indeed, one can check that it fails near $0$ for the function $x \mapsto e^{1/x}$.
On the other hand, it should be rather easy to obtain for sub-analytic functions, so this is another instance (along with the Jacobian Property)
of a generalization of a result from an analytic setting to an axiomatic one.
\end{remark}

\private{
Set $r = 0$ and $x = x_0 + d$, with $dc < x_0$. We consider $|e^{1/x} - e^{1/x_0}| < |c \cdot (e^{1/x_0})'\cdot d| = |c \cdot e^{1/x_0} \cdot x_0^{-2}\cdot d|$.
After dividing by $e^{1/x_0}$, we obtain $|e^{1/x-1/x_0}| = |e^{-d/(x\cdot x_0)} - 1| < cdx_0^{-2}$. In the LHS, we approximate $x$ by $x_0$, and then we approximate
the exponential by the degree 3 Taylor approx. Then one sees that the corollary fails.
}

In the proof, we use the following lemma:

\begin{lem}\label{lem:L2Lomin}
For any $A \subset K$, every finite $\cL(A)$-definable set $C \subset K$ is already $\Lomin(A)$-definable.
\end{lem}

\begin{proof}
Using the order, we reduce to the case where $C = \{a\}$ is a singleton.
$\cL(A)$-definability means that $a = f(0)$, for some $\cL(A)$-definable function $f\colon K \to K$.
By Lemma~\ref{lem:piecewise-omin}, $f(0) = g(0)$ for some $\Lomin(A)$-definable $g\colon K \to K$; this implies that $C$ is $\Lomin(A)$-definable.
\end{proof}

\begin{proof}[Proof of Corollary~\ref{cor:arch}]
Fix a $|K|^+$-saturated elementary extension $K' \succ K$ and let $\cO_{K'}$ be the convex closure of $K$ in $K'$. By Theorem~\ref{thm:Tcon}, $K'$ is $1$-h-minimal as an $\cL$-structure, for $\cL$ as in the theorem.
(Note that the saturation assumption implies $\cO_{K'} \ne K'$.)
In the following, we denote the valuation on $K'$ by $|\cdot|_v$, to distinguish it from the absolute value $|\cdot|$. We suppose that a family of functions $(f_q)_{q \in K^m}$ is given as in the corollary, and by abuse of notation, we also write $(f_q)_{q \in (K')^m}$ for the corresponding family in $K'$ (defined by the same formula).

Suppose that the corollary fails, i.e., that for every $c \in K_{>0}$ and for every formula $\psi$ that could potentially define the family $(C_q)_{q \in K^m}$, there exists a $q \in K^m$ and a pair $x_0, x$ for which the implication (\ref{eq:arch.ass}) $\Rightarrow$ (\ref{eq:arch.imp}) fails.
By our saturation assumption, there exist $q \in (K')^m$ and $x_0, x \in K'$ such that the implication fails for every $c \in K_{>0}$ and every $\psi$.
Using that $\cO_{K'}$ is the convex closure of $K$, this failure for every $c \in K_{>0}$ is equivalent to the conjunction
\begin{equation}\label{eq:v.ass}
|x - x_0|_v < \min_{a \in C_q}|x_0 -  a|_v
\end{equation}
and
\begin{equation}\label{eq:v.imp}
|f(x) -  T^{\le r}_{f_q,x_0}(x) |_v >  |f_q^{(r+1)}(x_0)\cdot (x-x_0)^{r+1}|_v.
\end{equation}
So we obtained: For every $\Lomin(q)$-definable set $C_q$, there exist $x_0, x \in K'$ in the same ball $1$-next to $C_q$ (by (\ref{eq:v.ass})) such that (\ref{eq:v.imp}) holds.
This contradicts Theorem~\ref{thm:high-ord}: \emph{A priori}, that theorem only provides a finite $\cL(q)$-definable set $C$, but $\Lomin(q)$-definability of that set then follows using Lemma~\ref{lem:L2Lomin}.
\end{proof}

One can expect that similar results in higher dimension can be obtained,
and that they may lead to finer versions of the preparation results in power-bounded real closed fields from \cite{DS.ominPrep,NguyenValette}, which are used to deduce the existence of Mostowski's Lipschitz stratifications.

\subsection{Comparison to V-minimality}\label{sec:comparison}

In \cite{HK}, Hrushovski--Kazhdan introduced the notion of $V$-minimal theories, which at first sight has the same goal as Hensel minimality, namely: to provide a powerful axiomatic framework for geometry in valued fields. However,
the relation between $V$-minimality and Hensel minimality is similar to the relation between strong minimality and o-minimality:
By working in a strongly minimal theory of (algebraically closed) fields, one obtains many useful results about geometry in real closed fields, but one cannot treat genuinely o-minimal languages like $\RR_{\mathrm{exp}}$. In a similar way, working in a $V$-minimal theory of (algebraically closed) valued fields does provide many useful insights about Henselian valued fields (as explained in \cite[Section~12]{HK}), but there are examples of Hensel minimal theories that cannot be treated in this way.

Concretely, since a V-minimal theory has not more structure on $\RV$ than a pure valued field, every definable function $K \to K$ ultimately grows like $x \mapsto x^r$ for some rational number $r$,
and this remains true if we expand the language by predicates on $\RV$ (following \cite[Section~12]{HK}). In contrast, Section~\ref{sec:Tcon} provides examples of $1$-h-minimal structures without this property, for example the $T$-convex structure obtained from the (power-bounded o-minimal) expansion of $\RR$ by one function $x \mapsto x^r$ for every real number $r$.

On the other hand, if we restrict to the context for which V-minimality has been designed, then it agrees with Hensel minimality. Moreover, in this case, $0$-h-minimality
and $1$-h-minimality agree. Let us recall the definition of V-minimality.

\begin{defn}[V-minimality; {\cite[Section~3.4]{HK}}]\label{defn:Vmin}
Fix a language $\cL \supset \Lval$ and a complete theory $\cT$ containing the theory $\operatorname{ACVF}_{0,0}$ of algebraically closed fields of equi-characteristic $0$.
The theory $\cT$ is called \emph{$V$-minimal} if for every model $K \models \cT$, we have the following:
\begin{enumerate}\setcounter{enumi}{-1}
 \item Every definable (with parameters) subset of $K$ is a finite boolean combination of points, open balls, and closed balls.
 \item Every $\cL(K)$-definable subset of $\RV^n$ is already
   $\Lval(K)$-definable (where $\Lval$ is the pure valued field language).
 \item Every definable (with parameters) family of nested closed balls in $K$ has non-empty intersection.
 \item For every $A \subset K$, if $X \subset K$ is an $A$-definable set which is the union of finitely many disjoint closed balls $B_i$, then $\acl_K(A) \cap B_i \ne \emptyset$ for every $i$.
\end{enumerate}
\end{defn}

\private{In the last item, HK write that $X$ should be an ``almost $A$-definable closed ball''; this is defined as: there exists an $A$-definable equivalence relation with finitely many classes such that $X$ is a union of equivalence classes....}

\begin{prop}[V-minimality]\label{prop:vmin}
Suppose that $\cT$ is a complete theory containing $\operatorname{ACVF}_{0,0}$, in a language $\cL \supset \Lval$, and suppose moreover that every $\cL(K)$-definable subset of $\RV^n$ is already $\Lval(K)$-definable.
Then the following are equivalent:
\begin{enumerate}[(i)]
 \item $\cT$ is V-minimal.
 \item $\cT$ is $0$-h-minimal.
 \item $\cT$ is $1$-h-minimal.
\end{enumerate}
\end{prop}

\begin{proof}[Proof of Proposition~\ref{prop:vmin}]
(iii) $\Rightarrow$ (ii) is trivial.

(i) $\Rightarrow$ (iii): We use the criteria from Theorem~\ref{thm:tame2vf}, so let $f\colon K \to K$ be $A$-definable,
for some $A \subset K \cup \RV$.

First note that by the Remark just above \cite[Lemma~3.30]{HK}, adding parameters from $K \cup \RV$ to the language preserves V-minimality, so using compactness,
the results from \cite{HK} hold uniformly in families parametrized by $K$ or $\RV$.

By dimension theory (e.g.\ \cite[Lemma~3.55]{HK}),
$f$ has only finitely many infinite fibers, i.e., Condition~(T2) from Theorem~\ref{thm:tame2vf} holds.
By applying \cite[Corollary~4.3]{HK} to all fibers $f^{-1}(b)$ of $f$ (where $b$ runs over $K$),
we find an $A$-definable map $\rho\colon K \to \RV^k$ (for some $k \ge 0$) such that for each $\xi\in \RV^k$, the restriction $f|_{\rho^{-1}(\xi)}$ is either constant or injective.
Apply \cite[Corollary~5.9]{HK} to each injective restriction $f|_{\rho^{-1}(\xi)}$ and refine the map $\rho$ accordingly, i.e.,
such that afterwards, $f$ is ``nice'' on each open ball contained in one fiber of $\rho$ in the sense of \cite[Definition~5.8]{HK}.
Finally, by applying \cite[Corollary~4.3]{HK} to the graph of $\rho$, we find a finite $A$-definable set $C$ $1$-preparing $\rho$ (namely, the image of the map $c$ provided by the corollary).
Then $f$ is nice on each ball $1$-next to $C$, and this implies Condition~(T1) from Theorem~\ref{thm:tame2vf}.

(ii) $\Rightarrow$ (i): We prove the conditions from Definition~\ref{defn:Vmin}:

(0) Let $X \subset K$ be definable, and let $C \subset K$ be a finite set preparing $X$. Then $X$ can be written as a union of the form
$\bigcup_{c \in C} X_c$, where $X_c = \{c + x \mid \rv(x) \in Z_c\}$
for suitable definable sets $Z_c \subset \RV$. Using the assumption that definable subsets of $\RV$ are already definable in the language $\Lval$,
we obtain that each $X_c$ is a finite boolean combination of points, open balls, and closed balls. This then also follows for $X$.

(1) holds by assumption.

(2) is a special case of Lemma~\ref{lem:intersection}.

(3) By 0-h-minimality, there exists a finite $A$-definable set $C$ preparing $X$. This set $C$ cannot be disjoint from any $B_i$, since for any $c \in C \setminus B_i$, the ball $1$-next to $c$ containing $B_i$ is strictly bigger than $B_i$.
\end{proof}

\subsection{Some open questions}

We finish the paper with a few questions; more questions are stated in the sequel \cite{CHRV}.
Probably the most
immediate question in this context is:

\begin{question}
Does $0$-h-minimality imply $1$-h-minimality, and, does $1$-h-minimality imply $\omega$-h-minimality? More generally: For which $\ell < \ell'$ does $\ell$-h-minimality imply $\ell'$-h-minimality?
\end{question}

If $1$-h-minimality is not equivalent to $\omega$-h-minimality, there are still the following questions, motivated by the results from Sections \ref{sec:Tcon} and \ref{sec:comparison}.

\begin{question}
Are the $\Tomin$-convex structures from Section~\ref{sec:Tcon} (with $\Tomin$ power-bounded) $\omega$-h-minimal?
Does V-minimality imply $\omega$-h-minimality?
\end{question}

\bibliographystyle{amsplain}
\bibliography{anbib}
\end{document}